\documentclass[11pt,a4paper]{amsart}

\usepackage{graphicx,amssymb}
\usepackage[top=30mm,bottom=25mm,left=20mm,right=20mm]{geometry} 
\usepackage[frame]{xy}
\xyoption{arc}
\xyoption{graph}
\xyoption{matrix}

\theoremstyle{plain}
\newtheorem{theorem}{Theorem}[section]
\newtheorem{lemma}[theorem]{Lemma}
\newtheorem{corollary}[theorem]{Corollary}
\newtheorem{proposition}[theorem]{Proposition}

\theoremstyle{definition}
\newtheorem{definition}[theorem]{Definition}
\newtheorem{example}[theorem]{Example}
\newtheorem{remark}[theorem]{Remark}

\makeatletter
\newif\if@my@Hidebox \@my@Hideboxtrue
\newbox\@my@Hidebox
\def\myHidebox{\setbox\@my@Hidebox\vbox\bgroup}
\def\endmyHidebox{\egroup
   \if@my@Hidebox
      \unvbox\@my@Hidebox\par
   \else\par
   \fi
}
\def\showmyHidebox{\@my@Hideboxtrue}
\def\hidemyHidebox{\@my@Hideboxfalse}
\def\showmyproof{\@my@Hideboxtrue}
\def\hidemyproof{\@my@Hideboxfalse}

\makeatother

\hidemyHidebox  
\showmyHidebox  

\def\EKLnewpage{}

\let\le\leqslant
\let\ge\geqslant
\def\Z{\mathbb Z}

\let\rla\rightleftharpoons

\def\St{\operatorname{St}}
\def\Lk{\operatorname{Lk}}

\def\supp{\operatorname{supp}}
\def\Lk{{\operatorname{Lk}}}

\def\diam{\operatorname{diam}}

\def\minTL{\mathcal L}
\def\Mod{\operatorname{Mod}}

\begin{document}

\title[Acylindricity of the action of RAAGs on extension graphs]
{Acylindricity of the action of right-angled \\
Artin groups on extension graphs}

\author{Eon-Kyung Lee and Sang-Jin Lee}
\address{Department of Mathematics and Statistics, Sejong University, Seoul, Korea}
\email{eonkyung@sejong.ac.kr}

\address{Department of Mathematics, Konkuk University, Seoul, Korea}
\email{sangjin@konkuk.ac.kr}
\date{\today}

\begin{abstract}
The action of a right-angled Artin group on its extension graph is known to be acylindrical
because the cardinality of the so-called $r$-quasi-stabilizer
of a pair of distant points is bounded above by a function of $r$.
The known upper bound of the cardinality is an exponential function of $r$.
In this paper we show that
the $r$-quasi-stabilizer is a subset of a cyclic group and its cardinality is bounded above by a linear function of $r$.
This is done by exploring lattice theoretic properties of group elements,
studying prefixes of powers
and extending the uniqueness of quasi-roots
from word length to star length.
We also improve the known lower bound for the minimal asymptotic translation length
of a right angled Artin group on its extension graph.

\medskip\noindent
{\em Keywords\/}:
Right-angled Artin groups, extension graphs,
acylindrical actions, translation lengths.\\
{\em 2010 Mathematics Subject Classification\/}: Primary 20F36; Secondary 20F65
\end{abstract}

\maketitle


\section{Introduction}
\label{sec:intro}

Throughout the paper $\Gamma$ denotes a finite simplicial graph,
not necessarily connected,
with vertex set $V(\Gamma)$ and edge set $E(\Gamma)$.
The \emph{right-angled Artin group} $A(\Gamma)$
with the underlying graph $\Gamma$ is the group generated by
$V(\Gamma)$ such that the defining relations are
the commutativity between adjacent vertices,
hence $A(\Gamma)$ has the group presentation
$$
A(\Gamma)=\langle\, v\in V(\Gamma)\mid v_iv_j=v_jv_i\ \
\mbox{for each $\{v_i,v_j\}\in E(\Gamma)$}\,\rangle.
$$
Right-angled Artin groups are important groups in geometric group theory,
which played a key role in Agol's proof of the
virtual Haken conjecture~\cite{Ago13,KM12,Wis21}.

The \emph{extension graph} $\Gamma^e$ is the graph such that
the vertex set $V(\Gamma^e)$ is the set of all elements of $A(\Gamma)$
that are conjugate to a vertex of $\Gamma$,
and two vertices $v_1^{g_1}$ and $v_2^{g_2}$
are adjacent in $\Gamma^e$ if and only if
they commute when considered as elements of $A(\Gamma)$.
(Here, $v^g$ denotes the conjugate $g^{-1}vg$.)
Therefore
\begin{align*}
V(\Gamma^e) & =\{v^g: v\in V(\Gamma),\ g\in A(\Gamma)\},\\
E(\Gamma^e) & =\{\, \{v_1^{g_1}, v_2^{g_2}\}\,:\,
v_1^{g_1}v_2^{g_2}=v_2^{g_2}v_1^{g_1}\ \mbox{in $A(\Gamma)$}\,\}.
\end{align*}
Extension graphs are usually infinite and locally infinite.
They are very useful
in the study of right-angled Artin groups
such as the embeddability problem
between right-angled Artin groups~\cite{KK13, KK14a, LL16, LL18, Kat18}
and the purely loxodromic subgroups
which are analogous to convex cocompact subgroups
of the mapping class groups of surfaces~\cite{KMT17}.
It is known that $\Gamma^e$ is a quasi-tree,
hence a $\delta$-hyperbolic graph~\cite{KK13}.

\begin{definition}[acylindrical action]
When a group $G$ acts on a path-metric space $(X, d)$
isometrically from the right,
the action is called \emph{acylindrical} if for any $r>0$,
there exist $R,N>0$ such that
whenever $x$ and $y$ are two points of $X$ with $d(x,y)\ge R$,
the cardinality of the set
$$
\xi(x,y;r)
=\{g\in G:
d(xg,x)\le r\mbox{ and }
d(yg,y)\le r\}$$
is at most $N$.
The set $\xi(x,y;r)$ is called the \emph{$r$-quasi-stabilizer}
of the pair of points $(x,y)$.
We sometimes use the notation $\xi_{(X,d)}(x,y;r)$
for the set $\xi(x,y;r)$.
Notice that $R$ and $N$ are functions of $r$.
When we need to specify the acylindricity constants
$R$ and $N$, we say that
the action is \emph{$(R,N)$-acylindrical}.
\end{definition}

There have been many works on properties and examples
of groups with an acylinrical action
on a geodesic hyperbolic metric space.
For example, see~\cite{Bow08, Osi16, DGO17}.

Let $d$ denote the graph metric of $\Gamma^e$.
The right-angled Artin group $A(\Gamma)$ acts
on $(\Gamma^e,d)$ isometrically from the right by conjugation, i.e.~
the image of the vertex $v^h$ under the action of
$g\in A(\Gamma)$ is $v^{hg}$.
The action of $A(\Gamma)$ on $\Gamma^e$ behaves much like
the action of the mapping class group $\mbox{Mod}(S)$ of a hyperbolic surface
$S$ on the curve graph $\mathcal C(S)$.
One of the fundamental properties is that
the action of $A(\Gamma)$ on $\Gamma^e$ is acylindrical,
which is shown by Kim and Koberda~\cite{KK14}.

\begin{theorem}[{\cite[Theorem 30]{KK14}}]
The action of $A(\Gamma)$ on $\Gamma^e$ is acylindrical.
\end{theorem}

More precisely, it is shown that the action is $(R,N)$-acylindrical with
\begin{align*}
R&=R(r) =D(2r+4D+7),\\
N&=N(r) =\left(V 2^{2V}\right)^{r+2D+1},
\end{align*}
where $D=\diam(\Gamma)$ is the diameter of $\Gamma$
and $V=|V(\Gamma)|$ is the cardinality of $V(\Gamma)$.
Notice that $N(r)$ is an exponential function of $r$.

\medskip
For a graph $\Gamma$, let $\bar\Gamma$ denote
the complement graph of $\Gamma$, i.e.~the graph
on the same vertices as $\Gamma$ such that
two distinct vertices are adjacent in $\bar\Gamma$
if and only if they are not adjacent in $\Gamma$.

For the reader's convenience, we give some remarks
on the cases where $|V(\Gamma)|$ is small
and where $\Gamma$ or $\bar\Gamma$ is disconnected.

The following are known for the extension graph
$\Gamma^e$~\cite[Lemma 3.5]{KK13}:
if $\Gamma$ is disconnected,
then $\Gamma^e$ has countably infinite number of path-components;
if $\bar\Gamma$ is disconnected, i.e.\ $\Gamma$ is a join,
then $\Gamma^e$ is also a join, hence $\diam(\Gamma^e)\le 2$;
if $|V(\Gamma)|=1$, then $|V(\Gamma^e)|=1$.
If $|V(\Gamma)|\in\{2,3\}$,
then either $\Gamma$ or $\bar\Gamma$ is disconnected.
In fact, $\Gamma^e$ is a connected graph
with infinite diameter if and only if $|V(\Gamma)|\ge 4$
and both $\Gamma$ and $\bar \Gamma$ are connected.
Therefore, when we consider the action of $A(\Gamma)$ on $\Gamma^e$,
it is natural to require that $|V(\Gamma)|\ge 4$ and
both $\Gamma$ and $\bar \Gamma$ are connected.

In the study of extension graphs, we use the star length
metric $d_*$ on $A(\Gamma)$. (See \S\ref{sec:star} for the definition
of star length.) The metric space $(A(\Gamma), d_*)$ is quasi-isometric
to the extension graph $(\Gamma^e, d)$.
If $|V(\Gamma)|=1$
or if $\bar\Gamma$ is disconnected, then $(A(\Gamma), d_*)$ has diameter
at most 2, which is not interesting.
Therefore, when we consider the action of $A(\Gamma)$ on $(A(\Gamma),d_*)$,
it is natural to require that $|V(\Gamma)|\ge 2$ and
$\bar \Gamma$ is connected (see Remark~\ref{rmk:st-ass}).

From the above discussions, the following settings are natural.
\begin{enumerate}
\item When we consider the action of $A(\Gamma)$ on $(\Gamma^e,d)$,
we will assume that $|V(\Gamma)|\ge 4$ and
both $\Gamma$ and $\bar\Gamma$ are connected.
\item When we consider the action of $A(\Gamma)$ on $(A(\Gamma), d_*)$,
we will assume that $|V(\Gamma)|\ge 2$ and $\bar\Gamma$ is connected.
\end{enumerate}

\medskip
The following is the main result of this paper,
which shows that we can take $N(r)$ as a linear function of $r$
and furthermore the quasi-stabilizer $\xi(x,y;r)$ is a subset of a cyclic group.

\medskip\noindent
\textbf{Theorem A} (Theorem~\ref{thm:ex})\ \
\emph{
Let $\Gamma$ be a finite simplicial graph
such that $|V(\Gamma)|\ge 4$
and both $\Gamma$ and $\bar\Gamma$ are connected.
Then the action of\/ $A(\Gamma)$ on $\Gamma^e$ is
$(R, N)$-acylindrical with
\begin{align*}
&R=R(r) =D(2V+7)(r+1)+10D,\\
&N=N(r) =2(V-2)r-1,
\end{align*}
where $D=\diam(\Gamma)$ and $V=|V(\Gamma)|$.
Moreover, for any $x,y\in V(\Gamma^e)$ with
$d(x,y)\ge R$, if\/ $\xi(x,y;r)\ne \{1\}$, then
there exists a loxodromic element $g\in A(\Gamma)$
such that
\begin{enumerate}
\item
$\xi(x,y;r)\subset\{1,g^{\pm 1}, g^{\pm 2},\ldots, g^{\pm k}\}$
for some $1\le k\le (V-2)r-1$;
\item
the Hausdorff distance between the $\langle g\rangle$-orbit of $x$ and
that of $y$ is at most $D(2r+7)$.
\end{enumerate}
}

\bigskip

The following is an easy example to come up with for $g\in\xi(x,y;r)$.
Let $g$ be a loxodromic element with a quasi-axis
$L=z^{\langle g\rangle}=\{z^{g^m}:m\in\Z\}$
for some $z\in V(\Gamma^e)$
such that $d(z^g,z)$ is sufficiently small.
If both $x$ and $y$ are close enough to $L$,
then $d(x^g,x)$ and $d(y^g,y)$ are
also small so that $g\in\xi(x,y;r)\setminus\{1\}$.
In this case, the Hausdorff distance between
the $\langle g\rangle$-orbits
$x^{\langle g\rangle}$ and
$y^{\langle g\rangle}$ is small.
Theorem~\ref{thm:ex} says that, in some sense,
this is the only case where $g\in \xi(x,y;r)\setminus\{1\}$ happens:
$g$ is loxodromic and
the Hausdorff distance between $x^{\langle g\rangle}$ and
$y^{\langle g\rangle}$ is small.
Moreover, by Theorem~\ref{thm:ex}(i), the set $\xi(x,y;r)\setminus\{1\}$
is purely loxodromic, that is, there is no elliptic element
that $r$-quasi-stabilizes a pair of sufficiently distant points.

\medskip
In order to prove Theorem A, we develop several tools
such as lattice theoretic properties of group elements,
decomposition of conjugating elements,
properties of prefixes of powers, and then
extend the uniqueness of quasi-roots in~\cite{LL22}
from word length to star length.
Using these tools, we also obtain a new
lower bound for the minimal asymptotic translation length
of the action of $A(\Gamma)$ on $\Gamma^e$.

\begin{definition}
When a group $G$ acts on a connected metric space $(X, d)$
by isometries from right,
the \emph{asymptotic translation length} of an element $g\in G$
is defined by
\begin{equation}\label{eq:tl}
\tau(g) =\tau_{(X,d)}(g)= \lim_{n\to\infty} \frac{d(xg^n, x)}{n},
\end{equation}
where $x\in X$.
This limit always exists,
is independent of the choice of $x\in X$,
and satisfies $\tau(g^n) = |n|\tau(g)$
and $\tau(h^{-1}gh) = \tau(g)$ for all $g, h\in G$ and $n\in\Z$.
If $\tau(g)>0$, $g$ is called \emph{loxodromic}.
If $\{d(xg^n, x)\}_{n=1}^\infty$ is bounded,
$g$ is called \emph{elliptic}.
If $\tau(g)=0$ and $\{d(xg^n, x)\}_{n=1}^\infty$ is unbounded,
$g$ is called \emph{parabolic}.
For a subgroup $H$ of $G$,
the \emph{minimal asymptotic translation length} of $H$
for the action on $(X,d)$ is defined by
\begin{equation}\label{eq:tl}
\minTL_{(X,d)}(H) = \min\{ \tau_{(X,d)}(h):
h\in H,~\tau_{(X,d)}(h)>0\}.
\end{equation}
\end{definition}

There have been many works on minimal asymptotic translation lengths
of the action of mapping class groups on curve graphs.
Let $S_g$ denote a closed orientable surface of genus $g$.
For the action of the mapping class group $\Mod(S_g)$
on the curve graph $\mathcal C(S_g)$,
Gadre and Tsai~\cite{GT11} proved that
$$
\minTL_{\mathcal C(S_g)}(\Mod(S_g)) \asymp \frac{1}{g^2},
$$
where $f(g)\asymp h(g)$ denotes that
there exist positive constants $A$ and $B$ such that
$A f(g)\le h(g)\le B f(g)$.
The braid group $\mbox{B}_n$ can be regarded
as the mapping class group
of the $n$-punctured disk $D_n$ fixing boundary pointwise.
The pure braid group $\mbox{PB}_n$ is the subgroup
of $\mbox{B}_n$ consisting of mapping classes that fix each puncture.
Kin and Shin~\cite{KS19} and Baik and Shin~\cite{BS20} showed that
$$
\minTL_{\mathcal C(D_n)}(\mbox{B}_n) \asymp \frac{1}{n^2},
\qquad
\minTL_{\mathcal C(D_n)}(\mbox{PB}_n) \asymp \frac{1}{n}.
$$
For the action of $A(\Gamma)$ on $\Gamma^e$,
it follows from a result of Kim and Koberda~\cite{KK14} that
$$\minTL_{(\Gamma^e,\,d)}(A(\Gamma))
\ge \frac1{2|V(\Gamma)|^2}.$$
Baik, Seo and Shin~\cite{BSS21} proved that
all loxodromic elements of $A(\Gamma)$ on $\Gamma^e$ have rational asymptotic translation
lengths with a common denominator.

In this paper, we show the following, where
the denominator of the lower bound is
improved from a quadratic function to a linear function of $|V(\Gamma)|$.

\medskip\noindent
\textbf{Theorem B} (Theorem~\ref{thm:mTL})\ \
\emph{
Let $\Gamma$ be a finite simplicial graph
such that $|V(\Gamma)|\ge 4$ and both $\Gamma$ and $\bar\Gamma$ are connected.
Then}
$$
\minTL_{(\Gamma^e,\,d)}(A(\Gamma))
\ge
\frac1{|V(\Gamma)|-2}
\,.
$$

\bigskip
In the remaining of this section,
we explain briefly our ideas 
and the structure of this paper.

\subsection{Idea for the acylindricity}
Let us first explain our idea for the acylindricity.
For $g\in A(\Gamma)$, let $\|g\|$ denote the word length of $g$
with respect to the generating set $V(\Gamma)^{\pm 1}$,
and let $d_\ell$ denote temporarily the word length metric defined by
$d_\ell(g,h)=\|gh^{-1}\|$ for $g,h\in A(\Gamma)$.
The right multiplication induces an isometric action of
$A(\Gamma)$ on $(A(\Gamma), d_\ell)$.
Since $\xi(x,y;r)=x^{-1}\xi(1,yx^{-1};r)x$ for any $x,y\in A(\Gamma)$,
it suffices to consider $r$-quasi-stabilizers
of the form $\xi(1,w;r)$ for the acylindricity.

Suppose that we are given $R>0$ large, $r>0$ small
and $w\in A(\Gamma)$ with $\|w\|=d_\ell(w,1)\ge R$.
Let $g\in\xi(1,w;r)\setminus\{1\}$.
Since $\|g\|=d_\ell(g,1)\le r$ and $\|wgw^{-1}\|=d_\ell(wg,w)\le r$, we have
$$
\|w\|\ge R,\quad
\|g\|\le r,\quad
\|wgw^{-1}\|\le r.
\leqno(*)
$$
In other words, $\|w\|$ is large
whereas $\|g\|$ and $\|wgw^{-1}\|$ are small.
This happens typically when
$$
w=ag^n,\qquad n\in\Z,\quad a\in A(\Gamma)
\leqno(**)
$$
with $\|a\|$ small and $|n|$ large.
In this case, $d_\ell(w,g^n)=d_\ell(ag^n,g^n)=\|a\|$ is small,
hence we can say that $w$ is ``close to a power of $g$''.

Even though it is clearly over-optimistic and false,
one may hope that the following hold:
given a triple $(R,r,w)$ as above (i.e.\ $R>0$ is large, $r>0$ is small
and $w\in A(\Gamma)$ with $\|w\|\ge R$),
\begin{enumerate}
\item if $(*)$ holds, then $(**)$ holds
for some $n\in \Z$  and $a\in A(\Gamma)$ with $\|a\|$ small;
\item only a small number of triples $(a,g,n)$
with $\|a\|$ small and $\|g\|\le r$
satisfy $(**)$.
\end{enumerate}
Of course, the above statements are not true
at least as they are written.
Moreover, the metric spaces $(A(\Gamma),d_\ell)$ and $(\Gamma^e,d)$ are not
quasi-isometric, hence the above statements do not imply
the acylindricity of $(\Gamma^e,d)$.
However, we will see that this approach in fact works
in the study of the acylindricity of the action of $A(\Gamma)$ on $(\Gamma^e,d)$
if we replace the word length metric with the star length metric.

\subsection{Lattice structure}
In~\S\ref{sec:lattice}, we collect basic combinatorial
group theoretic properties of right-angled Artin groups.
Those properties are stated using lattice theoretic notations.

The motivation comes from Garside groups which are a lattice theoretic generalization
of braid groups and finite type Artin groups.
For Garside groups, there are elegant tools especially for the word and conjugacy problems
and the asymptotic translation
length~\cite{Gar69, BS72, Del72, BKL98, DP99, DDGKM15, LL07,Lee07}.
Right-angled Artin groups are not Garside groups, except free abelian groups,
hence we cannot apply Garside theory to right angled-Artin groups.
However, some ideas from Garside theory are very useful in our approach.

For $g\in A(\Gamma)$, the \emph{support} of $g$, denoted $\supp(g)$,
is the set of generators that appear in a shortest word on $V(\Gamma)^{\pm 1}$
representing $g$.

For $g_1,g_2\in A(\Gamma)$, we say that $g_1$ and $g_2$ \emph{disjointly commute},
denoted $g_1\rla g_2$, if $\supp(g_1)\cap\supp(g_2)=\emptyset$ and
each $v_1\in \supp(g_1)$ commutes with each $v_2\in\supp(g_2)$.

Let $g=g_1\cdots g_k$ for some $g,g_1,\ldots,g_k\in A(\Gamma)$.
We say that the decomposition is \emph{geodesic} if $\|g\|=\|g_1\|+\cdots+\|g_k\|$.
If $g=g_1g_2$ is geodesic,
we say that $g_1$ is a \emph{prefix} of $g$,
denoted $g_1\le_L g$, and that $g$ is a \emph{right multiple} of $g_1$.

The relation $\le_L$ is a partial order on $A(\Gamma)$,
hence the notions of gcd $g_1\wedge_L g_2$ and lcm $g_1\vee_L g_2$ make sense.
Theorem~\ref{thm:lcm} shows that
for $g_1,g_2\in A(\Gamma)$, the gcd $g_1\wedge_L g_2$ always exists
and the lcm $g_1\vee_L g_2$ exists if and only if
$g_1$ and $g_2$ have a common right multiple.
Moreover, in this case, there exist $g_1', g_2'\in A(\Gamma)$
such that $g_i=(g_1\wedge_L g_2)g_i'$ for $i=1,2$,
$g_1'\rla g_2'$ and $g_1\vee_L g_2=(g_1\wedge_L g_2) g_1'g_2'$.

\subsection{Cyclic conjugations}
In \S\ref{sec:conj}, we study conjugations $g^u=u^{-1}gu$.
The decomposition $u^{-1}gu$ is not geodesic in general,
i.e.~$\|u^{-1}gu\|\ne\|u^{-1}\|+\|g\|+\|u\|$.

Let $g$ be cyclically reduced, i.e.\
the word length $\|g\|$ is minimal in its conjugacy class.
If $u\le_L g$, then $g=ug_1$ is geodesic for some $g_1\in A(\Gamma)$
and $g^u=u^{-1}(ug_1)u=g_1u$.
In other words, the conjugation of $g$ by $u$
moves the prefix $u$ to the right.
An iteration of this type of conjugations
is called a \emph{left cyclic conjugation}.
The \emph{right cyclic conjugation} is defined similarly.
The \emph{cyclic conjugation} is an iteration of
left and right cyclic conjugations.

Proposition~\ref{prop:cc1} shows that
for a cyclically reduced element $g$,
the conjugation $g^u$ is a left cyclic conjugation of $g$
if and only if $u\le_L g^n$ for some $n\ge 1$.

Theorem~\ref{thm:cnj1} shows that given $g,u\in A(\Gamma)$ with $g$ cyclically reduced,
there exists a unique geodesic decomposition $u=u_1u_2u_3$ such that
$u_1$ disjointly commutes with $g$;
$g^{u_2}$ is a cyclic conjugation;
$g^u=u_3^{-1} g^{u_2} u_3$ is geodesic,
i.e.~$\|u_3^{-1} g^{u_2} u_3\|=\|u_3^{-1}\|+\| g^{u_2}\|+\| u_3\|$.
Furthermore, there is a geodesic decomposition $u_2=u_2'u_2''$
such that $g^{u_2'}$ (resp.\ $g^{u_2''}$) is a left (resp.\ right)
cyclic conjugation and $u_2'\rla u_2''$.

\subsection{Star length}
An element $g\in A(\Gamma)$ is called a \emph{star-word}
if $\supp(g)$ is contained in the star of some vertex.
The \emph{star length}, denoted $\|g\|_*$, of $g$
is the minimum $\ell$ such that
$g$ can be written as a product of $\ell$ star-words.
Let $d_*$ denote the metric on $A(\Gamma)$
induced by the star length: $d_*(g_1,g_2)=\|g_1g_2^{-1}\|_*$.

The right multiplication induces an isometric action
of $A(\Gamma)$ on $(A(\Gamma), d_*)$.
The metric spaces $(A(\Gamma),d_*)$ and $(\Gamma^e,d)$ are quasi-isometric~\cite{KK14}.
It seems that, for some algebraic tools,
$(A(\Gamma),d_*)$ is easier to work with than $(\Gamma^e,d)$.

In \S\ref{sec:star}, we study basic properties of the star length
concerning the prefix order $\le_L$ and the geodesic decomposition
of group elements.
For example, Corollary~\ref{cor:st3} shows that
if $g_1g_2$ is geodesic, then
$\|g_1\|_* + \|g_2\|_* -2 \le  \|g_1g_2\|_*
\le \|g_1\|_* + \|g_2\|_*$.

\subsection{Prefixes of powers of cyclically reduced elements}
Recall that, for a cyclically reduced element $g$,
if $g^u$ is a left cyclic conjugation,
then $u\le_L g^m$ for some $m\ge 1$,
i.e.\ $u$ is a prefix of some power of $g$.
In \S\ref{sec:power}, we study prefixes of powers.
In particular, we show that if $g$ is cyclically reduced and non-split
and if $u\le_L g^m$ for some $m\ge 1$,
then $u=g^n a$ is geodesic for some $0\le n\le m$ and $a\in A(\Gamma)$
with $\|a\|_*\le \|g\|_*+1$ (see Corollary~\ref{cor:po3}).

\subsection{Asymptotic translation length}
In \S\ref{sec:tl}, we prove Theorem B
by using the results in~\S\ref{sec:power}.

\subsection{Uniqueness of quasi-roots}
An element $g$ is called a \emph{quasi-root} of $h$ if there is a decomposition
$$ h=ag^n b$$
for some $n\ge 1$ and $a,b\in A(\Gamma)$ such that $\|h\|=\|a\|+n\|g\|+\|b\|$.
It is called an \emph{$(A,B,r)$-quasi-root}
if $\|a\|\le A$, $\|b\|\le B$ and $\|g\|\le r$
and an \emph{$(A,B,r)^*$-quasi-root}
if $\|a\|_*\le A$, $\|b\|_*\le B$ and $\|g\|_*\le r$.
The conjugates $aga^{-1}$ and $b^{-1}gb$ are called the \emph{leftward-}
and the \emph{rightward-extraction} of the quasi-root $g$, respectively.

In~\cite{LL22}, it is shown that
if $\|h\|\ge A+B+(2|V(\Gamma)|+1)r$, then strongly non-split and primitive
$(A,B,r)$-quasi-roots of $h$
are unique up to conjugacy, and their leftward- and rightward-extractions are unique.
(See \S\ref{sec:star} and \S\ref{sec:quasiroot}
 for the definitions of strongly non-split elements and primitive
 elements.)

In \S\ref{sec:quasiroot}, we extend the above result
to $(A,B,r)^*$-quasi-roots:
if $\|h\|_*\ge 2A+2B+(2|V(\Gamma)|+3)r+2$, then primitive
$(A,B,r)^*$-quasi-roots of $h$ are unique up to conjugacy,
and their leftward- and rightward-extractions are unique.

\subsection{Proof of the acylindricity}
In \S\ref{sec:acyl},
we first compute the acylindricity constants for the action of $A(\Gamma)$
on $(A(\Gamma),d_*)$ (Theorem~\ref{thm:st}) by combining the results
from the previous sections.
Then we prove Theorem A using the quasi-isometry between
$(A(\Gamma),d_*)$ and $(\Gamma^e,d)$.

\subsection{Conventions and notations}

Throughout the paper,
all the group actions are right-actions.

For graphs $\Gamma_1$ and $\Gamma_2$,
the \emph{disjoint union} $\Gamma_1\sqcup \Gamma_2$
is the graph such that
\begin{align*}
V(\Gamma_1\sqcup\Gamma_2) &=V(\Gamma_1)\sqcup V(\Gamma_2),\\
E(\Gamma_1\sqcup\Gamma_2) &=E(\Gamma_1)\sqcup E(\Gamma_2).
\end{align*}
The \emph{join} $\Gamma_1 * \Gamma_2$ is the graph such that
$\overline{\Gamma_1 * \Gamma_2} =  \bar\Gamma_1 \sqcup \bar\Gamma_2$, hence
\begin{align*}
V(\Gamma_1*\Gamma_2) &=V(\Gamma_1)\sqcup V(\Gamma_2),\\
E(\Gamma_1*\Gamma_2) &=E(\Gamma_1)\sqcup E(\Gamma_2)\sqcup
\{\,\{v_1,v_2\}: v_1\in V(\Gamma_1),~v_2\in V(\Gamma_2)\,\}.
\end{align*}
A graph is called a {\em join} if it is the join of two nonempty graphs.
A subgraph that is a join is called a {\em subjoin}.

For $X\subset V(\Gamma)$, $\Gamma[X]$
denotes the subgraph of $\Gamma$ induced by $X$,
i.e.
$$
V(\Gamma[X])=X,\qquad
E(\Gamma[X])=\{\,\{v_1,v_2\}\in E(\Gamma): v_1,v_2\in X\}.
$$

For $g\in A(\Gamma)$, the subgraphs $\Gamma[\supp(g)]$ and $\bar\Gamma[\supp(g)]$
are abbreviated to $\Gamma[g]$ and $\bar\Gamma[g]$, respectively.

For $v\in V(\Gamma)$ and $X\subset V(\Gamma)$,
the sets $\Lk_\Gamma(v)$, $\St_\Gamma(v)$ and $\St_\Gamma(X)$
are defined as follows:
\begin{align*}
\Lk_\Gamma(v)  & = \{v_1\in V(\Gamma): \{v_1,v\}\in E(\Gamma)\},\\
\St_\Gamma(v)  & = \{ v\} \cup \Lk_\Gamma(v),\\
\St_\Gamma(X)  & = \bigcup_{v\in X} \St_\Gamma(v).
\end{align*}
They are called
the \emph{link} of $v$,
the \emph{star} of $v$ and
the \emph{star} of $X$,
respectively.
They will be written as $\Lk(v)$, $\St(v)$ and $\St(X)$ by omitting $\Gamma$
whenever the context is clear.

The \emph{path graph} $P_k=(v_1,v_2,\ldots,v_k)$
is the graph with $V(P_k)=\{v_1,\ldots,v_k\}$
and $E(P_k)=\{ \{v_i,v_{i+1}\}:1 \leqslant i \leqslant k-1\}$,
hence $P_k$ looks like
$\begin{xy}/r.7mm/:
(0,0)="a" *+!U{v_1} *{\bullet},
"a"+(10,0)="a" *+!U{v_2} *{\bullet},
"a"+(20,0)="a" *+!U{v_{k-1}} *{\bullet},
"a"+(10,0)="a" *+!U{v_k} *{\bullet},
(0,0)="a"; "a"+(15,0) **@{-},
"a"+(20,0) *{\cdots},
(25,0)="a"; "a"+(15,0) **@{-},
\end{xy}$.

A \emph{path} in a graph $\Gamma$ is a tuple $(v_1,v_2,\ldots,v_k)$
of vertices of $\Gamma$ such that $\{v_i,v_{i+1}\}\in E(\Gamma)$
for all $1 \leqslant i \leqslant k-1$.
(We do not assume that the vertices or the edges in the path
are mutually distinct.
Hence it means the \emph{walk} in the graph theoretical terminology.)

\EKLnewpage

\section{Lattice structures}
\label{sec:lattice}

In this section we study lattice structures in right-angled Artin groups.

An element of $V(\Gamma)^{\pm 1}=V(\Gamma)\cup V(\Gamma)^{-1}$
is called a \emph{letter}.
A \emph{word} means a finite sequence of letters.
For words $w_1$ and $w_2$, the notation $w_1\equiv w_2$ means that
$w_1$ and $w_2$ coincide as sequences of letters.
A word $w'$ is called a {\em subword} of a word $w$
if $w \equiv w_1 w' w_2$ for (possibly empty) words $w_1$ and  $w_2$.

Suppose that $g\in A(\Gamma)$ is expressed
as a word $w$ on $V(\Gamma)^{\pm 1}$.
The word $w$ is called \emph{reduced} if $w$ is a shortest
word among all the words representing $g$.
In this case, the length of $w$ is called the \emph{word length} of $g$
and denoted by $\|g \|$.

\begin{definition}[support]
For $g\in A(\Gamma)$, the \emph{support} of $g$, denoted $\supp(g)$,
is the set of generators that appear in a reduced word representing $g$.
It is known that $\supp(g)$ is well defined (by~\cite{HM95}),
i.e.\ it does not depend on the choice of a reduced word
representing $g$.
\end{definition}

\begin{definition}[disjointly commute]
We say that $g_1, g_2\in A(\Gamma)$ \emph{disjointly commute},
denoted $g_1\rla g_2$,
if $\supp(g_1)\cap\supp(g_2)=\emptyset$ and  each $v_1\in \supp(g_1)$ commutes with each $v_2\in\supp(g_2)$.
(In particular, the identity element $1\in A(\Gamma)$ disjointly commutes with
any $g\in A(\Gamma)$.)
\end{definition}

The notation $\Gamma[g]$ is an abbreviation of $\Gamma[\supp(g)]$,
the subgraph of $\Gamma$ induced by $\supp(g)$.
From a graph theoretical viewpoint,
$g_1\rla g_2$ means that $\supp(g_1)\cap\supp(g_2)=\emptyset$
and $\bar\Gamma[g_1 g_2] = \bar\Gamma[g_1] \sqcup\bar\Gamma[g_2]$
in the complement graph $\bar\Gamma$
(or equivalently $\Gamma[g_1 g_2] = \Gamma[g_1] *\Gamma[g_2]$ in the graph $\Gamma$).
Recall that $\St_{\bar\Gamma}(\supp(g))$ denotes the star of $\supp(g)$
in the complement graph $\bar\Gamma$.
The following lemma is now obvious.

\begin{lemma}\label{lem:st-s1}
For $g_1,g_2\in A(\Gamma)$, the following are equivalent:
\begin{itemize}
\item[(i)]
$g_1\rla g_2$ in $A(\Gamma)$;

\item[(ii)]
$\St_{\bar\Gamma}(\supp(g_1))\cap \supp(g_2)=\emptyset$.
\end{itemize}
\end{lemma}

Let $w$ be a (non-reduced) word on  $V(\Gamma)^{\pm 1}$.
A subword $v^{\pm1}w_1v^{\mp1}$ of $w$, where $v\in V(\Gamma)$,
is called a \emph{cancellation} of $v$ in $w$
if ${\supp}(w_1)\subset \St_{\Gamma}(v)$,
i.e.\ each $v_1\in\supp(w_1)$ commutes with $v$.
If, furthermore, no letter in $w_1$ is equal to $v$ or $v^{-1}$,
it is called an \emph{innermost cancellation} of $v$ in $w$.
It is known that the following are equivalent:
\begin{enumerate}
\item $w$ is a reduced word;
\item $w$ has no cancellation;
\item $w$ has no innermost cancellation.
\end{enumerate}

Abusing terminology, we do not distinguish between an element $g\in A(\Gamma)$
and a reduced word $w$ representing $g$ if there is no confusion.
For example, if there is a cancellation in $w_1w_2$, where each $w_i$ is a reduced word
representing an element $g_i$, then we just say that there is a cancellation in $g_1g_2$.

\begin{definition}[geodesic decomposition]\label{def:GD}
For $k\ge 1$ and $g,g_1,\ldots,g_k\in A(\Gamma)$, we say that the decomposition
$g=g_1\cdots g_k$ is \emph{geodesic}, or $g_1\cdots g_k$ is \emph{geodesic},
if $\|g\|=\|g_1\|+\cdots+\|g_k\|$.
\end{definition}

If $g_1\cdots g_k$ is geodesic, then
the following are obvious from the definition:
\begin{enumerate}
\item $g_k^{-1}g_{k-1}^{-1}\cdots g_1^{-1}$ is geodesic;
\item $g_pg_{p+1}\cdots g_q$ is geodesic for any $1\le p< q\le k$;
\item $\supp(g_1\cdots g_k)=\supp(g_1)\cup\cdots\cup\supp(g_k)$.
\end{enumerate}

\begin{definition}[prefix order] \label{def:PO}
Let $g=g_1g_2$ be geodesic for $g,g_1,g_2\in A(\Gamma)$.
We say that $g_1$ is a \emph{prefix} (or a \emph{left divisor}) of $g$, denoted $g_1\le_L g$,
and that $g$ is a \emph{right multiple} of $g_1$.
Similarly, we say that $g_2$ is a \emph{suffix} (or a \emph{right divisor}) of $g$,
denoted $g_2\le_R g$,
and that $g$ is a \emph{left multiple} of $g_2$.
\end{definition}

Clearly both $\le_L$ and $\le_R$ are partial orders on $A(\Gamma)$.
The following lemma shows their basic properties.
The proof is straightforward, hence we omit it.

\begin{lemma}\label{lem:PO}
Let $g,g_1,\ldots,g_n, h_1, h_2\in A(\Gamma)$.
\begin{enumerate}
\item $g_1\le_L g_2$ if and only if $g_1^{-1}\le_R g_2^{-1}$.
\item If $gg_1$ and $gg_2$ are geodesic, then $gg_1\le_L gg_2$ if and only if $g_1\le_L g_2$.
\item $g_1\cdots g_n$ is geodesic if and only if
$g_1\cdots g_k\le_L g_1\cdots g_{k+1}$ for all $1\le k\le n-1$.

\item
Suppose $g_1 g_2 = h_1 h_2$ such that both $g_1 g_2$ and $h_1 h_2$ are geodesic.
Then $g_1 \le_L h_1$ if and only if $h_2\le_R g_2$.
\end{enumerate}
\end{lemma}

\begin{definition}[gcd and lcm]\label{def:gcd_lcm}
For $g,h\in A(\Gamma)$,
the symbols $g\wedge_L h$ and $g\vee_L h$ denote
the \emph{greatest common divisor} (gcd)
and the \emph{least common multiple} (lcm) with respect to $\le_L$.
In other words, $g\wedge_L h$ is an element such that
(i) $g\wedge_L h\le_L g$ and  $g\wedge_L h\le_L h$;
(ii) if $u \le_L g$ and  $u\le_L h$ for some $u\in A(\Gamma)$, then $u\le_L g\wedge_L h$.
Similarly,  $g\vee_L h$ is an element such that
(i) $g\le_L g\vee_L h$ and  $h\le_L g\vee_L h$;
(ii) if $g\le_L u$ and  $h\le_L u$ for some $u\in A(\Gamma)$, then $g\vee_L h\le_L u$.

The symbols $g\wedge_R h$ and  $g\vee_R h$ denote the gcd and lcm respectively with respect to
$\le_R$.
\end{definition}

The elements $g\wedge_L h$ and $g\vee_L h$ are unique if they exist.
In Theorem~\ref{thm:lcm} we will show that $g\wedge_L h$ always exists
and that $g\vee_L h$ exists if and only if
$g$ and $h$ admit a common right multiple.

Note that $g$ and $h$ have no nontrivial common prefix
if and only if $g\wedge_L h=1$, i.e.\ the gcd $g\wedge_L h$ exists
and is equal to the identity.
Therefore even though we did not prove yet
the existence of $g\wedge_L h$ for arbitrary $g$ and $h$,
we can safely use the expression $g\wedge_L h=1$.

The following lemma is an easy consequence of the fact that
a word is reduced if and only if it has no innermost cancellation.

\begin{lemma}\label{lem:cc}
Let $u, g, g_1,\ldots,g_k\in A(\Gamma)$.
\begin{itemize}
\item[(i)]
Suppose that $g_1\cdots g_k$ is not geodesic.
Then there exist $1\le p<q\le k$ and $x\in V(\Gamma)^{\pm 1}$ such that
$$
x^{-1}\le_R g_p,\quad
x\le_L g_q,\quad
x\rla g_j~\mbox{for all $p<j<q$}.
$$
Furthermore, if both $g_1\cdots g_{k-1}$
and $g_2\cdots g_{k}$ are geodesic, then
$p=1$ and $q=k$.

\item[(ii)]
Suppose that for each $1\le p<q\le k$, either $g_pg_q$ is geodesic
or $g_pg_{j_1}\cdots g_{j_r}g_q$ is geodesic for some
$p<j_1<\cdots<j_r<q$.
Then $g_1\cdots g_k$ is geodesic.

\item[(iii)] Suppose that $gg$ is geodesic.
For any $n\ge 2$ and $a,b\in A(\Gamma)$, the following are equivalent:
\begin{itemize}
\item[(a)] $agb$ is geodesic;
\item[(b)] $a\underbrace{gg\cdots g}_nb$ is geodesic;
\item[(c)] $ag^nb$ is geodesic.
\end{itemize}
In particular, $g^n=gg\cdots g$ is geodesic
for any $n\ge 2$.

\item[(iv)] Suppose that $g_ig_i$ is geodesic for all $1\le i\le k$
and that $a_1g_1a_2g_2\cdots a_kg_k a_{k+1}$ is geodesic
for some $a_1,\ldots, a_{k+1}\in A(\Gamma)$.
Then $a_1g_1^{n_1}a_2g_2^{n_2}\cdots a_kg_k^{n_k} a_{k+1}$ is geodesic
for any $n_i\ge 1$.
\end{itemize}
\end{lemma}

\begin{myproof}
(i)\ \
Let $w_i$ be a reduced word representing $g_i$ for $i=1,\ldots, k$.
Since $g_1\cdots g_k$ is not geodesic,
the word $w\equiv w_1\cdots w_k$ is not reduced,
hence it has an innermost cancellation.
Since each $w_i$ is reduced, the cancellation must occur
between $x^{-1}$ in $w_p$ and  $x$ in $w_q$
for some $1\le p<q\le k$ and $x\in V(\Gamma)^{\pm 1}$.
Therefore $w_p$ and $w_q$ are of the form
$w_p\equiv w_p'x^{-1}w_p''$ and $w_q\equiv w_q'xw_q''$
such that $x$ disjointly commutes with
$w_p'', w_{p+1},\ldots, w_{q-1}, w_q'$,
hence $x^{-1}\le_R g_p$, $x\le_L g_q$
and $x\rla g_j$ for all $p<j<q$.

If either $p>1$ or $q<k$,
then either $g_2\cdots g_k$ or $g_1\cdots g_{k-1}$
is not geodesic.
Therefore if both $g_1\cdots g_{k-1}$
and $g_2\cdots g_{k}$ are geodesic, then
$p=1$ and $q=k$.

\smallskip
(ii)\ \
Assume that $g_1\cdots g_k$ is not geodesic.
By (i), there exist $1\le p<q\le k$ and $x\in V(\Gamma)^{\pm 1}$
such that
$x^{-1}\le_R g_p$,
$x\le_L g_q$
and $x\rla g_j$  for all $j$ with $p<j<q$.
Therefore none of $g_pg_q$ and  $g_pg_{j_1}\cdots g_{j_r}g_q$ ($p<j_1<\cdots<j_r<q$) is geodesic,
which contradicts the hypothesis.

\smallskip
(iii)\ \
(a) $\Rightarrow$ (b):\ \
Let $h_1=a$, $h_i=g$ for $i=2,\ldots,n+1$ and $h_{n+2}=b$.
Then $h_1,\ldots,h_{n+2}$ satisfy the hypothesis of (ii),
hence $h_1h_2\cdots h_{n+1}h_{n+2}=a\underbrace{g\cdots g}_n b$
is geodesic.

(b) $\Rightarrow$ (a):\ \
Since $a\underbrace{g\cdots g}_n b$ is geodesic, $ag$ and $gb$ are geodesic.
If $agb$ is not geodesic, then there exists $x\in V(\Gamma)^{\pm 1}$
such that $x^{-1}\le_R a$, $x\le_L b$ and $x\rla g$ by (i).
Hence $ag\cdots gb$ is not geodesic, which is a contradiction.

(b) $\Leftrightarrow$ (c):\ \
From (a) $\Rightarrow$ (b) with $a=b=1$,
$g^n=gg\cdots g$ is geodesic, i.e.\ $\|g^n\|=n\|g\|$.
Therefore $\|ag^nb\|=\|a\| + \|g^n\|+\|b\|$
if and only if $\|ag^nb\|=\|a\| + n\|g\|+\|b\|$,
i.e.\ $ag^nb$ is geodesic if and only if $a\underbrace{gg\cdots g}_nb$ is geodesic.

\smallskip
(iv)\ \
Applying (iii)
with $a=a_1$, $g=g_1$, $b=a_2g_2\cdots a_{k+1}$ and $n=n_1$,
we get that  $a_1g_1^{n_1}a_2g_2\cdots a_{k+1}$ is geodesic.
Then applying (iii) with
$a=a_1g_1^{n_1}a_2$, $g=g_2$, $b=a_3g_3\cdots a_{k+1}$ and $n=n_2$,
we get that   $a_1g_1^{n_1}a_2g_2^{n_2}a_3g_3\cdots a_{k+1}$ is geodesic.
Iterating this process, we get that
$a_1g_1^{n_1}a_2g_2^{n_2}\cdots a_kg_k^{n_k} a_{k+1}$ is geodesic.
\end{myproof}

\begin{lemma}\label{lem:geo}
Let $g_1,g_2\in A(\Gamma)$ and $x\in V(\Gamma)^{\pm 1}$.
\begin{itemize}
\item[(i)]
If\/ $g_1g_2$ is not geodesic, then
there exists $y\in V(\Gamma)^{\pm 1}$ such that
$y^{-1}\le_R g_1$ and $y\le_L g_2$.

\item[(ii)]
Let  $g_1g_2$ be geodesic.
If\/ $x\le_L g_1g_2$ and $x \not\le_L g_1$,
then $x\le_L g_2$ and $x\rla g_1$.

\item[(iii)]
Let $g_1g_2$ be geodesic.
If\/ $x\le_R g_1g_2$ and $x \not\le_R g_2$,
then $x\le_R g_1$ and $x\rla g_2$.
\end{itemize}
\end{lemma}

\begin{myproof}
(i)\ \
It follows from Lemma~\ref{lem:cc}(i) with $k=2$.

\smallskip
(ii)\ \
Since $x\not\le_L g_1$, the decomposition $x^{-1}\cdot g_1$ is geodesic.
Since $x\le_L g_1g_2$, the decomposition $x^{-1}\cdot g_1g_2$ is not geodesic.
Since both $x^{-1}\cdot g_1$ and $g_1\cdot g_2$ are geodesic,
there exists $y\in V(\Gamma)^{\pm 1}$ such that
$y^{-1}\le_R x^{-1}$, $y\le_L g_2$ and $y\rla g_1$
(by Lemma~\ref{lem:cc}(i)), hence $x=y$.
Therefore $x\le_L g_2$ and $x\rla g_1$.

\smallskip
(iii)\ \ The proof is analogous to (ii).
\end{myproof}

\begin{lemma}\label{lem:pre2}
Let $g\in A(\Gamma)$ and $x\ne y\in V(\Gamma)^{\pm 1}$
(possibly $y=x^{-1}$).
\begin{itemize}
\item[(i)]
If $x\le_L g$ and $y\le_R g$, then $g=xhy$ is geodesic for some $h\in A(\Gamma)$.
\item[(ii)]
If $x,y\le_L g$, then $x\rla y$ and $g=xyh$ is
geodesic for some $h\in A(\Gamma)$.
\item[(iii)]
If $x,y\le_R g$, then $x\rla y$ and $g=hxy$ is
geodesic for some $h\in A(\Gamma)$.
\end{itemize}
\end{lemma}

\begin{myproof}
(i)\ \
Since $y\le_R g$, $g=g'y$ is geodesic for some $g'\in A(\Gamma)$.
Since $x\le_L g=g'y$, if $x\not\le_L g'$, then $x\le_L y$ (by Lemma~\ref{lem:geo}(ii)),
which contradicts the hypothesis $x\ne y$.
Thus $x\le_L g'$, hence $g'=xh$ is geodesic for some $h\in A(\Gamma)$.
Therefore $g=g'y=xhy$ is geodesic.

\smallskip
(ii)\ \
Since $x\le_L g$, $g=xg'$ is geodesic for some $g'\in A(\Gamma)$.
Since $y\ne x$ (hence $y\not\le_L x$) and $y\le_L xg'$,
we have $y\rla x$ and $y\le_L g'$
(by Lemma~\ref{lem:geo}(ii)),
hence
$g'=yh$ is geodesic for some $h\in A(\Gamma)$.
Therefore $g=xg'=xyh$ is geodesic.

\par\smallskip(iii)\ \
The proof is analogous to (ii).
\end{myproof}

\begin{lemma}\label{lem:lcm1}
Let $g_1,g_2, h_1,h_2, h\in A(\Gamma)$ with both $g_1 g_2$ and $h_1 h_2$ geodesic.

\begin{enumerate}
\item
If $h\wedge_L g_1=h\wedge_L g_2=1$, then $h\wedge_L (g_1g_2)=1$.

\item
If $h\le_L g_1g_2$ and $h\wedge_L g_1=1$,
then $h\rla g_1$ and $h\le_L g_2$.

\item
Let $g_1 g_2=h_1 h_2$.
If $g_1\wedge_L h_1=g_2 \wedge_R h_2=1$,
then $g_1\rla h_1$, $g_1=h_2$ and $g_2=h_1$.
\end{enumerate}
\end{lemma}

\begin{figure}
$\begin{xy}  /r1.8mm/:
( 0,10)="A0", (20,10)="A1", (50,10)="A2",
( 5, 5)="B0", (25, 5)="B1",
(10, 0)="C0", (30, 0)="C1",
"A0"; "A1" **@{-} ?(.5) *+!D{g_1},
      "B0" **@{-} ?(.5) *!R{h_1},
"B0"; "B1" **@{-} ?(.5) *+!D{g_1},
      "C0" **@{-} ?(.5) *+!R{x},
"C0"; "C1" **@{-} ?(.5) *+!D{g_1},
"A1"; "A2" **@{-} ?(.35) *!D{g_2},
      "B1" **@{-} ?(.5) *!R{h_1},
"B1"; "A2" **@{-} ?(.3) *!D{g_2'},
      "C1" **@{-} ?(.5) *+!R{x},
"C1"; "A2" **@{-} ?(.3) *!DR{g_2''},
\end{xy}$
\caption{van Kampen diagram for Lemma~\ref{lem:lcm1}(ii) }
\label{fig:prefix1}
\end{figure}
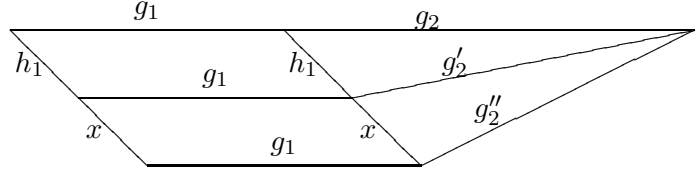

\begin{myproof}
(i)\ \ If $h\wedge_L (g_1g_2)\ne 1$, then there exists $x\in V(\Gamma)^{\pm 1}$ such that
$x\le_L h$ and $x\le_L g_1g_2$.
Since $x\le_L h$ and  $h\wedge_L g_1=h\wedge_L g_2=1$, we have
$x\not\le_L g_1$ and $x\not\le_L g_2$.
Since $x\le_L g_1g_2$ and $x\not\le_L g_1$,
we have $x\le_L g_2$ by Lemma~\ref{lem:geo}(ii),
which is a contradiction.

\smallskip
(ii)\ \
We use induction on $\|h\|$.

If $\|h\|=0$, there is nothing to prove.
If $\|h\|=1$, it holds by Lemma~\ref{lem:geo}(ii).

Assume $\|h\|\ge 2$. Then $h=h_1x$ is geodesic
for some $h_1\in A(\Gamma)$ and $x\in V(\Gamma)^{\pm 1}$.
See Figure~\ref{fig:prefix1}.
Notice that $h_1\le_L g_1g_2$ and $h_1 \wedge_L g_1=1$.
By the induction hypothesis, $h_1\rla g_1$ and $h_1\le_L g_2$,
hence  $g_2=h_1g_2'$ is geodesic for some $g_2'\in A(\Gamma)$.

Since $h_1\rla g_1$, we have  $g_1g_2=g_1h_1g_2'=h_1g_1g_2'$.
Since $g_1 g_2$ is geodesic, so is $h_1 g_1 g_2'$.
Since $h_1x=h\le_L g_1g_2= h_1g_1g_2'$ and
both $h_1x$ and $h_1g_1g_2'$ are geodesic,
we have $x\le_L g_1g_2'$.

Observe $x\not\le_L g_1$.
(If $x\le_L g_1$, then $x\rla h_1$ because $h_1\rla g_1$.
Since $h=h_1x=xh_1$ is geodesic, we have $x\le_L h$.
Thus $x$ is a common prefix of $g_1$ and $h$, which contradicts the hypothesis $h\wedge_L g_1=1$.)
By Lemma~\ref{lem:geo}(ii), we get $x\rla g_1$ and $x\le_L g_2'$,
hence $g_2'=xg_2''$ is geodesic for some $g_2''\in A(\Gamma)$.

Since $g_2=h_1g_2'=h_1xg_2''=hg_2''$ and since $hg_2''$ is geodesic,
we have $h\le_L g_2$.
On the other hand, since $h_1\rla g_1$ and $x\rla g_1$,
we have $h=h_1x\rla g_1$.

(iii) \ \
Since $g_1\le_L h_1 h_2$, $h_1\le_L g_1g_2$ and $g_1\wedge_L h_1=1$,
we have $g_1\rla h_1$, $g_1\le_L h_2$ and $h_1\le_L g_2$
(by (ii)).
Thus $h_2=g_1h_2'$ and $g_2 = h_1 g_2'$ are geodesic for some $g_2', h_2'\in A(\Gamma)$.

Observe $g_1 h_1 g_2'=g_1 g_2=h_1 h_2=h_1 g_1h_2'=g_1h_1 h_2'$,
which implies $g_2'=h_2'$.
Since $g_2\wedge_R h_2=1$, we have $g_2'=h_2'=1$, hence $g_1=h_2$ and $g_2=h_1$.
\end{myproof}

The following is the main result of this section.

\begin{theorem}\label{thm:lcm}
For $g_1,g_2\in A(\Gamma)$,
the gcd $g_1\wedge_L g_2$ always exists
and the lcm $g_1\vee_L g_2$ exists if and only if $g_1$ and $g_2$ have a common right multiple.

More precisely,
if\/ $g_0$ is a maximal common prefix of $g_1$ and $g_2$,
hence $g_1=g_0g_1'$ and $g_2=g_0g_2'$ are geodesic for some $g_1',g_2'\in A(\Gamma)$
with $g_1'\wedge_L g_2'=1$,
then the following hold.
\begin{itemize}
\item[(i)]
$g_1$ and $g_2$ have a common right multiple if and only if $g_1'\rla g_2'$.
In this case, $g_1\vee_L g_2$ exists and
$g_1\vee_L g_2=g_1g_2'=g_2g_1'=g_0g_1'g_2'$.
In particular, $\supp(g_1\vee_L g_2)=\supp(g_1)\cup\supp(g_2)$.

\item[(ii)]
$g_1\wedge_L g_2=g_0$.
\end{itemize}
\end{theorem}

\begin{figure}
$$
\begin{xy}  /r1.5mm/:
(0,0)="A0",
"A0"+(15,  0)="A1",
"A1"+(20,  0)="A2",
"A2"+(15,  0)="A3",
"A1"+(10, 10)="T1",
"A1"+(10,-10)="B1",
"A0";
"T1" **@{-} ?(.6) *!DR{g_1},
"A1" **@{-} ?(.7) *!D {g_0},
"B1" **@{-} ?(.6) *!UR{g_2},
"A1";
"T1" **@{-} ?(.5) *!UL{g_1'},
"B1" **@{-} ?(.5) *!DL{g_2'},
"A2";
"T1" **@{-} ?(.5) *!UR{h_1'},
"B1" **@{-} ?(.5) *!DR{h_2'},
"A3";
"T1" **@{-} ?(.6) *!DL{h_1},
"A2" **@{-} ?(.7) *!D {h_0},
"B1" **@{-} ?(.6) *!UL{h_2},
\end{xy}
$$
\caption{van Kampen diagram for Theorem~\ref{thm:lcm}}
\label{fig:gcdlcm}
\end{figure}
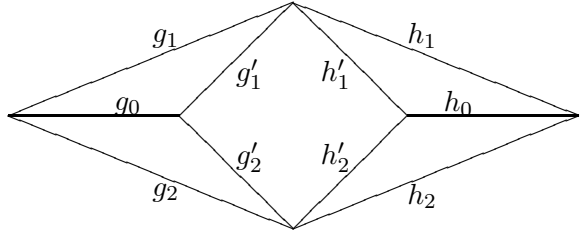

\begin{myproof}
(i)\ \
Assume $g_1'\rla g_2'$.
Then $g_1'g_2'$ is geodesic (otherwise
there exists $x\in V(\Gamma)^{\pm 1}$
such that $x^{-1}\le_R g_1'$ and $x\le_L g_2'$ by Lemma~\ref{lem:geo}(i),
hence $g_1'$ and $g_2'$ do not disjointly commute).
Since $g_0g_1'$, $g_0g_2'$ and $g_1'g_2'$ are all geodesic,
$g_0g_1'g_2'$ is geodesic
(by Lemma~\ref{lem:cc}(ii)).
Therefore $g_0g_1'g_2'=g_1g_2'=g_2g_1'$ is a common right multiple
of $g_1$ and $g_2$.

Conversely, assume that $g_1$ and $g_2$ have a common right multiple $h$.
Then $h=g_1h_1=g_2h_2$ are geodesic for some $h_1,h_2\in A(\Gamma)$.
We need to show that $g_1'\rla g_2'$.

Let $h_0$ be a maximal common suffix of $h_1$ and $h_2$.
Then $h_1=h_1'h_0$ and $h_2=h_2'h_0$ are geodesic for some $h_1',h_2'\in A(\Gamma)$
with $h_1' \wedge_R h_2'=1$.
See Figure~\ref{fig:gcdlcm}.
Notice that  $g_1'h_1'=g_2'h_2'$ and that $g_1',g_2',h_1',h_2'$ satisfy
the hypotheses of Lemma~\ref{lem:lcm1}(iii).
Therefore $g_1'\rla g_2'$.

Lemma~\ref{lem:lcm1}(iii) also claims $g_1'=h_2'$ and $g_2'=h_1'$,
hence $h=g_1 h_1 = g_0 g_1' h_1' h_0 = g_0 g_1' g_2' h_0$.
Therefore $g_0 g_1'g_2'$ is a prefix of any common
right multiple $h$ of $g_1$ and $g_2$, namely,
$g_1 \vee_L g_2= g_0g_1'g_2'$.
It follows immediately that
$\supp(g_1\vee_L g_2)=\supp(g_1)\cup\supp(g_2)$.
Since $g_0g_1'g_2' = g_0g_2'g_1'$, we have
$g_1\vee_L g_2 =g_0g_1'g_2' =g_1g_2'=g_2g_1'$.

\smallskip\noindent
(ii)\ \
Let $u_0$ be a common prefix of $g_1$ and $g_2$.
Since $g_1$ and $g_2$ are common right multiples of $g_0$ and $u_0$,
the lcm $g_0\vee_L u_0$ exists (by (i)) and is a prefix of both $g_1$ and $g_2$,
hence $g_0\vee_L u_0$ is a common prefix of $g_1$ and $g_2$.
Since $g_0\vee_L u_0$ is a right multiple of $g_0$
and since $g_0$ is a maximal common prefix of $g_1$ and $g_2$,
we have $g_0=g_0\vee_L u_0$, hence $u_0\le_L g_0$.
Therefore $g_0=g_1\wedge_L g_2$.
\end{myproof}

Obviously we can replace $(\wedge_L,\vee_L)$ in Theorem~\ref{thm:lcm}  with
$(\wedge_R,\vee_R)$ as follows.

\begin{theorem}\label{thm:lcm2}
For $g_1,g_2\in A(\Gamma)$,
the gcd $g_1\wedge_R g_2$ always exists
and the lcm $g_1\vee_R g_2$ exists if and only if
$g_1$ and $g_2$ have a common left multiple.

More precisely,
if\/ $g_0$ is a maximal common suffix of $g_1$ and $g_2$,
hence $g_1=g_1'g_0$ and $g_2=g_2'g_0$ are geodesic
for some $g_1',g_2'\in A(\Gamma)$
with $g_1'\wedge_R g_2'=1$,
then the following hold.
\begin{itemize}
\item[(i)]
$g_1$ and $g_2$ have a common left multiple if and only if
$g_1'\rla g_2'$.
In this case, $g_1\vee_R g_2$ exists and
$g_1\vee_R g_2=g_2'g_1=g_1'g_2=g_1'g_2'g_0$.
In particular, $\supp(g_1\vee_R g_2)=\supp(g_1)\cup\supp(g_2)$.

\item[(ii)]
$g_1\wedge_R g_2=g_0$.
\end{itemize}
\end{theorem}

Observe that the gcds $g_1\wedge_L
g_2$ and $g_1\wedge_R g_2$ exist for any $g_1,g_2\in A(\Gamma)$
by the above theorems.

The following lemma is obvious, hence we omit the proof.
\begin{lemma}\label{lem:gcd1}
Let $g_1,g_2\in A(\Gamma)$.
\begin{itemize}
\item[(i)]
$(g_1\wedge_L g_2)^{-1}=g_1^{-1}\wedge_R g_2^{-1}$.

\item[(ii)]
If $g_1=g_0g_1'$ and $g_2=g_0g_2'$ are geodesic,
then $g_1\wedge_L g_2=g_0(g_1'\wedge_L g_2')$.
In particular, if $g_1'\wedge_L g_2'=1$, then $g_1\wedge_L g_2=g_0$.
\item[(iii)]
If $g_1\le_L g_2$, then $(h\wedge_L g_1)\le_L (h\wedge_L g_2)$
for any $h\in A(\Gamma)$.
\item[(iv)]
The statements analogous to (ii) and (iii) also hold for $(\le_R, \wedge_R)$.
\end{itemize}
\end{lemma}

\begin{lemma}\label{lem:gcd2}
Let $g_1,g_2,h\in A(\Gamma)$ with $g_1g_2$ geodesic.
\begin{enumerate}
\item If $h\rla g_1$, then $h\wedge_L (g_1g_2)=h\wedge_L g_2$.
\item If $\supp(h)\cap\supp(g_2)=\emptyset$, then $h\wedge_L (g_1g_2)=h\wedge_L g_1$.
\item If $h\le_L g_1 g_2$ and $h\rla g_1$, then $h\le_L g_2$.
\item If $h\le_L g_1 g_2$ and $\supp(h)\cap\supp(g_2)=\emptyset$, then $h\le_L g_1$.
\item The statements analogous to (i)--(iv)
also hold for $(\le_R, \wedge_R)$.
\end{enumerate}
\end{lemma}

\begin{myproof}
(i)\ \
Let $h_0=h\wedge_L g_2$. Then $h=h_0h'$ and $g_2=h_0g_2'$ are
geodesic for some $h',g_2'\in A(\Gamma)$
with $h'\wedge_L g_2'=1$.
Since $h\rla g_1$ and $h=h_0h'$ is geodesic,
we have $h_0\rla g_1$ and $h'\rla g_1$,
hence $h'\wedge_L g_1=1$.

Notice that $g_1g_2'$ is geodesic because $g_1g_2(=g_1h_0g_2')=h_0g_1g_2'$ is geodesic.
Since  $h'\wedge_L g_1=h'\wedge_L g_2'=1$,
we have $h'\wedge_L (g_1g_2')=1$ (by Lemma~\ref{lem:lcm1}(i)).
Therefore by Lemma~\ref{lem:gcd1}.(ii)
\begin{align*}
h\wedge_L (g_1g_2)
&= (h_0 h') \wedge_L (g_1 h_0 g_2') = (h_0 h') \wedge_L (h_0 g_1 g_2')\\
&= h_0 (h' \wedge_L (g_1 g_2')) = h_0 =h\wedge_L g_2.
\end{align*}

\smallskip
(ii)\ \
Let $h_0=h\wedge_L g_1$. Then $h=h_0h'$ and $g_1=h_0g_1'$ are
geodesic for some $h',g_1'\in A(\Gamma)$
with $h'\wedge_L g_1'=1$.
Since $\supp(h)\cap\supp(g_2)=\emptyset$ and since $h=h_0h'$ is geodesic,
we have $\supp(h')\cap\supp(g_2)=\emptyset$,
hence $h'\wedge_L g_2=1$.

Notice that $g_1'g_2$ is geodesic because $g_1g_2=h_0g_1'g_2$ is geodesic.
Since  $h'\wedge_L g_1'=h'\wedge_L g_2=1$,
we have $h'\wedge_L (g_1'g_2)=1$ (by Lemma~\ref{lem:lcm1}(i)).
Therefore by Lemma~\ref{lem:gcd1}.(ii)
$$h\wedge_L (g_1g_2) = (h_0 h')\wedge_L (h_0g_1'g_2)
= h_0 (h' \wedge_L (g_1' g_2))=h_0 =h\wedge_L g_1.$$

\smallskip
(iii) and (iv) are direct consequences of (i) and (ii), respectively.

(v)\ \
The proof is analogous to (i)--(iv).
\end{myproof}

\begin{corollary}\label{cor:clo}
Suppose that a set $C\subset A(\Gamma)$ satisfies the following conditions.
\begin{enumerate}
\item[(P1)]
$C$ is prefix-closed, i.e.\
if $g\in C$ and $h\le_L g$, then $h\in C$ .

\item[(P2)]
For $g\in A(\Gamma)$ and $x,y\in V(\Gamma)^{\pm 1}$ such that both $gx$ and $gy$ are geodesic,
if $gx, gy\in C$ and $x\rla y$,
then $gxy\in C$.
\end{enumerate}
Then $C$ is lcm-closed,
i.e.~if $g_1,g_2\in C$
and $g_1\vee_L g_2$ exists, then $g_1\vee_L g_2\in C$.
\end{corollary}

\begin{myproof}
Let $g_1,g_2\in C$ such that $g_1\vee_L g_2$ exists.
Let $g_0=g_1\wedge_L g_2$.
Then
$$
g_1=g_0g_1'\quad\mbox{and}\quad
g_2=g_0g_2'
$$
are geodesic for some $g_1',g_2'\in A(\Gamma)$.
By Theorem~\ref{thm:lcm},
$g_1'\rla g_2'$ and $g_1\vee_L g_2=g_0g_1'g_2'$.

We use induction on $\|g_1'\|+\|g_2'\|$.
If $\|g_1'\|=0$ or $\|g_2'\|=0$,
then $g_1\vee_L g_2$ is either $g_2$ or $g_1$, respectively,
hence there is nothing to prove.
If $\|g_1'\|=\|g_2'\|=1$, then $g_1\vee_L g_2=g_0g_1'g_2'\in C$ by (P2).
Therefore we may assume $\|g_1'\|+\|g_2'\|\ge 3$ and $\|g_1'\|,\ \|g_2' \|\ge 1$.

Then $g_1'=g_1'' x_1$ and $g_2'=g_2'' x_2$ are geodesic for some $g_1'',g_2''\in A(\Gamma)$
and $x_1,x_2\in V(\Gamma)^{\pm1}$.
Thus
$$g_1=g_0g_1'' x_1\quad\mbox{and}\quad
g_2=g_0g_2'' x_2$$
are geodesic, where $g_1'' x_1 \rla g_2'' x_2$.

Since $g_1, g_2\in C$, we have $g_0 g_1'', g_0g_2''\in C$ by (P1),
hence by the induction hypothesis we have
$$g_1\vee_L (g_0g_2'') = g_0 g_1'' g_2'' x_1 \in C \quad\mbox{and}\quad
g_0g_1''\vee_L g_2 = g_0 g_1'' g_2'' x_2 \in C.$$
Therefore  $g_1\vee_L g_2=g_0g_1''g_2'' x_1 x_2 \in C$ by (P2).
\end{myproof}

\EKLnewpage

\section{Cyclic conjugations}
\label{sec:conj}

\begin{definition}[cyclically reduced]
An element $g\in A(\Gamma)$ is called \emph{cyclically reduced}
if it has the minimal word length in its conjugacy class.
\end{definition}

Servatius~\cite[Proposition on p.\,38]{Ser89} showed that
every $g\in A(\Gamma)$ has a unique geodesic decomposition
$$g=u^{-1}hu$$
with $h$ cyclically reduced.
The following lemma shows that
$u$ is determined from $g$ by $u=g\wedge_R g^{-1}$.

\begin{lemma}\label{lem:u-cu}
Let $g, h, u\in A(\Gamma)$.
\begin{enumerate}
\item
If $g=u^{-1}hu$ is geodesic with $h$ cyclically reduced,
then $u=g\wedge_R g^{-1}$.
\item
$g$ is cyclically reduced if and only if $g\wedge_R g^{-1}=1$.
\end{enumerate}
\end{lemma}

\begin{myproof}
(i)\ \
We have two geodesic decompositions
$g=u^{-1}hu$ and $g^{-1}=u^{-1}h^{-1}u$.
By Lemma~\ref{lem:gcd1}, it suffices to show
$(u^{-1}h)\wedge_R(u^{-1}h^{-1})=1$
or equivalently $(hu)\wedge_L (h^{-1}u)=1$.

Assume $(hu)\wedge_L (h^{-1}u)\ne1$.
Then there exists $x\in V(\Gamma)^{\pm 1}$ with
$x\le_L hu$ and $x\le_L h^{-1}u$.

If $x\not\le_L h$, then $x\rla h$ (hence $x\rla h^{-1}$)
and $x\le_L u$ (by Lemma~\ref{lem:geo}(ii)).
Let $u=xu_1$ be geodesic for some $u_1\in A(\Gamma)$.
Then $g=u^{-1}hu=u_1^{-1}x^{-1}hxu_1=u_1^{-1}hu_1$.
This contradicts that $g=u^{-1}hu$ is geodesic.
Therefore $x\le_L h$.
By the same reason, $x\le_L h^{-1}$,
hence $x^{-1}\le_R h$.

Since $x\le_L h$ and $x^{-1}\le_R h$,
$h=xh_1x^{-1}$ is geodesic for some $h_1$ (by Lemma~\ref{lem:pre2}(i)),
which contradicts that $h$ is cyclically reduced.
Therefore $(hu)\wedge_L (h^{-1}u)=1$.

\smallskip
(ii) It follows from (i).
\end{myproof}

\begin{definition}[starting set]
For $g\in A(\Gamma)$, the {\em starting set} $S(g)$ of $g$ is defined as
$$S(g)=\{ x\in V(\Gamma)^{\pm 1} : x\le_L g\}.$$
\end{definition}

\EKLnewpage

\begin{lemma}\label{lem:cred1}
The following hold.
\begin{enumerate}
\item
For any $g\in  A(\Gamma)$, the following are equivalent.
\begin{enumerate}
\item $g$ is cyclically reduced.

\item There is no geodesic decomposition such as $g=u^{-1}hu$,
where $u,h\in A(\Gamma)$ with $u\ne 1$.
\item For any geodesic decomposition $g=g_1g_2$,
$g_2g_1$ is geodesic.
\item
For any $g_1\in A(\Gamma)$ with $g_1\le_L g$, $gg_1$ is geodesic.

\item $g^n=gg\cdots g$ is geodesic (i.e.\ $\|g^n\|=n\|g\|$) for some $n\ge 2$.
\item $g^n=gg\cdots g$ is geodesic (i.e.\ $\|g^n\|=n\|g\|$) for all $n\ge 2$.

\item $g^n$ is cyclically reduced for some $n\ge 2$.
\item $g^n$ is cyclically reduced for all $n\ge 2$.
\end{enumerate}

\item
Let $g_1\cdots g_k$ be geodesic
(i.e.\ $\|g_1\cdots g_k\|=\|g_1\|+\cdots+\|g_k\|$).
Then $g_1^{n_1}\cdots g_k^{n_k}$ is geodesic
(i.e.\ $\|g_1^{n_1}\cdots g_k^{n_k}\|=\|g_1^{n_1}\|+\cdots+\|g_k^{n_k}\|$)
for any positive integers $n_i$.

\item
For any $g\in A(\Gamma)$ and $n\ge 2$,
$\supp(g^n)=\supp(g)$ and $S(g^n)=S(g)$.
\end{enumerate}
\end{lemma}

\begin{myproof}
(i)\ \
The equivalences between (a), (b), (c), (e), (f)
are easy to prove. For example, see~\cite[Lemma 2.1]{LL22}.
We show the remaining equivalences
assuming the known equivalences.

\smallskip
(a) $\Rightarrow$ (d): Assume that $g_1\le_L g$ but $gg_1$ is not geodesic.
Then there exists a letter $x\in V(\Gamma)^{\pm 1}$ such that
$x\le_R g$ and $x^{-1}\le_L g_1$ (by Lemma~\ref{lem:geo}(i)).
Since $g_1 \le_L g$, we have $x^{-1}\le_L g$, hence $x\le_R g^{-1}$.
Now $x\le_R g\wedge_R g^{-1}$, hence $g\wedge_R g^{-1}\ne 1$.
By Lemma~\ref{lem:u-cu}(ii), $g$ is not cyclically reduced.

\smallskip
(d) $\Rightarrow$ (e): Since $g\le_L g$, $gg$ is geodesic.

\smallskip
(f) $\Rightarrow$ (h):
Let $n\ge 2$.
Since $\|g^{2n}\|=2n\|g\|$ and $\|g^n\|=n\|g\|$, we have $\|g^{2n}\|=2\|g^n\|$,
hence $g^n\cdot g^n$ is geodesic.
Because (a) and (e) are equivalent,
$g^n$ is cyclically reduced.

\smallskip
(h) $\Rightarrow$ (g): It is obvious.

\smallskip
(g) $\Rightarrow$ (a):
Assume that $g$ is not cyclically reduced.
Then $g=u^{-1}hu$ is geodesic for some $u,h\in A(\Gamma)$
such that $u\ne 1$ and $h$ is cyclically reduced~\cite{Ser89}.
Observe that $hh$ is geodesic (by (a) $\Leftrightarrow$ (f)).
Hence $g^n$ has a geodesic decomposition
$u^{-1} h^{n}u$ for any $n\ge 2$ (by Lemma~\ref{lem:cc}).
Therefore $g^n$ is not cyclically reduced for any $n\ge 2$
(by (a) $\Leftrightarrow$ (b)).

\medskip
(ii)\ \
Let $g_i=u_i^{-1}h_iu_i$ be a geodesic decomposition of $g_i$
with $h_i$ cyclically reduced for $i=1,\ldots, k$.
Then
\begin{align*}
g_1\cdots g_k
&= u_1^{-1}h_1u_1\cdots u_k^{-1}h_ku_k,\\
g_1^{n_1}\cdots g_k^{n_k}
&= u_1^{-1}h_1^{n_1}u_1 \cdots u_k^{-1}h_k^{n_k}u_k.
\end{align*}
In particular,
$u_1^{-1}h_1u_1\cdots u_k^{-1}h_ku_k$
is geodesic because $g_1\cdots g_k$ and each $g_i=u_i^{-1}h_iu_i$ are geodesic.
Notice that each $h_i h_i$ is geodesic by (i).
Applying Lemma~\ref{lem:cc}(iv), we get that
$u_1^{-1}h_1^{n_1}u_1 \cdots u_k^{-1}h_k^{n_k}u_k$ is geodesic.
Therefore $g_1^{n_1}\cdots g_k^{n_k}$ is geodesic.

\medskip
(iii)\ \
Let $g=u^{-1}hu$ be a geodesic decomposition of $g$
with $h$ cyclically reduced. Then
\begin{equation}\label{eq:g3}
g^n
= u^{-1}h^n u
= u^{-1}h\cdots h u
\end{equation}
are each geodesic decompositions of $g^n$ (by Lemma~\ref{lem:cc}(iii)).

Notice that $\supp(u)=\supp(u^{-1})$
and that if $g_1\cdots g_k$ is a geodesic decomposition, then
$\supp(g_1\cdots g_k)=\supp(g_1)\cup\cdots\cup\supp(g_k)$.
Therefore $\supp(g)=\supp(u)\cup\supp(h)=\supp(g^n)$ from \eqref{eq:g3}.

Observe that $x\le_L h^n$ if and only if $x\le_L h$:
if $x\le_L h$, then it is obvious that $x\le_L h^n$;
if $x\not\le_L h$, then $x\not\le_L h^n$
(otherwise, $x\le_L h^n=h\cdot h^{n-1}$ implies
$x\rla h$ and $x\le_L h^{n-1}$, which contradicts that
$\supp(h)=\supp(h^{n-1})$).

Since $g=u^{-1}hu$ is geodesic,
$x\le_L g$ if and only if one of the following holds
(by Lemma~\ref{lem:geo}(ii)):
(i) $x\le_L u^{-1}$;
(ii) $x\le_L h$ and $x\rla u^{-1}$;
(iii) $x\le_L u$ and $x\rla u^{-1}h$.
Notice that (iii) cannot happen.
Since $x\le_L h$ if and only if $x\le_L h^n$,
we can conclude that $x\le_L g$ if and only if $x\le_L g^n$.
Therefore $S(g)=S(g^n)$.
\end{myproof}

\begin{definition}[cycling, cyclic conjugation]
Let $g\in A(\Gamma)$ be cyclically reduced.
\begin{enumerate}
\item
For a letter $x\in V(\Gamma)^{\pm 1}$, the conjugation $g^x=x^{-1}gx$
is called a \emph{left} (resp.\ \emph{right}) \emph{cycling}
if $x\le_L g$ (resp.\ $x^{-1}\le_R g$).
Left and right cyclings are collectively called \emph{cyclings}.
\item
For an element $u\in A(\Gamma)$, the conjugation $g^u=u^{-1}gu$ is
called a \emph{cyclic conjugation} of $g$ by $u$ if
$\|g^u\|=\|g\|$ and $\supp(u)\subset\supp(g)$.
A cyclic conjugation $g^u$ is called a \emph{left} (resp.\ \emph{right})
\emph{cyclic conjugation} if $gu$ (resp.\ $u^{-1}g$) is geodesic.
\end{enumerate}
\end{definition}

For $g\in A(\Gamma)$ and $x\in V(\Gamma)^{\pm 1}$,
if $g^x$ is a left cycling, i.e. $x\le_L g$,
then $g=xh$ is geodesic for some $h\in A(\Gamma)$ and
$g^x=x^{-1}gx=hx$ is geodesic.
Therefore the left cycling $g^x$ is obtained
from $g=xh$ by moving the first letter $x$ to the last.
Similarly, if $g^x$ is a right cycling, then $g^x$ is
obtained from $g=hx^{-1}$ by moving the last letter $x^{-1}$
to the first.

If $g^x$ is a cycling,
then it is easy to see that $\|g^x\|=\|g\|$ and $\supp(x)\subset\supp(g)$,
hence $g^x$ is a cyclic conjugation.
Conversely, we will show in Lemma~\ref{lem:cy2}
that a cyclic conjugation $g^u$ is obtained by iterated application of cyclings.

If $g\in A(\Gamma)$ is cyclically reduced and $g^u$ is a cyclic conjugation,
then $\|g^u\|=\|g\|$, hence $g^u$ is also cyclically reduced.

\begin{lemma}\label{lem:cy}
Let $g\in A(\Gamma)$ and $x,y\in V(\Gamma)^{\pm 1}$ with $g$ cyclically reduced.

\begin{itemize}

\item[(i)]
The conjugation $g^x$ cannot be both a left cycling and a right cycling.

\item[(ii)] Let $y\ne x^{-1}$.
If $g^x$ and $(g^x)^y$ are cyclings of different type, then $x\rla y$.

\item[(iii)]
Let $x\rla y$.
If both $g^x$ and $(g^x)^y$ are cyclings, then so are $g^y$ and $(g^y)^x$.

\item[(iv)]
Let $x\rla y$.
If both $g^x$ and $g^y$ are cyclings, then so are $(g^x)^y$ and $(g^y)^x$.
\end{itemize}
In (iii) and (iv), the types of cyclings depend only on the conjugating letters.
For example, if $g^x$ is a left cycling, then $(g^y)^x$ is also a left cycling, and so on.
\end{lemma}

\begin{myproof}

\smallskip(i)\ \
If $g^x$ is both a left cycling and a right cycling,
then $x\le_L g$ and $x^{-1}\le_R g$,
hence $g=xhx^{-1}$ is geodesic for some $h\in A(\Gamma)$
(by Lemmas~\ref{lem:pre2}(i)).
Thus $g$ is not cyclically reduced
(by Lemma~\ref{lem:cred1}(i)).

\smallskip(ii)\ \
Assume that $g^x$ is a left cycling and $(g^x)^y$ is a right cycling.
(An analogous argument applies to the case
where $g^x$ is a right cycling and $(g^x)^y$ is a left cycling.)

Since $g^x$ is a left cycling, we have $x\le_L g$,
hence $g=xh$ is geodesic for some $h\in A(\Gamma)$.
Notice that $g^x=hx$ is geodesic.
Since $(g^x)^y$ is a right cycling, $y^{-1}\le_R g^x=hx$.

Since $y^{-1}\ne x$, we have $x\rla y^{-1}$ (by Lemma~\ref{lem:geo}(iii)) and hence $x\rla y$.

\smallskip(iii) and (iv)\ \
Assume that $g^x$ and $(g^x)^y$ are left cyclings,
hence $x\le_L g$ and $y\le_L g^x$.
Since $x\le_L g$, $g=xh_1$ is geodesic for some $h_1\in A(\Gamma)$, hence $g^x=h_1x$ is also geodesic.
Since $y\le_L g^x=h_1x$ and $y\rla x$ (hence $y\not\le_L x$),
we have $y\le_L h_1$,
hence $h_1=yh_2$ is geodesic for some $h_2\in A(\Gamma)$.
Now we know that
$$g=xh_1=xyh_2=yxh_2$$ and
$g^y=xh_2y$ are all geodesic,
hence $y\le_L g$ and $x\le_L g^y$.
This means that $g^y$ and $(g^y)^x$ are left cyclings.

For the other cases, it is easy to see that
$g$ has a geodesic decomposition as one of
$xyh$, $xhy^{-1}$, $yhx^{-1}$ and $hx^{-1}y^{-1}$
depending on the types of cyclings,
from which the conclusions follow.
\end{myproof}

\begin{lemma}\label{lem:cy2}
Let $g, u, u_1,u_2\in A(\Gamma)$ with $g$ cyclically reduced.

\begin{itemize}
\item[(i)] The following are equivalent:
\begin{itemize}
\item[(a)]
$g^u$ is a cyclic (resp.\ left cyclic, right cyclic) conjugation;

\item[(b)]
there exists a reduced word $w_0\equiv y_1\cdots y_k$ representing $u$
such that $(g^{y_1\cdots y_{i-1}})^{y_i}$ is a cycling
(resp.\ left cycling, right cycling)
for all $1\le i\le k$;

\item[(c)]
for any reduced word $w\equiv x_1\cdots x_k$ representing $u$,
$(g^{x_1\cdots x_{i-1}})^{x_i}$ is a cycling
(resp.\ left cycling, right cycling)
for all $1\le i\le k$.
\end{itemize}
In particular, if $g^u$ is a cyclic conjugation, then $\supp(g^u)=\supp(g)$.

\item[(ii)] Let $u=u_1u_2$ be geodesic.
Then $g^u$ is a cyclic (resp.\ left cyclic, right cyclic) conjugation if and only if
both $g^{u_1}$ and $(g^{u_1})^{u_2}$ are
cyclic (resp.\ left cyclic, right cyclic) conjugations.

\item[(iii)]
If $g^{u_1}$ and $g^{u_2}$ are
cyclic (resp.\ left cyclic, right cyclic) conjugations
and $u_1\vee_L u_2$ exists,
then $g^{u_1\vee_L u_2}$ is also a
cyclic (resp.\ left cyclic, right cyclic) conjugation.

\item[(iv)]
Let $u_1 \rla u_2$.
Suppose that $g^{u_1}$ and $g^{u_2}$ are a left cyclic conjugation and a right cyclic conjugation, respectively.
Then $(g^{u_2})^{u_1}$ and $(g^{u_1})^{u_2}$ are a left cyclic conjugation and a right cyclic conjugation,
respectively.
Moreover, $u_2^{-1}gu_1$ is geodesic.

\item[(v)]
Suppose that $g^u$ is a cyclic conjugation.
Then there is a geodesic decomposition $u=u_1u_2$
such that $u_1\rla u_2$ and
$g^{u_1}$ (resp.~$g^{u_2}$) is a left
(resp.~right) cyclic conjugation.
Moreover, $u_2^{-1}gu_1$ is geodesic.
\end{itemize}
\end{lemma}

\begin{myproof}
The statements (i)--(iii) concern three types of cyclic conjugations.
We prove only the case of cyclic conjugations.
The other cases (i.e.\ left and right cyclic conjugations)
can be proved analogously.

We use the following claim.

\smallskip\noindent\textbf{Claim 1.}\ \
If $g^u$ is a cyclic conjugation for some $u\in A(\Gamma)\backslash\{1\}$,
then there exists $y_1\in V(\Gamma)^{\pm 1}$ and $u_1\in A(\Gamma)$
such that $u=y_1u_1$ is geodesic,
$g^{y_1}$ is a cycling and $(g^{y_1})^{u_1}$ is
a cyclic conjugation.

\begin{proof}[Proof of Claim 1]
Since $\|g^u\|=\|g\|$, the decomposition $u^{-1}gu$ is not geodesic.
If both $u^{-1}g$ and $gu$ are geodesic,
then there exists $x\in V(\Gamma)^{\pm1}$ such that
$x^{-1}\le_R u^{-1}$, $x\le_L u$ and $x\rla g$ (by Lemma~\ref{lem:cc}(i)).
However, the relation $x\rla g$ is impossible
because $x\le_L u$ and $\supp(u)\subset\supp(g)$.
Hence either $u^{-1}g$ or $gu$ is not geodesic,
i.e.\ there exists $y_1\in V(\Gamma)^{\pm 1}$ such that
either $y_1^{-1}\le_R u^{-1}$ and $y_1\le_L g$
or $y_1\le_L u$ and $y_1^{-1}\le_R g$.
This means that $g^{y_1}$ is a cycling and that $u=y_1u_1$ is geodesic for some $u_1\in A(\Gamma)$.
Therefore $(g^{y_1})^{u_1}$ is a cyclic conjugation because
$\|(g^{y_1})^{u_1}\| =\|g^u\| =\|g\| =\|g^{y_1}\|$
and $\supp(u_1)\subset \supp(u)\subset\supp(g)=\supp(g^{y_1})$.
\end{proof} 

\smallskip(i)\ \
We may assume $u\neq 1$ because otherwise there is nothing to prove.

(a) $\Rightarrow$ (b):\ \
Suppose that $g^u$ is a cyclic conjugation.
By Claim 1, there is a geodesic decomposition $u=y_1u_1$
such that $g^{y_1}$ is a cycling and $(g^{y_1})^{u_1}$ is
a cyclic conjugation.
Applying Claim 1 again  to $(g^{y_1})^{u_1}$,
we have a geodesic decomposition $u_1=y_2u_2$ such that
$(g^{y_1})^{y_2}$ is a cycling
and $(g^{y_1y_2})^{u_2}$ is a cyclic conjugation.
Iterating this process, we get a desired reduced word
$w_0\equiv y_1\ldots y_k$.

(b) $\Rightarrow$ (c):\ \
Let $w\equiv x_1\cdots x_k$ be a reduced word representing $u$.
Notice that the word $w_0\equiv y_1\cdots y_k$ can be transformed
into the word $w\equiv x_1\cdots x_k$ by using only commutation relations.
Therefore each $(g^{x_1\cdots x_{i-1}})^{x_i}$ is a cycling
(by Lemma~\ref{lem:cy}(iii)).

(c) $\Rightarrow$ (a):\ \
Let $w\equiv x_1\cdots x_k$ be a reduced word representing $u$,
where $x_i=v_i^{\epsilon_i}$, $v_i\in V(\Gamma)$ and $\epsilon_i=\pm 1$
for all $1\le i\le k$.
Then, for each $1\le i\le k$,
$(g^{x_1\cdots x_{i-1}})^{x_i}$ is a cycling, hence
$$\|g^{x_1\cdots x_{i-1}}\|=\|g^{x_1\cdots x_i}\|\quad\mbox{and}\quad  v_i\in \supp(g^{x_1\cdots x_{i-1}})=\supp(g^{x_1\cdots x_i}).$$
Thus $\|g^{x_1\cdots x_k}\|=\|g\|$ and
$\{v_1,\ldots,v_k\}\subset \supp(g^{x_1\cdots x_{k}})= \supp(g)$.
Therefore $\|g^u\|=\|g\|$ and $\supp(u)\subset \supp(g)$,
hence  $g^u$ is a cyclic conjugation.

\smallskip(ii)\ \
Let $u_1=x_1\cdots x_j$ and $u_2=x_{j+1}\cdots x_k$ be geodesic decompositions,
where $x_1\cdots x_{k} \in V(\Gamma)^{\pm 1}$.
Then $u=x_1\cdots x_k$ is also geodesic because $u = u_1 u_2$ is geodesic.
By (i),
$g^u$ is a cyclic conjugation if and only if
$(g^{x_1\cdots x_{i-1}})^{x_i}$ is a cycling for each $1\le i\le k$,
and this happens if and only if
both $g^{u_1}$ and $(g^{u_1})^{u_2}$ are cyclic conjugations.

\smallskip(iii)\ \
Let $C(g)$ be the set of all $u\in A(\Gamma)$ such that $g^u$ is a
cyclic conjugation.
Then $C(g)$ satisfies (P1) in Corollary~\ref{cor:clo}
by (ii) in this lemma. Therefore it suffices to show
that $C(g)$ satisfies (P2) in Corollary~\ref{cor:clo}.

Let $ux, uy\in C(g)$ (i.e.\ $g^{ux}$ and $g^{uy}$ are cyclic conjugations)
such that $ux$ and $uy$ are geodesic and $x\rla y$,
where $u\in A(\Gamma)$ and $x,y\in V(\Gamma)^{\pm 1}$.
Then both $(g^u)^x$ and $(g^u)^y$ are cyclings of $g^u$ (by (ii)).
By Lemma~\ref{lem:cy}(iv), $(g^{ux})^y$ is a cycling,
hence $g^{uxy}$ is a cyclic conjugation (by (ii)).
Therefore $uxy\in C(g)$, hence
$C(g)$ satisfies (P2) in Corollary~\ref{cor:clo}.

\smallskip(iv)\ \
Notice that $u_1 \vee_L u_2 = u_1 u_2 = u_2 u_1$ and that both $u_1 u_2$ and $u_2 u_1$ are geodesic,
because $u_1 \rla u_2$.
Since both $g^{u_1}$ and $g^{u_2}$ are cyclic conjugations, so are
$g^{u_1 u_2}$, $(g^{u_1})^{u_2}$ and $(g^{u_2})^{u_1}$ (by (ii) and (iii)).

Let us show that the cyclic conjugation $(g^{u_2})^{u_1}$ is a left cyclic conjugation,
i.e.\ the decomposition $g^{u_2} u_1$ is geodesic.
(The proof for $(g^{u_1})^{u_2}$ is analogous.)
Observe
$$u_2^{-1}gu_1 = g^{u_2} u_1 u_2^{-1}.$$
Since both $u_2^{-1}g$ and $g u_1$ are geodesic and since $u_1 \rla u_2$,
$u_2^{-1}gu_1$ is geodesic (by Lemma~\ref{lem:cc}(ii)).
Since $\|g^{u_2}\| =\| g\|$, the decomposition $g^{u_2} u_1 u_2^{-1}$ is also geodesic.
Therefore $g^{u_2} u_1$ is geodesic.

\smallskip(v)\ \
We use induction on $\|u\|$.
If $\|u\|=1$, there is nothing to prove.

Suppose that $u=u'x$ is geodesic for some $u'\in A(\Gamma)\setminus\{1\}$
and $x\in V(\Gamma)^{\pm 1}$.
Then $g^{u'}$ is a cyclic conjugation and $(g^{u'})^x$ is a cycling (by (ii)).
Suppose that $(g^{u'})^x$ is a left cycling.
(The proof is analogous for the case where $(g^{u'})^x$ is a right cycling.)
By the induction hypothesis, we have a geodesic decomposition $u'=u_1'u_2'$ such that
$u_1'\rla u_2'$ and $g^{u_1'}$ (resp.\ $g^{u_2'}$)
is a left (resp.\ right) cyclic conjugation.

\smallskip\noindent
\textbf{Claim 2.}\ \
$x \rla u_2'$, and $u=u_1' x u_2'$ is geodesic.

\begin{proof}[Proof of Claim 2]
Let $u_2' = y_1\cdots y_k$ be geodesic, where $y_1,\ldots, y_k\in V(\Gamma)^{\pm 1}$.
Then $u$ has the following three geodesic decompositions:
$$u=u'x=u_1' u_2' x = u_1' y_1\cdots y_k x.$$

Let $h_0=g^{u_1'}$ and $h_i=g^{u_1'y_1\cdots y_i}$ for $1\le i\le k$.
Then each $h_i$ is cyclically reduced (by (ii)), and $h_i=h_{i-1}^{y_i}$.
Since $(g^{u_1'})^{u_2'}$ is a right cyclic conjugation (by (iv)),
each $h_{i-1}^{y_i}=(g^{u_1'y_1\cdots y_{i-1}})^{y_i}$ is a
right cycling (by (i)).

Since $u=u_1' y_1\cdots y_k x$ is geodesic, we have $y_k\neq x^{-1}$.
We know that $h_{k-1}^{y_k}$ is a right cycling and
that $(h_{k-1}^{y_k})^x=(g^{u'})^x$ is a left cycling,
hence $x\rla y_k$ (by Lemma~\ref{lem:cy}(ii)).
Therefore $u= u_1' y_1\cdots y_{k-1} x y_k$
and $h_{k-1}^x$ is a left cycling (by Lemma~\ref{lem:cy}(iii)).

Applying the above argument to the right cyclings
$h_{k-2}^{y_{k-1}},\ldots, h_0^{y_1}$
in this order iteratively,
we obtain $x\rla y_i$ for all $1\le i\le k$.
Therefore $x \rla u_2'$ and hence $u=u_1' u_2' x = u_1' x u_2'$.
Since $u_1' u_2' x$ is geodesic, so is $u_1' x u_2'$.
\end{proof} 

Let $u_1=u_1'x$ and $u_2=u_2'$.
Then $u=u_1u_2$ is geodesic, $u_1\rla u_2$ and
$g^{u_1}$ (resp.~$g^{u_2}$) is a left
(resp.~right) cyclic conjugation.
Moreover, $u_2^{-1}gu_1$ is geodesic (by (iv))
\end{myproof}

For a cyclically reduced $g\in A(\Gamma)$,
if $u\le_L g$, then $g^u$ is obviously a left cyclic conjugation.
The following proposition is concerned with the opposite direction.

\begin{proposition}\label{prop:cc1}
Let $g, u\in A(\Gamma)$ with $g$ cyclically reduced.
Then the following are equivalent.
\begin{enumerate}
\item $g^u$ is a left (resp.\ right) cyclic conjugation.
\item $u\le_L g^n$ (resp.\ $u^{-1}\le_R g^n$) for some $n\ge 1$.
\end{enumerate}
\end{proposition}

\begin{myproof}
We prove the equivalence only for the left cyclic conjugation.
The proof for the right cyclic conjugation is analogous.
We may assume $\|g\|\ge 2$ and $\|u\|\ge 1$
(otherwise it is obvious).

\smallskip\noindent
(ii) $\Rightarrow$ (i):
We may assume $n\ge 2$ (otherwise it is obvious).
We proceed by induction on $\|u\|$.
If $\|u\|=1$, then $u$ is a letter.
In this case, $u\le_L g^n$ implies  $u\le_L g$ (by Lemma~\ref{lem:cred1}(iii)),
hence $g^u$ is a left cycling.

Suppose $\|u\|\ge 2$.
Then $u=x u_1$ is geodesic
for some $x\in V(\Gamma)^{\pm 1}$ and $u_1\in A(\Gamma)\backslash\{1\}$.
Since $xu_1=u\le_L g^n$,
we get $x\le_L g^n$ and hence $x\le_L g$ (by Lemma~\ref{lem:cred1}(iii)).
Therefore $g=xg_1$ is geodesic for some $g_1\in A(\Gamma)$,
and $g^x=g_1x$ is also geodesic.
Since both $g\cdots g$ and $g=xg_1$ are geodesic,
$xg_1 xg_1 \cdots xg_1 x$ is geodesic,
hence the following three decompositions are all geodesic.
$$g^n x= xg_1 xg_1 \cdots xg_1 x= x(g^x)^n$$
Since $xu_1=u\le_L g^n\le_L g^nx=x(g^x)^n$,
we have $u_1\le_L (g^x)^n$ (by Lemma~\ref{lem:PO}(ii)).
By the induction hypothesis, $(g^{x})^{u_1}$ is a left cyclic conjugation.
And $g^{x}$ is also a left cyclic conjugation because $x\le_L g$.
Therefore $g^u$ is a left cyclic conjugation (by Lemma~\ref{lem:cy2}(ii)).

\smallskip\noindent
(i) $\Rightarrow$ (ii): As before, we use induction on $\|u\|$.
If $\|u\|=1$, $g^u$ is a left cycling, hence $u\le_L g$.

Suppose $\|u\|\ge 2$.
Then $u= xu_1$ is geodesic for some $x\in V(\Gamma)^{\pm}$
and $u_1\in A(\Gamma)\backslash\{1\}$.
Since $g^u$ is a left cyclic conjugation, both $g^x$ and
$(g^x)^{u_1}$ are left cyclic conjugations (by Lemma~\ref{lem:cy2}(ii)).
Since $g^x$ is a left cyclic conjugation and $x$ is a letter, we have $x\le_L g$.
Since $(g^x)^{u_1}$ is a left cyclic conjugation,
$u_1\le_L (g^x)^n$ for some $n$ by the induction hypothesis.
Using a similar argument as above,
we get that $x(g^x)^n$ and $g^n x$ are geodesic, hence
$u=xu_1\le_L x(g^x)^n=g^n x\le_L g^{n+1}$.
\end{myproof}

From the above proposition, $g^u$ is a right cyclic conjugation if and only if $\left( g^{-1}\right)^{u}$ is a left cyclic conjugation.

\begin{theorem}\label{thm:cnj1}
Let $g, u\in A(\Gamma)$ with $g$ cyclically reduced.
Then there exists a unique geodesic decomposition $u=u_1u_2u_3$ such that
\begin{itemize}
\item[(i)] $u_1$ disjointly commutes with $g$;
\item[(ii)] $g^{u_2}$ is a cyclic conjugation;
\item[(iii)] $g^u=u_3^{-1}g^{u_2}u_3$ is geodesic,
i.e.\ $\|g^u\|=\|u_3^{-1}\|+\|g^{u_2}\|+\|u_3\|=\|g\|+2\|u_3\|$.
\end{itemize}
Moreover, the following hold:
$u_1$ is the maximal prefix of $u$ that disjointly commutes with $g$;
$u_2$ is the maximal prefix of $u$ such that $g^{u_2}$ is a cyclic conjugation;
$u_3=g^u\wedge_R (g^u)^{-1}$.
In particular, $u_1\rla u_2$.
\end{theorem}

\begin{myproof}
We first prove the existence of the decomposition $u=u_1u_2u_3$.

If $u_1'$ and $u_1''$ are prefixes of $u$
such that $u_1'\rla g$ and $u_1''\rla g$,
then $u_1'\vee_L u_1''$ exists
(because $u_1'$ and $u_1''$ have a common right multiple $u$).
Observe that $u_1'\vee_L u_1''$ is also a prefix of $u$
and also disjointly commutes with $g$
(by Theorem~\ref{thm:lcm}).
Therefore there exists a unique maximal prefix  $u_1$  of $u$
that disjointly commutes with $g$.

If $u_2'$ and $u_2''$ are prefixes of $u$
such that $g^{u_2'}$ and $g^{u_2''}$
are cyclic conjugations,
then $u_2'\vee_L u_2''$ exists
(because $u_2'$ and $u_2''$ have a common right multiple $u$)
and is also a prefix of $u$,
and $g^{u_2'\vee_L u_2''}$ is also a cyclic conjugation
(by Lemma~\ref{lem:cy2}(iii)).
Therefore there exists a unique maximal prefix  $u_2$  of $u$
such that $g^{u_2}$ is a cyclic conjugation.

Notice that $u_1\rla u_2$
because $\supp(u_2)\subset \supp(g)$ and $u_1\rla g$.
Thus $u_1\vee_L u_2=u_1u_2$ is a prefix of $u$
and $u_1 u_2$ is geodesic,
hence $u=u_1u_2u_3$ is geodesic for some $u_3\in A(\Gamma)$.
Observe
$$
g^u
=u_3^{-1}u_2^{-1}u_1^{-1}gu_1u_2u_3
=u_3^{-1}u_2^{-1}gu_2u_3
=u_3^{-1} g^{u_2} u_3.$$

Let us show that $u_3^{-1}g^{u_2}u_3$ is geodesic.

If $g^{u_2} u_3$ is not geodesic, then
there exists $x\in V(\Gamma)^{\pm 1}$
such that $x\le_L u_3$ and $x^{-1}\le_R g^{u_2}$ (by Lemma~\ref{lem:geo}(i)),
hence $(g^{u_2})^x$ is a cyclic conjugation.
Notice that $u_2x$ is geodesic.
By Lemma~\ref{lem:cy2}(ii), $g^{u_2x}$ is also a cyclic conjugation,
hence $x\in\supp(g)$, which implies $x\rla u_1$.
Consequently, $u_2 x\rla u_1$ and hence $u_2 x\le_L u$.
This contradicts the maximality of $u_2$.
Therefore $g^{u_2} u_3$ is geodesic.
Similarly $u_3^{-1}g^{u_2}$ is geodesic.

Since both  $u_3^{-1}g^{u_2}$ and $g^{u_2} u_3$ are geodesic,
if $u_3^{-1}g^{u_2}u_3$ is not geodesic,
then there exists $x\in V(\Gamma)^{\pm 1}$ such that
$x\le_L u_3$, $x^{-1}\le_R u_3^{-1}$ and $x\rla g^{u_2}$
(by Lemma~\ref{lem:cc}(i)).
Since $\supp(g)=\supp(g^{u_2})$, we have $x\rla g$,
hence $x\rla u_2$ and $u_1x\rla g$.
Notice that $u_1x$ is geodesic.
Since $u_1x$ is a prefix of $u$, this contradicts the maximality of $u_1$.
Therefore $u_3^{-1}g^{u_2}u_3$ is geodesic.

Since $g^u=u_3^{-1}g^{u_2}u_3$ is geodesic
such that $g^{u_2}$ is cyclically reduced,
$u_3$ satisfies the formula
$u_3=g^u\wedge_R (g^u)^{-1}$ (by Lemma~\ref{lem:u-cu}(i)).

So far we have shown that $u=u_1u_2u_3$ is a desired decomposition.
We will now show the uniqueness of the decomposition.
Let $u=u_1'u_2'u_3'$ be another geodesic decomposition
satisfying the conditions (i), (ii) and (iii) of the theorem.

Since $u_2'$ and $u_3'$ satisfy the conditions (ii) and (iii),
we have $u_3'=g^u\wedge_R (g^u)^{-1}$ (by Lemma~\ref{lem:u-cu}(i)),
hence
$$u_3=u_3',\quad
g^{u_2}=g^{u_2'}\quad \mbox{and}\quad
u_1u_2=u_1'u_2'.
$$

Since both $u_1$ and $u_1'$ are prefixes of $u$ that disjointly commute with $g$,
so is $u_1\vee_L u_1'$ (by Theorem~\ref{thm:lcm}).
By the maximality of $u_1$, we have $u_1\vee_L u_1'\le_L u_1$, hence $u_1'\le_L u_1$.

Similarly, since $g^{u_2}$ and $g^{u_2'}$ are cyclic conjugations,
so is $g^{u_2\vee_L u_2'}$ (by Lemma~\ref{lem:cy2}).
By the maximality of $u_2$, we have $u_2\vee_L u_2'\le_L u_2$, hence $u_2'\le_L u_2$.

Since $u_1'\le_L u_1$, $u_2'\le_L u_2$ and $u_1u_2=u_1'u_2'$,
we have $u_1=u_1'$ and $u_2=u_2'$.
\end{myproof}

The following seems to be well known to experts.

\begin{corollary}\label{lem:cnj2}
Let $g_1, g_2\in A(\Gamma)$ be cyclically reduced.
If $g_1$ and $g_2$ are conjugate, then they are cyclically conjugate.
\end{corollary}

\begin{myproof}
Since $g_1$ and $g_2$ are conjugate,
$g_2=g_1^u$ for some $u\in A(\Gamma)$.
Let $u=u_1u_2u_3$ be the geodesic decomposition
for $g_1^u$ as in Theorem~\ref{thm:cnj1}.
Since $u_1\rla g_1$, we may assume $u_1=1$.
Since $u_3^{-1} g_1^{u_2}u_3$ is a geodesic decomposition
of $g_2$ and since $\|g_2\|=\|g_1\|=\|g_1^{u_2}\|$,
we have $u_3=1$.
Therefore $u=u_2$, hence $g_1$ is cyclically conjugate to $g_2$.
\end{myproof}

\EKLnewpage

\section{Star length}
\label{sec:star}

Star lengths of elements of $A(\Gamma)$,
introduced in~\cite{KK14},
induce a metric $d_*$ on $A(\Gamma)$ such that
the metric space $(A(\Gamma),d_*)$ is quasi-isometric
to the extension graph $(\Gamma^e,d)$,
preserving the right action of $A(\Gamma)$.
In this section, we study basic properties of star lengths.

It is known that
the centralizer $Z(v)$ of $v\in V(\Gamma)$ in $A(\Gamma)$ is generated by
the vertices in $\St_\Gamma(v)$.

\begin{definition}[star-word, star length]
An element in the centralizer $Z(v)$ of some vertex $v$
is called a \emph{star-word}.
The \emph{star length} of $g\in A(\Gamma)$,
denoted $\|g\|_*$,
is the minimum $\ell$ such that
$g$ is written as a product of $\ell$ star-words.
Let $d_*$ denote the right-invariant metric on $A(\Gamma)$
induced by the star length: $d_*(g_1,g_2)=\|g_1g_2^{-1}\|_*$.
\end{definition}

The following example illustrates that the decompositions into
star-words are not unique.

\begin{example}
Let $\Gamma = \bar P_5$, where $P_5= (v_1,\ldots,v_5)$ is a path graph,
and let the underlying right-angled Artin group here be $A(\Gamma)$,
hence $v_iv_j=v_jv_i$ whenever $|i-j|\ge 2$.
Let $g=v_1v_3v_5v_2v_4$.
The following shows various decompositions of $g$
into two star-words.
\begin{align*}
g& =(v_1v_3v_5)(v_2v_4)
=(v_1v_3v_5v_2)(v_4)
=(v_1v_3v_5v_4)(v_2)\\
&=(v_1v_3v_2)(v_5v_4)
=(v_3v_5v_4)(v_1v_2)
=(v_3v_5)(v_4v_1v_2).
\end{align*}
Notice that all the parenthesized words are star-words.
For example, $v_1v_3v_5\in Z(v_i)$ for $i=1,3,5$,
$v_4\in Z(v_i)$ for $i=1,2,4$,
$v_1v_3v_5v_2\in Z(v_5)$ and so on.
Since $\supp(g)=\{v_1,\ldots,v_5\}$ is not contained
in $\St(v_i)$ for any $1\le i\le 5$, we have $\|g\|_*=2$.
\end{example}

The group $A(\Gamma)$ acts on $(A(\Gamma),d_*)$ by right multiplication
$w\mapsto wg$.
Recall that $A(\Gamma)$ acts on $(\Gamma^e,d)$ by conjugation
$v^w\mapsto v^{wg}$.
For any $v\in V(\Gamma)$, the following map is equivariant.
$$
\phi_{v}: A(\Gamma)\to \Gamma^e,\qquad \phi_v(w)=v^{w}
$$

\begin{lemma}\cite[Lemma 19]{KK14}
\label{lem:KK}
Let $\Gamma$ be connected and let $D=\diam(\Gamma)$.
The following holds between the metric
$d$ on $\Gamma^e$
and the star length $\|\cdot\|_*$ on $A(\Gamma)$:
for any $g\in A(\Gamma)$ and $v\in V(\Gamma)$,
\begin{align*}
\|g\|_*-1\le d(v^g,v)\le D(\|g\|_*+1).
\end{align*}
\end{lemma}

Notice that
$d(\phi_v(g),\phi_v(h))=d(v^g,v^h)=d(v^{gh^{-1}},v)$
and $d_*(g,h)=\|gh^{-1}\|_*$.
Therefore the above lemma implies that
$d_*(g,h)-1\le d(v^g,v^h)\le D(d_*(g,h)+1)$,
and hence that $\phi_v$ is a quasi-isometry.
The above lemma also yields the following corollary
for the asymptotic translation length.

\begin{corollary}\label{cor:QI}
Let $\Gamma$ be connected and let $D=\diam(\Gamma)$.
For every $g\in A(\Gamma)$,
$$
\tau_{(A(\Gamma),d_*)}(g) \le
\tau_{(\Gamma^e,d)}(g) \le
D\tau_{(A(\Gamma),d_*)}(g).
$$
\end{corollary}

\begin{myproof}
Notice that
\begin{align*}
\tau_{(A(\Gamma),d_*)}(g)
&=\lim_{n\to\infty} \frac{d_*(g^n,1)}n
=\lim_{n\to\infty} \frac{\|g^n\|_*}n,\\
\tau_{(\Gamma^e,d)}(g)
&=\lim_{n\to\infty} \frac{d(v^{g^n},v)}n,
\end{align*}
where $v$ is any vertex of $\Gamma$.
By Lemma~\ref{lem:KK},
$$
\frac{\|g^n\|_*-1}n\le
\frac{d_*(v^{g^n},v)}n\le
\frac{D(\|g^n\|_*+1)}n.
$$
By taking $n$ to infinity, we get the desired inequalities.
\end{myproof}

The following lemma shows basic properties of star length.
\begin{lemma}\label{lem:st0}
Let $g_1,g_2,g_3, g, h\in A(\Gamma)$.
\begin{itemize}
\item[(i)] If $g_1g_2g_3$ is geodesic,
then $\|g_1g_3\|_*\le \|g_1g_2g_3\|_*$.
In particular, if $g\le_L h$ or $g\le_R h$, then $\|g\|_*\le \|h\|_*$.
\item[(ii)] $\|g^m\|_*\le \|g^n\|_*$ for all $1\le m\le n$.
\item[(iii)] If $g\rla h$ and $h\ne 1$, then $\|g\|_*\le 1$.
\end{itemize}
\end{lemma}

\begin{myproof}
Let us denote $g\preccurlyeq_0 h$ if a reduced word representing $g$
can be obtained by deleting some letters
from a reduced word representing $h$.
For example, if $v_i$'s are distinct vertices,
then $v_1v_3\preccurlyeq_0 v_1v_2v_3v_4$.
It is proved in \cite[Lemma 20(i)]{KK14} that if $g\preccurlyeq_0 h$,
then $\|g\|_*\le \|h\|_*$.

\smallskip\noindent
(i)\ \
Since $g_1g_2g_3$ is geodesic,
we have $g_1g_3\preccurlyeq_0 g_1g_2g_3$,
hence $\|g_1g_3\|_*\le \|g_1g_2g_3\|_*$.

\smallskip\noindent
(ii)\ \ Let $g=u^{-1}hu$ be geodesic such that $h$ is cyclically reduced.
Then $g^k=u^{-1}\underbrace{h\cdots h}_k u$ is also geodesic
for all $k\ge 1$ (by Lemma~\ref{lem:cc}(iii)).
Therefore $g^m\preccurlyeq_0 g^n$,
hence $\|g^m\|_*\le \|g^n\|_*$.

\smallskip\noindent
(iii)\ \
Since $h\ne 1$, there is a vertex $v\in\supp(h)$.
Then $g\in Z(v)$, namely $\|g\|_*\le 1$.
\end{myproof}

\begin{lemma}\label{lem:st1}
Suppose that $g_1,g_2\in A(\Gamma)$ have a common right multiple
and that none of them is a prefix of the other, i.e.\
$g_1\not\le_L g_2$ and $g_2\not\le_L g_1$.
Then
$\|g_1^{-1}g_2 \|_*\le 2$
and $\|g_1\|_*-\|g_2\|_*\in\{0,\pm 1\}$.
\end{lemma}

\begin{myproof}
Let $g_i=(g_1\wedge_L g_2) g_i'$ for $i=1,2$.
Since $g_1$ and $g_2$ have a common right multiple,
$g_1'\rla g_2'$ (by Theorem~\ref{thm:lcm}).
Since $g_1\not\le_L g_2$ and $g_2\not\le_L g_1$,
both $g_1'$ and $g_2'$ are nontrivial,
hence $\|g_1'\|_*=\|g_2'\|_*=1$ (by Lemma~\ref{lem:st0}(iii)).
Therefore
$$
\|g_1^{-1}g_2\|_* =\|g_1'^{-1}g_2'\|_*
\le \|g_1'\|_* + \|g_2'\|_* = 1+1=2.
$$
Furthermore, for each $i=1,2$,
$$\|g_1\wedge_L g_2\|_*\le \|g_i\|_* \le \|g_1\wedge_L g_2 \|_* + \|g_i'\|_*
= \|g_1\wedge_L g_2\|_*+1,$$
hence $ \|g_i\|_*=\|g_1\wedge_L g_2\|_* + \epsilon_i$, where $\epsilon_i\in\{0,1\}$.
Therefore
$\|g_1\|_*-\|g_2\|_* = \epsilon_1-\epsilon_2\in\{0,\pm 1\}$.
\end{myproof}

\begin{corollary}\label{cor:st2}
Let $g_1, g_2, h\in A(\Gamma)$ with $g_1g_2$ geodesic.
If $h\le_L g_1g_2$ and $\|g_1\|_*\ge \|h\|_*+2$,
then $h\le_L g_1$.
\end{corollary}

\begin{myproof}
Observe that $g_1\not\le_L h$
(otherwise $\|g_1\|_*\le \|h\|_*$).
Assume $h\not\le_L g_1$.
Since $g_1$ and $h$ have a common right multiple, say $g_1g_2$,
we have $\|g_1\|_*-\|h\|_* \in\{0,\pm 1\}$
(by Lemma~\ref{lem:st1}).
This contradicts that $\|g_1\|_*\ge \|h\|_*+2$.
\end{myproof}

\begin{corollary}\label{cor:st3}
Let $g_1,g_2\in A(\Gamma)$.
If\/ $g_1g_2$ is geodesic, then
$$
\|g_1\|_* + \|g_2\|_* -2 \le
\|g_1g_2\|_*\le \|g_1\|_* + \|g_2\|_* .
$$
\end{corollary}

\begin{myproof}
Let $r=\|g_1\|_*$, $s=\|g_2\|_*$
and $t=\|g_1g_2\|_*$.
Then it is obvious that $t\le r+s$, hence it suffices to show $t\ge r+s-2$.
Since $g_2\le_R g_1g_2$, we have $t=\|g_1g_2\|_*\ge \|g_2\|_*=s$
(by Lemma~\ref{lem:st0}).
We may assume $r\ge 3$ because otherwise $t\ge s\ge r+s-2$.
Let
$$
g_1g_2=w_1w_2\cdots w_t
$$
be a geodesic decomposition of $g_1g_2$ into star-words.
Then $w_1\cdots w_{r-2}\le_L g_1$ (by Corollary~\ref{cor:st2}),
hence $g_2\le_R w_{r-1}\cdots w_t$.
Therefore $s=\|g_2\|_*\le
\|w_{r-1}\cdots w_t\|_*=t-r+2$,
namely $t\ge r+s-2$.
\end{myproof}

The following example shows that
the upper and lower bounds in the above corollary are sharp.

\begin{example}
Let $\Gamma = \bar P_5$, where $P_5 = (v_1,\ldots,v_5)$,
and let the underlying right-angled Artin group here be $A(\Gamma)$.

(i)
Let $g_1 = v_1 v_2$ and $g_2 = v_3 v_4$.
Then $g_1 \rla v_5$ and $g_2 \rla v_1$ and hence $\|g_1\|_* = \|g_2\|_* =1$.
Since $g_1g_2=v_1v_2v_3v_4\not\in Z(v_i)$ for any $1\le i\le 5$,
we have $\|g_1g_2 \|_* \ge 2$.
Since $\|g_1g_2\|_* \le \|g_1\|_* + \|g_2\|_*=2$,
we have $\|g_1g_2\|_* =\|g_1\|_* + \|g_2\|_*$
in this case.

(ii)
Let $g_1=g_2=v_2v_3v_4$.
Then $g_1g_2=v_2v_3v_4\cdot v_2v_3v_4=v_2v_3v_2\cdot v_4v_3v_4$.
Since $v_2v_3v_2\in Z(v_5)$ and $v_4v_3v_4\in Z(v_1)$,
we have $\|v_2v_3v_2\|_*=\|v_4v_3v_4\|_*=1$.
It is easy to see that $\|g_1\|_*=\|g_2\|_*=\|g_1g_2\|_*=2$.
Therefore $\|g_1g_2\|_* = \|g_1\|_* + \|g_2\|_*-2$ in this case.
\end{example}

The following is an immediate consequence of Lemma~\ref{lem:st0}(ii)
and Corollary~\ref{cor:st3}.

\begin{corollary}\label{cor:st4}
Let $g\in A(\Gamma)$ be cyclically reduced.
Then  $\{\|g^n\|_* \}_{n=0}^\infty$ is an increasing sequence
such that the following hold.
\begin{enumerate}
\item
If $\|g\|_* = 1$, then $\|g^n\|_* = 1$ for all $n\ge 1$.

\item
If $\|g\|_* = 2$, then $\|g^{n-1}\|_* \le \|g^{n}\|_* \le \|g^{n-1}\|_* +2$
for all $n\ge 1$.

\item
If $\|g\|_* \ge 3$, then $\|g^n\|_* \ge \|g^{n-1}\|_* +1$
and hence $\|g^n\|_* \ge n+2$
for all $n\ge 1$.
\end{enumerate}
\end{corollary}

\begin{corollary}\label{cor:st5}
Let $g, u\in A(\Gamma)$ with $g$ cyclically reduced.
If $\|g\|_*\ge 3$ and
$g\not\le_L u\le_L g^n$ for some $n\ge 1$,
then $u\le_L g^2$.
\end{corollary}

\begin{myproof}
We may assume $n\ge 3$ and $u\not\le_L g$
(otherwise it is obvious).
Since $g$ and $u$ have a common right multiple, say $g^n$, there exist
$g'$ and $u'$ such that
$gg'=uu'=g\vee_L u\le_L g^n$ and $u'\rla g'$
(by Theorem~\ref{thm:lcm}),
where $gg'$ and $uu'$ are geodesic.
Since $gg'\le_L g^{n}=g g^{n-1}$, (by Lemma~\ref{lem:PO})
$$
g'\le_L g^{n-1}=g\cdot g^{n-2}.
$$
Since $u\not\le_L g$ and $g\not\le_L u$,
both $g'$ and $u'$ are nontrivial,
hence $\|g'\|_*=\|u'\|_*=1$.
Since $\|g'\|_*=1$ and $\|g\|_*\ge 3$, we get
$g'\le_L g$ (by Corollary~\ref{cor:st2}).
Therefore $u\le_L uu'=gg'\le_L g^2$.
\end{myproof}

\begin{lemma}\label{lem:st5}
Let $g_1,g_2,g_3\in A(\Gamma)$ be such that
both $g_1g_2$ and $g_2g_3$ are geodesic.
If\/ $\|g_2\|_*\ge 2$,
then $g_1g_2g_3$ is geodesic.
\end{lemma}

\begin{myproof}
Assume that $g_1g_2g_3$ is not geodesic.
Since $g_1g_2$ and $g_2g_3$ are geodesic, there exists $x\in V(\Gamma)^{\pm1}$ such that
$x^{-1}\le_R g_1$, $x\le_L g_3$ and $x\rla g_2$ (by Lemma~\ref{lem:cc}(i)).
Observe that $x\rla g_2$ implies $\|g_2\|_*\le 1$,
which contradicts the hypothesis $\|g_2\|_*\ge 2$.
\end{myproof}

We introduce the notion of strongly non-split elements.
We will see (in Lemma~\ref{lem:loxo} and Remark~\ref{rmk:loxo})
that if $|V(\Gamma)|\ge 4$ and both $\Gamma$ and $\bar\Gamma$ are connected,
then a cyclically reduced element $g\in A(\Gamma)$ is strongly non-split
if and only if $g$ is loxodromic on the extension graph $\Gamma^e$.

\begin{definition}[non-split, strongly non-split]
Let $g\in A(\Gamma)\setminus\{ 1\}$.
\begin{enumerate}
\item
$g$ is called \emph{split}
if $g$ has a nontrivial geodesic decomposition $g=g_1 g_2$ with $g_1 \rla g_2$.

\item
$g$ is called \emph{non-split} if it is not split.

\item
$g$ is called \emph{strongly non-split} if $g$ is non-split and
$g\not\rla v$ for any $v\in V(\Gamma)$.
\end{enumerate}
\end{definition}

\begin{table}
\begin{tabular}{lll}\hline
\qquad $\rule{0pt}{12pt} g\in A(\Gamma)$
& \qquad\qquad $\Gamma$
& \qquad\qquad $\bar\Gamma$\\\hline
$\rule{0pt}{12pt} g$ is split
& $\Gamma[g]$ is a join
& $\bar\Gamma[g]$ is disconnected\\
$\rule{0pt}{12pt}g$ is non-split
& $\Gamma[g]$ is not a join
& $\bar\Gamma[g]$ is connected\\
$\rule{0pt}{12pt}g$ is strongly non-split \mbox{}\quad\mbox{}
& $\Gamma[g]$ is not contained in \mbox{}\quad\mbox{}
& $\bar\Gamma[g]$ is connected and\\
& a subjoin of $\Gamma$
& $\St_{\bar\Gamma}(\supp(g))= V(\Gamma)$
\\\hline
\end{tabular}

\medskip
\caption{Equivalent conditions for $g\in A(\Gamma)$
to be split, non-split and strongly non-split}
\label{tab:loxo}
\end{table}

It is easy to see that $g\in A(\Gamma)$ is split
if and only if $\Gamma[g]$ is a join
(equivalently, $\bar\Gamma[g]$ is disconnected).
Similarly, one can characterize
the property of being non-split and strongly non-split
using the graphs $\Gamma[g]$ and $\bar\Gamma[g]$
as shown in Table~\ref{tab:loxo}.

From definition, the existence of a strongly non-split element implies that
$\bar\Gamma$ is connected.

\begin{remark}\label{rmk:sp-po} 
Let $n\ge 2$ and $g,h\in A(\Gamma)\setminus\{1\}$.
Observe that strongly non-splitness of an element depends only on its support.
Note that $\supp(g^{-1})=\supp(g)=\supp(g^n)$ (by Lemma~\ref{lem:cred1}),
and that if either $g\le_L h$ or $g\le_R h$,
then $\supp(g)\subset\supp(h)$.
Therefore

\begin{enumerate}
\item $g$ is strongly non-split if and only if $g^{-1}$ is strongly non-split;
\item $g$ is strongly non-split if and only if $g^n$ is strongly non-split;
\item if $g$ is strongly non-split and either $g\le_L h$ or $g\le_R h$,
then $h$ is also strongly non-split.
\end{enumerate}
\end{remark}

Strongly non-splitness is related to the star length as follows.

\begin{lemma}\label{lem:SNS}
Let $g\in A(\Gamma)\setminus\{1\}$.
\begin{enumerate}
\item
If\/ $\|g\|_*\ge 3$, then $g$ is strongly non-split.

\item
$g$ is strongly non-split with $| \supp(g) | \ge 2$ if and only if
$g$ is non-split with $\|g\|_*\ge 2$.
\end{enumerate}
\end{lemma}

\begin{myproof}
(i)\ \
Assume that $g$ is not strongly non-split.
If $g$ is split, then clearly $\|g\|_* \le 2$.
If $g$ is non-split but not strongly non-split, then there is $v\in V(\Gamma)\backslash\supp(g)$
with $v\rla g$, hence $\|g\|_* = 1$.
In either case, $\|g\|_* \le 2$.

(ii)\ \
Suppose that $g$ is strongly non-split with $| \supp(g) | \ge 2$.
Then $g$ is non-split by definition.
Assume $\|g\|_*= 1$.
Then there exists $v\in V(\Gamma)$ with $\supp(g)\subset Z(v)$.
Since $g$ is strongly non-split, $v\in\supp(g)$.
Since $| \supp(g) | \ge 2$ and $\supp(g)\subset Z(v)$,
$g=v^n g_1$ is geodesic for some $n\neq 0$ and $g_1\in A(\Gamma)\backslash\{ 1\}$
with $g_1\rla v$.
Namely, $g$ is split, which is a contradiction.
Therefore $\|g\|_*\ge 2$.

Conversely, suppose that $g$ is non-split with $\|g\|_*\ge 2$.
Then $| \supp(g) | \ge 2$ and there does not exit $v\in V(\Gamma)\backslash\supp(g)$
with $v\rla g$.
Therefore $g$ is strongly non-split.
\end{myproof}

\EKLnewpage

\section{Prefixes of powers of cyclically reduced elements}
\label{sec:power}

In this section, we study prefixes of powers of cyclically reduced elements.
The main result is Theorem~\ref{thm:po3},
which plays important roles in the study of
the asymptotic translation length
and the acylindricity of the action of $A(\Gamma)$ on $\Gamma^e$.

\begin{lemma}\label{lem:po1}
Let $u, g_1,g_2,\ldots, g_m\in A(\Gamma)$.
If $g_1g_2\cdots g_m$ is geodesic, then
for each $1\le k\le m$ there exists a geodesic decomposition
$g_k=a_kb_k$ such that
\begin{enumerate}
\item
$u\wedge_L (g_1\cdots g_k)=a_1\cdots a_k$;

\item
$a_k \rla b_j$ for all $1\le j\le k-1$;

\item
$a_1\cdots a_kb_1\dots b_k$ is a geodesic decomposition of $g_1\cdots g_k$.
\end{enumerate}
\end{lemma}

\begin{myproof}
The relation $u\wedge_L (g_1\cdots g_k)=a_1\cdots a_k$
determines the elements $a_k$ inductively for $k=1,\ldots,m$.
Then the relation $g_k=a_kb_k$ determines the elements $b_k$
for all $1\le k\le m$.
Therefore we get elements
$a_1,\ldots,a_m, b_1,\ldots,b_m$
such that
$u\wedge_L (g_1\cdots g_k)=a_1\cdots a_k$ and
$g_k=a_kb_k$ for all $1\le k\le m$.

Since $g_1\cdots g_m$ is geodesic,
$g_1\cdots g_k\le_L g_1\cdots g_{k+1}$ for each $1\le k\le m-1$
(by Lemma~\ref{lem:PO}), hence
$$
a_1\cdots a_k=u\wedge_L(g_1\cdots g_k)
\le_L u\wedge_L(g_1\cdots g_{k+1})=a_1\cdots a_{k+1}.
$$
Therefore $a_1\cdots a_m$ and hence each $a_1\cdots a_k$
are geodesic (by Lemma~\ref{lem:PO} again).

For each $1\le k\le m$, let $u_k\in A(\Gamma)$ be the element such that
$u=a_1\cdots a_k u_k$.
Then each $a_1\cdots a_k u_k$ is geodesic
because $a_1\cdots a_k\le_L u$.

\medskip\noindent
\textbf{Claim.}\ \ For each $1\le k\le m$,
\begin{itemize}
\item[(a)] $a_k\le_L g_k$, hence $g_k=a_kb_k$ is geodesic;
\item[(b)] $a_k\rla b_j$ for all $1\le j\le k-1$;
\item[(c)] $a_1\cdots a_k b_1\cdots b_k$ is a geodesic
decomposition of $g_1\cdots g_k$.
\end{itemize}

\begin{proof}[Proof of Claim]
We use induction on $k$.

For $k=1$, (a) and (c) hold because $a_1=u\wedge_L g_1 \le_L g_1$ and $g_1=a_1b_1$,
and (b) is vacuously true.

Assume that the claim holds for some $1\le k<m$.
We now have the following geodesic decompositions at hand:
\begin{align*}
u                   &=(a_1\cdots a_k)u_k,\\
g_1\cdots g_k       &=(a_1\cdots a_k)(b_1\dots b_k),\\
g_1\cdots g_{k+1}   &=(a_1\cdots a_k)(b_1\cdots b_k) g_{k+1}.
\end{align*}

Since $u\wedge_L (g_1\cdots g_k)=a_1\cdots a_k$, we have
$u_k\wedge_L (b_1\cdots b_k)=1$ (by Lemma~\ref{lem:gcd1}).

Since $u\wedge_L (g_1\cdots g_{k+1})=a_1\cdots a_{k+1}$,
we have $a_{k+1}=u_k\wedge_L (b_1\cdots b_k g_{k+1})$, hence
$$
a_{k+1}\le_L u_k\quad\mbox{and}\quad
a_{k+1}\le_L (b_1\cdots b_k) g_{k+1}.
$$

Since $a_{k+1}\le_L u_k$, we have
$a_{k+1}\wedge_L (b_1\cdots b_k)\le_L
u_k\wedge_L (b_1\cdots b_k)=1$ (by Lemma~\ref{lem:gcd1}).

Since $a_{k+1}\le_L (b_1\cdots b_k) g_{k+1}$
and $a_{k+1}\wedge_L (b_1\cdots b_k)=1$,
we have
$$
a_{k+1}\rla b_1\cdots b_k
\quad\mbox{and}\quad
a_{k+1}\le_L g_{k+1}
$$
(by Lemma~\ref{lem:lcm1}(ii)).
In particular, $a_{k+1}\rla b_j$ for all $1\le j\le k$.
Therefore (a) and (b) hold for $k+1$.

Since $a_{k+1}\rla b_j$ for all $1\le j\le k$, we have
\begin{align*}
g_1 \cdots g_kg_{k+1}
&=(a_1\cdots a_k)(b_1\cdots b_k) (a_{k+1}b_{k+1})\\
&=(a_1\cdots a_{k+1})(b_1\cdots b_{k+1}).
\end{align*}
The above three decompositions are all geodesic
because $g_1\cdots g_k g_{k+1}$, $g_{k+1}=a_{k+1}b_{k+1}$
and $g_1\cdots g_k=a_1\cdots a_k b_1\cdots b_k$ are all geodesic.
Therefore (c) holds for $k+1$.
\end{proof}
The above claim completes the proof.
\end{myproof}

In the following, we frequently use the notation $\St_{\bar\Gamma[g]}(X)$,
for $X\subset\supp(g)$, which  denotes the star of $X$ in
$\bar\Gamma[g]=\bar\Gamma[\supp(g)]$.
Hence, $v\in\St_{\bar\Gamma[g]}(X)$
if and only if either $v\in X$
or $v\in\supp(g)$ and $\{v,v_1\}$ is an edge in $\bar\Gamma$
for some $v_1\in X$.
Therefore
$$\St_{\bar\Gamma[g]}(X)=\St_{\bar\Gamma}(X)\cap\supp(g).$$

When $g_1=\cdots=g_m$ in Lemma~\ref{lem:po1}, we have the following.

\begin{corollary}\label{cor:po2}
Let $m\ge 1$ and $g,u\in A(\Gamma)$ with $g$ cyclically reduced.
Then for each $1\le k\le m$ there exists a geodesic decomposition
$g=a_kb_k$ such that
\begin{itemize}
\item[(i)]
$u\wedge_L g^k=a_1a_2\cdots a_k$,
$a_k \rla b_j$ for all $1\le j<k$ and
$a_1\cdots a_kb_1\dots b_k$ is a geodesic decomposition of $g^k$;

\item[(ii)]
$\{a_k\}_{k=1}^m$ is descending with respect to $\le_L$ such that
\begin{eqnarray*}
& 1\le_L a_m\le_L \cdots \le_L a_2\le_L a_1\le_L g,\\
&\St_{\bar\Gamma[g]}(\supp(a_{k+1}))  \subset \supp(a_k);
\end{eqnarray*}

\item[(iii)]
$\{b_k\}_{k=1}^m$ is ascending with respect to $\le_R$ such that
\begin{eqnarray*}
&1\le_R b_1\le_R b_2\le_R\cdots \le_R  b_m\le_R g,\\
&\St_{\bar\Gamma[g]}(\supp(b_k))  \subset \supp(b_{k+1}).
\end{eqnarray*}
\end{itemize}
In (ii) and (iii), we let $a_{m+1}=1$ and $b_{m+1}=g$ for notational convenience.
\end{corollary}

\begin{myproof}
Since $g$ is cyclically reduced, $gg\cdots g$ is geodesic.
By Lemma~\ref{lem:po1}, there exists a geodesic decomposition
$g=a_kb_k$ for $1\le k\le m$ satisfying (i).

Since $g=a_kb_k=a_{k+1}b_{k+1}$,
we have $a_{k+1}\le_L a_kb_k$.
Since $a_{k+1}\rla b_k$,
we have $a_{k+1}\le_L a_k$ (by Lemma~\ref{lem:gcd2}(iv)),
hence $\{a_k\}_{k=1}^m$ is descending with respect to $\le_L$.
Since $a_1\le_L g$ and $1\le_L a_m$,
$$
1\le_L a_m\le_L \cdots \le_L a_2\le_L a_1\le_L g.
$$
Since $g=a_kb_k$, it follows immediately from the above inequalities
that the sequence $\{b_k\}_{k=1}^m$ is ascending with respect to $\le_R$
such that
$$1\le_R b_1\le_R b_2\le_R\cdots \le_R  b_m\le_R g.$$

Since $\supp(g)=\supp(a_j)\cup\supp(b_j)$ for $j=k,k+1$,
\begin{align*}
\supp(g)-\supp(a_{k+1}) &\subset \supp(b_{k+1}),\\
\supp(g)-\supp(b_k)     &\subset \supp(a_k).
\end{align*}

Since $a_{k+1}\rla b_k$, by Lemma~\ref{lem:st-s1}
\begin{align*}
&\supp(b_k)\cap \St_{\bar\Gamma}(\supp(a_{k+1}))=\emptyset,\\
&\supp(a_{k+1})\cap \St_{\bar\Gamma}(\supp(b_{k}))=\emptyset.
\end{align*}

Hence
\begin{align*}
\St_{\bar\Gamma[g]}(\supp(a_{k+1}))
&=\supp(g)\cap \St_{\bar\Gamma}(\supp(a_{k+1})) \\
&\subset \supp(g)-\supp(b_k)\subset \supp(a_k),\\
\St_{\bar\Gamma[g]}(\supp(b_k))
&=\supp(g)\cap \St_{\bar\Gamma}(\supp(b_k)) \\
&\subset \supp(g)-\supp(a_{k+1})\subset \supp(b_{k+1}).
\end{align*}
Therefore (ii) and (iii) are proved.
\end{myproof}


\begin{theorem}\label{thm:po3}
Let $m\ge 2$ and $g, u\in A(\Gamma)$ with $g$ cyclically reduced and non-split.
If
$$
g\not\le_L u\not\le_L g^{m-1}   \quad\mbox{and}\quad
u\le_L g^m,
$$
then the following hold.
\begin{enumerate}
\item
$m\le\diam(\bar\Gamma[g])$.
In particular, $m\le | \supp(g) |-1\le | V(\Gamma) |-1$.

\item
There is a geodesic decomposition
$g=g_mg_{m-1}\cdots g_1g_0$
such that
\begin{enumerate}
\item
$g_k\ne 1$ for all $0\le k\le m$;

\item
$g_i\rla g_j$ whenever $|i-j|\ge 2$;

\item
$u\wedge_L g^k=(g_m\cdots g_1) (g_m\cdots g_2)\cdots (g_m \cdots g_k)$
for all $1\le k\le m$.
\end{enumerate}
In particular,
$u=u\wedge_L g^m=(g_m\cdots g_1)
(g_m\cdots g_2)\cdots (g_m g_{m-1})(g_m).$

\item
$\|u\|_*\le \|g\|_*+1$.

\item
If\/ $\|g\|_*\ge 3$, then $m=2$.
(Equivalently, if $m\ge 3$, then $\|g\|_*\le 2$.)
\end{enumerate}
\end{theorem}

\begin{myproof}
For each $k\ge1$, let $g=a_kb_k$ be the geodesic decomposition
given by Corollary~\ref{cor:po2}.
Then
\begin{itemize}
\item
$u\wedge_L g^k=a_1\cdots a_k$ is geodesic and
$a_k\rla b_j$ for all $1\le j< k$;
\item
$\{ a_k\}_{k=1}^{\infty}$ is descending with respect to $\le_L$ such that
$1\le_L \cdots \le_L  a_2\le_L a_1\le_L g$;
\item
$\{ b_k\}_{k=1}^{\infty}$ is ascending with respect to $\le_R$ such that
$1\le_R b_1\le_R b_2\le_R\cdots\le_R g$.
\end{itemize}

The following claim is a result of the hypothesis
that $g\not\le_L u\not\le_L g^{m-1}$ and
$u\le_L g^m$.

\medskip\noindent
\textbf{Claim 1.}\ \
For each $1\le k\le m$, $a_k \not\in\{ 1, g, a_{k+1}\}$ and hence $b_k \not\in\{ 1, g, b_{k+1}\}$.
For each $k> m$, $a_{k}=1$ and hence $b_{k}=g$.

\begin{proof}[Proof of Claim 1]
For each $k\ge 1$,
$a_1\cdots a_k\le_L u$ and $a_1\cdots a_k\le_L g^k$
(because $a_1\cdots a_k=u\wedge_L g^k$).
Furthermore, $u=u\wedge_L g^m=a_1\cdots a_m$ (from the hypothesis $u\le_L g^m$).
Therefore
$$
a_1\le_L u,\quad
a_1\cdots a_{m-1}\le_L g^{m-1},\quad
a_1\cdots a_m=u.
$$

If $a_1=g$, then $g\le_L u$, which contradicts the hypothesis $g\not\le_L u$.
If $a_m=1$, then $u=a_1\cdots a_{m-1}a_m=a_1\cdots a_{m-1}\le_L g^{m-1}$,
which contradicts the hypothesis $u\not\le_L g^{m-1}$.
Thus $a_1\ne g$ and $a_m\ne 1$.
Therefore, for each $1\le k\le m$, we get $a_k \not\in\{ 1, g\}$ (because $1\le_L a_m\le_L a_k\le_L a_1 \le_L g$)
and hence $b_k\not\in\{ 1, g\}$ (because $g=a_kb_k$).

Assume that $a_k=a_{k+1}$ for some $1\le k\le m$.
Then $a_k\rla b_k$ because $a_{k+1}\rla b_k$.
Since $g=a_kb_k$ and both $a_k$ and $b_k$ are nontrivial,
this contradicts that $g$ is non-split.
Therefore $a_k\ne a_{k+1}$
and hence $b_k\ne b_{k+1}$ for all $1\le k\le m$.

Let $j\ge 1$.
Since $u\le_L g^m$ and $g$ is cyclically reduced,
we have $u\le_L g^{m+j}$,
hence $u\wedge_L g^m=u=u\wedge_L g^{m+j}$.
Therefore $a_1\cdots a_m=a_1\cdots a_ma_{m+1}\cdots a_{m+j}$,
hence $a_{m+1}\cdots a_{m+j}=1$.
Since the decomposition $a_{m+1}\cdots a_{m+j}$ is geodesic,
we have $a_{m+j}=1$.
Namely, for all $k> m$, $a_k=1$ and hence $b_k=g$.
\end{proof}

Define $\{g_k\}_{k=0}^m$ by
$g_0=b_1$ and $g_k=a_{k+1}^{-1}a_k$ (hence $a_k=a_{k+1}g_k$) for $1\le k\le m$.
Then $a_k=a_{k+1}g_k$ is geodesic for all $1\le k\le m$
because $a_{k+1}\le_L a_k$.

\medskip\noindent
\textbf{Claim 2.}\ \
The decomposition $g=g_mg_{m-1}\cdots g_1g_0$ is geodesic
such that
\begin{itemize}
\item[(a)] $g_k\ne 1$ for all $0\le k\le m$;
\item[(b)] $g_i\rla g_j$ whenever $|i-j|\ge 2$;
\item[(c)] $u\wedge_L g^k=(g_m\cdots g_1) (g_m\cdots g_2)\cdots (g_m \cdots g_k)$
for all $1\le k\le m$;
\item[(d)] $a_k=g_mg_{m-1}\cdots g_k$ and $b_k=g_{k-1}g_{k-2}\cdots g_0$
for all $1\le k\le m$.
\end{itemize}

\begin{proof}[Proof of Claim 2]
Since $a_k=a_{k+1}g_k$ is geodesic,
$\|g_k\|=\|a_{k}\|-\|a_{k+1}\|$ for all $1\le k\le m$.
Since $g_0=b_1$, $a_{m+1}=1$ and $g=a_1 b_1$ is geodesic,
\begin{align*}
\|g_0\|& + \|g_1\|+\cdots + \|g_m\|\\
& = \|b_1\|+ (\|a_1\|-\|a_{2}\|)
 + \cdots + (\|a_m\|-\|a_{m+1}\|) \\
& = \|b_1\|+ \|a_1\|-\|a_{m+1}\| = \|g\|.
\end{align*}
Consequently, $\| g\| = \|g_0\| + \|g_1\|+\cdots + \|g_m\|$.

For $1\le k\le m$,
$a_k=a_{k+1}g_k=a_{k+2}g_{k+1}g_k
=\cdots = a_{m+1}g_m\cdots g_k
= g_m\cdots g_k$
because $a_{m+1}=1$.
Therefore we have the following decompositions:
\begin{align*}
a_k &= g_m\cdots g_k\quad\text{for all}\ 1\le k\le m,\\
g &=a_1b_1=(g_m\cdots g_1)g_0=g_m\cdots g_0,\\
b_k &= a_k^{-1}g=g_{k-1}\cdots g_0\quad\text{for all}\  1\le k\le m.
\end{align*}

Observe that $g =g_m\cdots g_0$ is geodesic because $\| g\| = \|g_0\| + \|g_1\|+\cdots + \|g_m\|$.

The decompositions for $a_k$ and $b_k$ in the above  prove (d).

For each $1\le k\le m$,
$u\wedge_L g^k=a_1a_2\cdots a_k
=(g_m\cdots g_1)(g_m\cdots g_2)\cdots (g_m\cdots g_k)$.
This proves (c).

By Claim 1, $g_0=b_1\ne 1$ and $g_k=a_{k+1}^{-1}a_k \ne 1$ for all $1\le k\le m$.
This proves (a).

For each $(i,j)$ with $0\le j<j+2\le i\le m$, we know that $a_i\rla b_{j+1}$.
Since $a_i=g_m\cdots g_i$ and $b_{j+1}=g_j\cdots g_0$ are geodesic,
we have $g_i\le_R a_i$ and $g_j\le_L b_{j+1}$,
hence $g_i\rla g_j$.
This proves (b).
\end{proof}

Recall from Claim 2 that both $g_0=b_1$ and $g_m=a_m$ are nontrivial.

\medskip\noindent
\textbf{Claim 3.}\ \
For any path $(v_0, v_1,\ldots, v_{r-1}, v_r)$ in  $\bar\Gamma[g]$
such that $v_0 \in \supp(g_0)=\supp(b_1)$
and $v_r\in\supp(g_m)=\supp(a_m)$,
we have $m\le r$. In particular, $m\le\diam(\bar\Gamma[g])$.

\begin{proof}[Proof of Claim 3]
Using induction on $k$, we first show that
$$
v_k\in\supp(b_{k+1})
$$
for all $0\le k\le \min\{m-1,r-1\}$.
By the hypothesis of the claim, $v_0\in\supp(b_1)$.
Assume that $v_k\in\supp(b_{k+1})$ for some $0\le k\le \min\{m-2,r-2\}$.
Since $\{v_k,v_{k+1}\}$ is an edge in $\bar\Gamma[g]$,
we have $v_{k+1}\in \St_{\bar\Gamma[g]}(v_k)$.
Since $v_k\in\supp(b_{k+1})$ by induction hypothesis,
$\St_{\bar\Gamma[g]}(v_k)\subset \St_{\bar\Gamma[g]}(\supp(b_{k+1}))$,
hence $v_{k+1}\in\St_{\bar\Gamma[g]}(\supp(b_{k+1}))$.
By Corollary~\ref{cor:po2},
$\St_{\bar\Gamma[g]}(\supp(b_{k+1})) \subset \supp(b_{k+2})$, hence
$v_{k+1}\in \supp(b_{k+2})$.

If $m>r$, then $a_m\rla b_r$ (by Corollary~\ref{cor:po2}).
Since $v_r\in\supp(a_m)$ and $v_{r-1}\in \supp(b_r)$,
we have $v_r\rla v_{r-1}$,
which contradicts that $\{v_{r-1},v_r\}$ is an edge in $\bar\Gamma$.
Therefore $m\le r$.

Since $\bar\Gamma[g]$ is connected
and both $g_0$ and $g_m$ are nontrivial (by Claim 2),
we may assume that $(v_0,\ldots,v_r)$ is a shortest path
from $v_0\in\supp(g_0)$ to $v_r\in\supp(g_m)$ in $\bar\Gamma[g]$,
hence $r\le\diam(\bar\Gamma[g])$.
Therefore $m\le r\le\diam(\bar\Gamma[g])$.
\end{proof}

Claim 3 proves (i) and Claim 2 proves (ii).

Since $g_0 \neq 1$ and $g_j \rla g_0$ for all $j\ge 2$, we have
$\|(g_m\cdots g_2)\cdots (g_m g_{m-1}) g_m\|_*\le 1$.
Since $g_m\cdots g_1 \le_L g$,
we have $\|g_m\cdots g_1\|_*\le \|g\|_*$.
Therefore
$$
\|u\|_*\le \|g_m\cdots g_1\|_*+ \|(g_m\cdots g_2)\cdots
(g_m g_{m-1}) g_m\|_* \le \|g\|_*+1.
$$
This proves (iii).

Assume $m\ge 3$.
Since $g_0\ne 1$, $g_m\ne 1$, $g_m\cdots g_2\rla g_0$
and $g_1g_0\rla g_m$,
we have $\|g_m\cdots g_2\|_* \le 1$ and $\|g_1g_0\|_* \le 1$.
Therefore $\|g\|_*\le \|g_m\cdots g_2\|_* + \|g_1g_0\|_*\le 2$.
This proves (iv).
\end{myproof}

\begin{remark}\label{rmk:uu}
From the disjoint commutativity $g_i\rla g_j$ for $|i-j|\ge 2$,
the following decompositions are geodesic for all $1\le k\le m$.
\begin{align*}
g^k
&= (g_m\cdots g_0) (g_m\cdots g_0)
\cdots (g_m\cdots  g_0)\\
&= ((g_m\cdots g_1)\cdot (g_0)) ((g_m\cdots g_2)\cdot (g_1 g_0))
\cdots ((g_m\cdots g_k)\cdot (g_{k-1}\cdots g_0))\\
&=(g_m\cdots g_1)(g_m\cdots g_2)\cdots (g_m\cdots g_k)
\cdot (g_0) (g_1g_0)\cdots (g_{k-1}\cdots g_0)\\
&= (u\wedge_L g^k)  (g_0)(g_1g_0)\cdots (g_{k-1}\cdots g_0)
\end{align*}
In particular, $g^m=uu'$ is geodesic, where
\begin{align*}
u &= (g_m\cdots g_1)(g_m\cdots g_2)\cdots (g_m),\\
u' &= (g_0)(g_1g_0)\cdots (g_{m-1}\cdots g_0).
\end{align*}
\end{remark}

The following example shows that the upper bounds
$m\le \diam(\bar\Gamma[g])$ and $m\le |V(\Gamma)|-1$
in Theorem~\ref{thm:po3}(i) are sharp.

\begin{example}
Let $\Gamma=\bar P_4$, where $P_4 = (v_1,\ldots, v_4)$ is a path graph,
and let the underlying right-angled Artin group here be $A(\Gamma)$.
Let $g=v_1^2v_2v_3v_4$
and $u=v_1^2v_2v_3v_1^2v_2v_1$.
Then $g$ is clearly cyclically reduced and non-split.
It is easy to see that $g\not\le_L u$ and $u\not\le_L g^2$.
(By Lemma~\ref{lem:PO},
if $g\le_L u$ then $v_4\le_L v_1^2v_2v_1$,
and if $u\le_L g^2$ then $v_1\le_L v_4v_3v_4$.)
On the other hand, $u\le_L g^3$ because
\begin{align*}
g^3 &= (v_1^2v_2v_3v_4)(v_1^2v_2v_3v_4)(v_1^2v_2v_3v_4)\\
&= (v_1^2v_2v_3\cdot v_4)(v_1^2v_2\cdot v_3v_4)(v_1\cdot v_1v_2v_3v_4)\\
&= (v_1^2v_2v_3\cdot v_1^2v_2\cdot v_1)(v_4\cdot v_3v_4 \cdot v_1v_2v_3v_4)\\
&= u(v_4\cdot v_3v_4 \cdot v_1v_2v_3v_4).
\end{align*}
In the notation of Theorem~\ref{thm:po3},
$$m=3=\diam(\bar\Gamma)=\diam(\bar\Gamma[g]) =|\supp(g)|-1 = | V(\Gamma) |-1.$$
Thus the bounds of $m$ in Theorem~\ref{thm:po3}(i) are sharp.
\end{example}

\begin{corollary}\label{cor:po3}
Let $g, u\in A(\Gamma)$ with $g$ cyclically reduced and non-split.
If $u\le_L g^m$ for some $m\ge 1$,
then $u=g^k a$ is geodesic for some $0\le k\le m$ and $a\in A(\Gamma)$
with $\|a\|_* \le \|g\|_* +1$.
\end{corollary}

\begin{myproof}
We may assume that $u\not\le_L g^{m-1}$.
Let $k = \max\{ l\ge 0 : g^l \le_L u\}$.
Then $0\le k\le m$ and $u=g^k a$ is geodesic
for some $a\in A(\Gamma)$ with $g\not\le_L a \not\le_L g^{m-k-1}$
and $a\le_L g^{m-k}$.

If $m-k\le 1$, then it is obvious that
$\|a\|_* \le \|g\|_* \le \|g\|_* +1$.

If $m-k\ge 2$, then $\|a\|_* \le \|g\|_* +1$
by Theorem~\ref{thm:po3}.
\end{myproof}

For a cyclically reduced element $g\in A(\Gamma)$,
we have seen in Corollary~\ref{cor:st4} that
the sequence $\{\|g^n\|_*\}_{n=0}^\infty$ is increasing.
In particular, if $\|g\|_*\ge 3$, then $\|g^{n+1}\|_*\ge \|g^n\|_*+1$ for all $n\ge 0$.
However, if $\|g\|_*=2$, then it may happen that
$\|g\|_*=\|g^2\|_*=\cdots=\|g^n\|_*=2$ for some $n\ge 2$.
The following proposition finds $m$ with $\|g^m\|_*\ge 3$
when $\|g\|_*= 2$.

\begin{proposition}\label{pro:po4}
Let $g\in A(\Gamma)$ be cyclically reduced
and non-split with  $\|g\|_*=2$.
Then the following hold.
\begin{itemize}
\item[(i)] Let $m\ge 2$.
If either $m\ge |V(\Gamma)|-2$ or $m\ge \diam(\bar\Gamma[g])+1$,
then $\|g^m\|_* \ge 3$.
\item[(ii)] Let $m\ge 2$.
If\/ $\|g^m\|_*= 2$, then $m\le |V(\Gamma)|-3$ and $m\le \diam(\bar\Gamma[g])$.
\item[(iii)] If\/ $|V(\Gamma)|\le 4$, then  $\|g^2\|_*\ge 3$.
\end{itemize}
\end{proposition}

\begin{myproof}
The statements (i) and (ii) are equivalent, and (iii) follows from (i).
Therefore we prove only (ii).

Since $g$ is non-split, $\bar\Gamma[g]$ is connected.
Suppose that $\|g^m\|_*=2$ for some $m\ge 2$.
Then there is a geodesic decomposition
$$g^m=uu'$$
for some $u,u'\in A(\Gamma)$ with
$\|u\|_*=\|u'\|_*=1$.

\medskip\noindent
\textbf{Claim 1.}\ \
$g\not\le_L u\not\le_L g^{m-1}$ and $u\le_L g^m$.

\begin{proof}[Proof of Claim 1]
Since $g^m=uu'$ is geodesic, $u\le_L g^m$.
Since $\|u\|_*=1$ and $\|g\|_*=2$, $g\not\le_L u$.
If $u\le_L g^{m-1}$, then $g\le_R u'$ (by Lemma~\ref{lem:PO}(iv)),
which is impossible because
$\|u'\|_*=1$ and $\|g\|_*=2$.
Therefore $u\not\le_L g^{m-1}$.
\end{proof}

By Claim 1, we can apply Theorem~\ref{thm:po3},
hence $m\le \diam(\bar\Gamma[g])$,
which is the second inequality of (ii).

\medskip
By Theorem~\ref{thm:po3} and Remark~\ref{rmk:uu},
there is a geodesic decomposition
$g=g_mg_{m-1}\cdots g_0$
such that $g_i\ne 1$ for all $0\le i\le m$,
$g_i\rla g_j$ whenever $|i-j|\ge 2$ and
\begin{align*}
u&=(g_m\dots g_1)(g_m\cdots g_2)\cdots (g_m),\\
u'&=(g_0) (g_1g_0)\cdots (g_{m-1}\cdots g_0).
\end{align*}
Since $\|u\|_*=\|u'\|_*=1$, there exist vertices $x, y\in V(\Gamma)$
such that $u\in Z(x)$ and $u'\in Z(y)$,
where $Z(\cdot)$ denotes the centralizer.
Since $\|uu'\|_*=2$, $x\ne y$.
Notice that
\begin{align*}
&\supp(g_1),\ldots,\supp(g_m)\subset Z(x),\\
&\supp(g_0),\ldots,\supp(g_{m-1})\subset Z(y).
\end{align*}

\medskip\noindent
\textbf{Claim 2.}\ \
There is a path
$(x, v_0, v_1, v_2, \ldots, v_{r-1}, v_r, y)$
in $\bar\Gamma$ such that
\begin{enumerate}
\item[(a)] $v_0\in\supp(g_0)$ and $v_r\in\supp(g_m)$;
\item[(b)] the subpath $(v_0,\ldots,v_r)$ is a shortest path from
$v_0$ to $v_r$ in $\bar\Gamma[g]$;
\item[(c)] all the vertices on the path are mutually distinct.
\end{enumerate}

\begin{proof}[Proof of Claim 2]
If either $\supp(g)\subset Z(x)$ or $\supp(g)\subset Z(y)$, then $\|g\|_*=1$,
hence $\supp(g)\setminus Z(x)\ne\emptyset$ and $\supp(g)\setminus Z(y)\ne\emptyset$.

Choose any vertices $v_x\in\supp(g)\setminus Z(x)$
and $v_y\in\supp(g)\setminus Z(y)$,
equivalently, $v_x,v_y\in \supp(g)$ such that
$\{v_x,x\}, \{v_y,y\}\in E(\bar\Gamma)$.
Since $v_x\in \supp(g)\setminus Z(x)=\left(\bigcup_{k=0}^m\supp(g_k)\right)\setminus Z(x)$
and $\bigcup_{k=1}^m\supp(g_k)\subset Z(x)$,
we have  $v_x\in\supp(g_0)$.
Similarly, $v_y\in\supp(g_m)$.
Furthermore, $v_x\ne v_y$ and $\{v_x,v_y\}\not\in E(\bar\Gamma)$
because $g_0\rla g_m$.

Since $\{v_x,x\}, \{v_y,y\}\in E(\bar\Gamma)$ and
$\bar\Gamma[g]$ is connected,
there is a path in $\bar\Gamma$ from $x$ to $y$
$$
(x, v_0=v_x, v_1, v_2, \ldots, v_{r-1}, v_r=v_y, y)
$$
such that $v_k\in \supp(g)$ for all $0\le k\le r$.
Observe that $v_0=v_x\in\supp(g_0)$ and $v_r=v_y\in\supp(g_m)$.

We may assume that it is a shortest path among all the paths
from $x$ to $y$ such that $v_k\in \supp(g)$ for all $0\le k\le r$.
Then the subpath $(v_0,\cdots, v_r)$ must be a shortest path from $v_0$ to $v_r$
in $\bar\Gamma[g]$, hence $v_0,\ldots, v_r$ are mutually distinct.

If $x=v_j$ for some $0\le j\le r-1$,
then the path $(x, v_{j+1},\ldots, v_r, y)$ is shorter than the original one,
all of whose middle vertices belong to $\supp(g)$.
This contradicts that $(x,v_0,\cdots,v_r,y)$ is a shortest path among such paths.
If $x=v_r$ then $\{x,v_x\}\not\in E(\bar\Gamma)$
(because $x=v_r\in\supp(g_m)$, $v_x=v_0\in\supp(g_0)$ and $g_0\rla g_m$).
This is a contradiction.
Therefore $x\ne v_j$ for any $0\le j\le r$.
Similarly, $y\ne v_j$ for any $0\le j\le r$.
Since $x\ne y$, all the vertices on the path
$(x, v_0, v_1, v_2, \ldots, v_{r-1}, v_r, y)$ are mutually distinct.
\end{proof}

Since the $(r+3)$ points on the path in Claim 2 are mutually distinct,
$| V(\Gamma) |\ge r+3$.
By Claim 3 in the proof of Theorem~\ref{thm:po3}, we have $m\le r$.
Therefore $| V(\Gamma) |\ge r+3\ge m+3$.
This proves the first inequality of (ii), hence (ii) is proved.
\end{myproof}

The following example illustrates that
the upper bounds $m\le | V(\Gamma) |-3$ and
$m\le \diam(\bar\Gamma[g])$
in Proposition~\ref{pro:po4}(ii) are sharp.

\begin{example}
Let $\Gamma = \bar P_6$, where $P_6 = (v_0, v_1, \ldots, v_5)$
is a path graph,
and let the underlying right-angled Artin group here be $A(\Gamma)$.
Let $g=v_1v_2v_3v_4$. It is easy to see that $\|g\|_*=2$.
Since $\supp(g)=\{v_1,v_2,v_3,v_4\}$, $\diam(\bar\Gamma[g])=3$.
Observe
\begin{align*}
g^3 &= (v_1v_2v_3v_4)\cdot(v_1v_2v_3v_4)\cdot(v_1v_2v_3v_4)\\
&= (v_1v_2v_3\cdot v_4)\cdot (v_1v_2\cdot v_3v_4)\cdot (v_1\cdot v_2v_3v_4)\\
&= (v_1v_2v_3\cdot v_1v_2\cdot v_1)(v_4\cdot v_3v_4\cdot v_2v_3v_4).
\end{align*}
Let $u=v_1v_2v_3v_1v_2v_1$ and $u'=v_4v_3v_4v_2v_3v_4$.
Then $g^3=uu'$ is geodesic.
Since $u\in Z(v_5)$ and $u'\in Z(v_0)$,
we have $\|u\|_*=\|u'\|_*=1$, hence $\|g^3\|_*=2$.
Notice that $3=| V(\Gamma) |-3=\diam(\bar\Gamma[g])$.
\end{example}

\EKLnewpage

\section{Asymptotic translation length}
\label{sec:tl}

In this section, we study asymptotic translation lengths
of elements of $A(\Gamma)$ on $(A(\Gamma),d_*)$ and on $(\Gamma^e,d)$,
and then find a lower bound of the minimal asymptotic translation
length of $A(\Gamma)$ on $\Gamma^e$.

\begin{proposition}\label{pro:tl}
If $g\in A(\Gamma)$ is cyclically reduced
and non-split with  $\|g\|_*\ge 2$,
then
$$
\tau_{(A(\Gamma),d_*)}(g)\ge\frac1{\max\{2,\,|V(\Gamma)|-2\}}\,.
$$
\end{proposition}

\begin{myproof}
Let $\tau_*$ denote $\tau_{(A(\Gamma),d_*)}$,
and let $V=|V(\Gamma)|$.

Notice that if $\|g\|_*\ge 3$, then  $\|g^n\|_*\ge n+2$ for all $n\ge 1$ (by Corollary~\ref{cor:st4}), hence
$$\tau_*(g) =\lim_{n\to\infty}\frac{\|g^n\|_*}{n}
\ge\lim_{n\to\infty}\frac{n+2}{n} \ge 1.$$

Suppose $\|g\|_*=2$.
By Proposition~\ref{pro:po4}, if $V\le 4$ then $\|g^2\|_*\ge 3$,
and if $V\ge 5$ then $\|g^{V-2}\|_* \ge 3$.
Therefore $\|g^{\max\{2,V-2\}}\|_*\ge 3$.
From the above discussion, $\tau_*(g^{\max\{2,V-2\}})\ge 1$
and hence $\tau_*(g)\ge \frac1{\max\{2,V-2\}}$.
\end{myproof}

\begin{remark}\label{rmk:st-ass}
When we study the action of $A(\Gamma)$ on $(A(\Gamma),d_*)$,
we will assume that
``$|V(\Gamma)|\ge 2$ and $\bar\Gamma$ is connected''
because otherwise $\|g\|_*\le 2$ for all $g\in A(\Gamma)$
and hence $(A(\Gamma),d_*)$ has diameter at most 2:
if $|V(\Gamma)|=1$, then $\|g\|_*\le 1$ for all $g\in A(\Gamma)$;
if $\bar\Gamma$ is disconnected (i.e.\ $\Gamma$ is a join),
then $\|g\|_*\le 2$ for all $g\in A(\Gamma)$.
\end{remark}

\begin{lemma}
\label{lem:loxo}
Suppose that $|V(\Gamma)|\ge 2$ and $\bar\Gamma$ is connected.
The following are equivalent
for a cyclically reduced element $g\in A(\Gamma)$.

\begin{enumerate}
\item $g$ is strongly non-split and $| \supp(g) | \ge 2$.
\item $g$ is non-split and $\|g\|_*\ge 2$.
\item $\|g^n\|_*\ge 3$ for some $n\ge 1$.
\item $g$ is loxodromic on $(A(\Gamma),d_*)$,
i.e.~$\tau_{(A(\Gamma),d_*)}(g)>0$.
\item $\tau_{(A(\Gamma),d_*)}(g)\ge \frac{1}{\max\{2,\,|V(\Gamma)|-2\}}$.
\end{enumerate}
\end{lemma}

\begin{myproof}[Proof of Lemma~\ref{lem:loxo}]
(i)$\Leftrightarrow$(ii) follows from Lemma~\ref{lem:SNS}.

(ii)$\Rightarrow$(v) follows from Proposition~\ref{pro:tl}.

(v)$\Rightarrow$(iv)
and (iv)$\Rightarrow$(iii) are obvious.

(iii)$\Rightarrow$(i):
Since $\|g^n\|_*\ge 3$, $g^n$ is strongly non-split (by Lemma~\ref{lem:SNS}),
hence $g$ is also strongly non-split (see Remark~\ref{rmk:sp-po}).
It is obvious that $|\supp(g)|\ge 2$.
\end{myproof}

\begin{remark}\label{rmk:loxo}
Suppose that $|V(\Gamma)|\ge 4$ and
both $\Gamma$ and $\bar\Gamma$ are connected.
Then the condition $|\supp(g)|\ge 2$
of Lemma~\ref{lem:loxo}(i) is not necessary
because all strongly non-split elements $g$
must have $|\supp(g)|\ge 2$.
Moreover, $g$ is loxodromic on $(A(\Gamma),d_*)$
if and only if it is loxodromic on $(\Gamma^e,d)$ by Corollary~\ref{cor:QI}.
Therefore, if $|V(\Gamma)|\ge 4$ and
both $\Gamma$ and $\bar\Gamma$ are connected,
then (i) and (iv) in the above lemma are equivalent to
the following (i$'$) and (iv$'$) respectively.
\begin{itemize}
\item[(i$'$)] $g$ is strongly non-split.
\item[(iv$'$)] $g$ is loxodromic on $(\Gamma^e,d)$,
i.e.~$\tau_{(\Gamma^e,d)}(g)>0$.
\end{itemize}
\end{remark}

Kim and Koberda~\cite[Lemma 33]{KK14} showed that if
$g \in A(\Gamma)$ is cyclically reduced
and strongly non-split,
then $\|g^{2n|V(\Gamma)|^2}\|_*\ge n$
for all $n\ge 1$. Therefore (by Corollary~\ref{cor:QI})
$$\tau_{(\Gamma^e,d)}(g)
\ge \tau_{(A(\Gamma),d_*)}(g) \ge \frac{1}{2|V(\Gamma)|^2}.
$$
From this, a lower bound of the minimal asymptotic
translation length of $A(\Gamma)$ on $\Gamma^e$ follows:
$$
\minTL_{(\Gamma^e,\,d)}(A(\Gamma)) \ge
\frac1{2|V(\Gamma)|^2}.
$$
We improve the denominator of the lower bound
from a quadratic function to a linear function
of $|V(\Gamma)|$ as follows.

\begin{theorem}\label{thm:mTL}
Let $\Gamma$ be a finite simplicial graph
such that $|V(\Gamma)|\ge 4$ and
both $\Gamma$ and $\bar\Gamma$ are connected.
Then
$$
\minTL_{(\Gamma^e,\,d)}(A(\Gamma))\ge
\frac1{|V(\Gamma)|-2}\,.
$$
\end{theorem}

\begin{myproof}
Since $|V(\Gamma)|\ge 4$,
$|V(\Gamma)|-2=\max\{2,\,|V(\Gamma)|-2\}$.
Let $g\in A(\Gamma)$ be loxodromic on $(\Gamma^e,d)$ and hence on $(A(\Gamma),d_*)$
(by Remark~\ref{rmk:loxo}).
We may assume that $g$ is cyclically reduced
because asymptotic translation lengths are
invariant under conjugation.
By Corollary~\ref{cor:QI} and Lemma~\ref{lem:loxo},
$$
\tau_{(\Gamma^e,d)}(g)
\ge \tau_{(A(\Gamma),d_*)}(g)
\ge \frac{1}{|V(\Gamma)|-2}.
$$
Therefore $\minTL_{(\Gamma^e,\,d)}(A(\Gamma))\ge \frac1{|V(\Gamma)|-2}$.
\end{myproof}

\EKLnewpage

\section{Uniqueness of quasi-roots}
\label{sec:quasiroot}

The notion of quasi-roots in $A(\Gamma)$ was introduced in~\cite{LL22},
where the quasi-roots are defined using word length.
The uniqueness up to conjugacy was established by using
the normal form of elements introduced by Crisp, Godelle and Wiest~\cite{CGW09}.
In this section, we extend the uniqueness of quasi-roots
from word length to star length.

\begin{definition}(quasi-root)\label{def:qr}
An element $g\in A(\Gamma)\backslash\{1\}$ is called a \emph{quasi-root} of $h\in A(\Gamma)$
if there is a decomposition
$$h=ag^{n}b$$
for some $n\ge 1$ and $a,b\in A(\Gamma)$
such that $\|h\|=\|a\|+\|b\|+n\|g\|$.
The decomposition is called a \emph{quasi-root decomposition} of $h$.
The conjugates $aga^{-1}$ and $b^{-1}gb$ are called
the \emph{leftward-extraction} and the \emph{rightward-extraction}
of the quasi-root $g$, respectively.
We consider the following two cases.
\begin{enumerate}
\item
$g$ is called an \emph{$(A,B,r)$-quasi-root} of $h$
if $\|a\|\le A,~\|b\|\le B$
and $\|g\|\le r$.

\item
$g$ is called an \emph{$(A,B,r)^*$-quasi-root} of $h$
if $\|a\|_*\le A,~\|b\|_*\le B$
and $\|g\|_*\le r$.
\end{enumerate}
\end{definition}

In the above definition,
the condition $\|h\|=\|a\|+\|b\|+n\|g\|$
implies $\|g^n\|=n\|g\|$,
hence $g$ is cyclically reduced when $n\ge 2$.

Notice that if $g_1=aga^{-1}$ and $g_2=b^{-1}gb$ are respectively
the leftward- and the rightward-extractions of $g$,
then we have decompositions $h=g_1^n ab=abg_2^n$,
which are not necessarily geodesic.

\begin{definition}[primitive]
An element $g\in A(\Gamma)\backslash\{1\}$ is called \emph{primitive}
if $g$ is not a nontrivial power of another element,
i.e.\ $g=h^n$ never holds for any $n\ge 2$ and $h\in A(\Gamma)$.
\end{definition}

The following proposition is
Proposition 3.5 in~\cite{LL22} written
in the setting of this paper.
It shows a kind of uniqueness property of quasi-roots in right-angled Artin groups.

\begin{proposition}[{\cite[Proposition 3.5]{LL22}}]\label{pro:qr-wd}
Let $h\in A(\Gamma)$, $A,B\ge 0$  and  $r\ge 1$.
If
$$\|h\|\ge A+B+(2V+1)r,$$
where $V=|V(\Gamma)|$,
then strongly non-split and primitive $(A,B,r)$-quasi-roots of $h$
are conjugate to each other,
and moreover, their leftward- and rightward-extractions are unique.
\end{proposition}

In other words, Proposition~\ref{pro:qr-wd} shows that if
$$
h=a_1g_1^{n_1}b_1=a_2g_2^{n_2}b_2
$$
are two quasi-root decompositions of $h$
such that for each $i=1,2$, $g_i$ is strongly non-split and primitive,
$$
\|a_i\|\le A,\quad
\|b_i\|\le B,\quad
\|g_i\|\le r,\atop
\|h\|\ge A+B+(2V+1)r,
$$
then $g_1$ and $g_2$ are conjugate, and moreover,
$a_1g_1a_1^{-1}=a_2g_2a_2^{-1}$ and $b_1^{-1}g_1b_1=b_2^{-1}g_2b_2$.

\smallskip
The following theorem is the main result of this section.
It is a star length version of Proposition~\ref{pro:qr-wd},
which plays an important role in the proof of Theorem~\ref{thm:ex}.
We remark that the word length and
the star length are quite independent,
hence the word length version does not naively extend to
a star length version.
We exploit lattice structures developed in~\S\ref{sec:lattice}.

We also remark that in the following theorem
since $\|h\|_* \ge 2A+2B+(2V+3)r+2\ge 3r+2\ge 5$,
the existence of such an element $h$ implies that
$|V(\Gamma)|\ge 2$ and $\bar\Gamma$ is connected
(as explained in Remark~\ref{rmk:st-ass}).

\EKLnewpage

\begin{theorem}\label{thm:qr}
Let $h\in A(\Gamma)$, $A,B\ge 0$  and  $r\ge 1$.
If
$$\|h\|_* \ge 2A+2B+(2V+3)r+2,$$
where $V=|V(\Gamma)|$, then primitive $(A,B,r)^*$-quasi-roots of $h$
are conjugate to each other, and moreover,
their leftward- and rightward-extractions are unique.
In other words, if
$$
h=a_1g_1^{n_1}b_1=a_2g_2^{n_2}b_2
$$
are two quasi-root decompositions of $h$
such that for each $i=1, 2$, $g_i$ is primitive and
$$
\|a_i\|_*\le A,\quad
\|b_i\|_*\le B,\quad
\|g_i\|_*\le r,\atop
\|h\|_* \ge 2A+2B+(2V+3)r+2,
$$
then $g_1$ and $g_2$ are conjugate to each other such that
$$
a_1g_1a_1^{-1}=a_2g_2a_2^{-1}\quad\mbox{and}\quad
b_1^{-1}g_1b_1=b_2^{-1}g_2b_2.
$$
\end{theorem}

\begin{myproof}
Let $i=1, 2$.

\medskip\noindent
\textbf{Claim 1.}\ \
$n_i\ge 4$ and $g_i$ is cyclically reduced and
strongly non-split with $\|g_i\|_*\ge 2$.

\begin{proof}[Proof of Claim 1]
If $n_i\le 3$, then
\begin{align*}
\|h\|_*
&=\|a_i g_i^{n_i} b_i \|_*
\le \|a_i\|_* +n_i \|g_i\|_* +\|b_i \|_*
\le A+3r+B.
\end{align*}
This contradicts the hypothesis $\|h\|_*\ge 2A+2B+(2V+3)\, r+2$.
Therefore $n_i\ge 4$ and hence
$g_i$ is cyclically reduced (see the paragraph following Definition~\ref{def:qr}).

Observe that $\|g_i^{n_i}\|_*\ge 5$ because
\begin{align*}
\|g_i^{n_i}\|_*
&\ge \|h\|_*-\|a_i\|_*-\|b_i\|_*
\ge A+B+(2V+3)r+2 \ge 3r+2\ge 5.
\end{align*}
Therefore $g_i$ is strongly non-split
and $\|g_i\|_*\ge 2$ (by Lemma~\ref{lem:loxo}).
\end{proof}

Let $\alpha_i$ and $\beta_i$ be integers defined by
\begin{align*}
\alpha_i & =\min\{k\ge 1: \|g_i^k\|_*\ge A+2\},\\
\beta_i & =\min\{k\ge 1: \|g_i^k\|_*\ge B+2\}.
\end{align*}
The numbers $\alpha_i$ and $\beta_i$ are well-defined
because the sequence $\{\|g_i^k\|_*\}_{k=1}^\infty$ is
increasing such that
$\lim_{k\to\infty}\|g_i^k\|_*=\infty$
(by Claim 1, Lemmas~\ref{lem:st0} and~\ref{lem:loxo}).
Since $\|g_i^k\|_*-\|g_i^{k-1}\|_*
\le \|g_i\|_*\le r$ for all $k\ge 1$,
we get
\begin{align*}
A+2 & \le \|g_i^{\alpha_i}\|_* \le A+1+r ,\\
B+2 & \le \|g_i^{\beta_i}\|_*  \le B+1+r .
\end{align*}

\bigskip
\noindent\textbf{Claim 2.}\ \
$n_i-\alpha_i-\beta_i \ge 2V+1$.

\begin{proof}[Proof of Claim 2]

Observe that
\begin{align*}
\|g_i^{n_i}\|_* &\ge \|h\|_*-\|a_i\|_*-\|b_i\|_*\ge \|h\|_* -A-B, \\
\|g_i^{n_i}\|_*
 - \|g_i^{\alpha_i+\beta_i}\|_*
&\ge \|g_i^{n_i}\|_* - \|g_i^{\alpha_i}\|_* - \|g_i^{\beta_i}\|_*\\
&\ge (\|h\|_* -A-B) - (A+1+r) - (B+1+r)\\
& = \|h\|_* -2A-2B-2r-2 \ge (2V+1)r.
\end{align*}
Since
$\{\|g_i^k\|_*\}_{k=1}^\infty$ is increasing
and $\|g_i^{n_i}\|_* - \|g_i^{\alpha_i+\beta_i}\|_* \ge (2V+1)r>0$,
we have $n_i> \alpha_i+\beta_i$.
Since
\begin{align*}
(n_i-\alpha_i-\beta_i)\|g_i\|_*
&\ge \|g_i^{n_i-\alpha_i-\beta_i}\|_*
\ge\|g_i^{n_i}\|_* - \|g_i^{\alpha_i+\beta_i}\|_*\\
&\ge (2V+1)r
\ge (2V+1) \|g_i \|_*,
\end{align*}
we get $n_i-\alpha_i-\beta_i \ge 2V+1$ as desired.
\end{proof}

Let $a_0=a_1\wedge_L a_2$ and $b_0=b_1\wedge_R b_2$.
Then we have geodesic decompositions
$$
\left\{\begin{array}{l}
a_1=a_0a_1',\\
a_2=a_0a_2',
\end{array}\right.
\qquad
\left\{\begin{array}{l}
b_1=b_1'b_0,\\
b_2=b_2'b_0
\end{array}\right.
$$
for some $a_1', a_2', b_1',b_2'\in A(\Gamma)$
with $a_1'\wedge_L a_2'=1$ and $b_1'\wedge_R b_2'=1$.
Observe that
\begin{align*}
\|a_i'\|_* &\le  \|a_i\|_*\le A,\qquad\\
\|b_i'\|_* &\le  \|b_i\|_*\le B.
\end{align*}
Since $h=a_1g_1^{n_1}b_1=a_2g_2^{n_2}b_2$,
we have
$h=a_0(a_1'g_1^{n_1}b_1')b_0=a_0(a_2'g_2^{n_2}b_2')b_0$.

Let $h_0=a_0^{-1}hb_0^{-1}$.
Then $h_0$ has the following two geodesic decompositions.
\begin{equation}\label{E:qr}
h_0 =a_1'g_1^{n_1}b_1'=a_2'g_2^{n_2}b_2'
\end{equation}
On the other hand, since $a_1$ and $a_2$ have a common right multiple, say $h$, we have
$a_1'\rla a_2'$ (by Theorem~\ref{thm:lcm}).
Since  $a_1'\le_L  a_2'g_2^{n_2}b_2'$ and $a_2'\le_L  a_1'g_1^{n_1}b_1'$,
we have (by Lemma~\ref{lem:gcd2})
$$
a_1'\le_L g_2^{n_2}b_2'  \quad\mbox{and}\quad a_2'\le_L g_1^{n_1}b_1'.
$$

Let $A'=\|a_1'\|+\|a_2'\|$ and $B'=\|b_1'\|+\|b_2'\|$.

\bigskip
\noindent\textbf{Claim 3.}\ \
$\|h_0\|\ge A'+B'+(2V+1) \|g_i\|$.

\begin{proof}[Proof of Claim 3]
We know that  $n_1-\alpha_1>0$ (by Claim 2) and
$a_2'\le_L g_1^{n_1}b_1'=g_1^{\alpha_1}\cdot g_1^{n_1-\alpha_1}b_1'$.
Since $g_1^{\alpha_1}\cdot g_1^{n_1-\alpha_1}b_1'$ is geodesic
and $\|g_1^{\alpha_1}\|_*\ge A+2\ge \|a_2'\|_*+2$,
we have $a_2'\le_L g_1^{\alpha_1}$ (by Corollary~\ref{cor:st2}).
Similarly, $b_2' \le_R g_1^{\beta_1}$.
(In other words, $g_1^{n_1}$ has a geodesic decomposition
$g_1^{n_1}=g_1^{\alpha_1}\cdot g_1^{n_1-\alpha_1-\beta_1}\cdot g_1^{\beta_1}$
such that $a_2'\le_L g_1^{\alpha_1}$ and $b_2'\le_R g_1^{\beta_1}$.)
Therefore
\begin{align*}
&\|a_2'\|\le \|g_1^{\alpha_1}\|=\alpha_1\|g_1\|,\qquad\\
&\|b_2'\|\le \|g_1^{\beta_1 }\|=\beta_1\|g_1\|.
\end{align*}
Since $h_0=a_1'g_1^{n_1}b_1'=a_2'g_2^{n_2}b_2'$, we get
\begin{align*}
\|h_0\|-(A'+B')
&= (\|a_1'\|+n_1\|g_1\|+\|b_1'\|) - (\|a_1'\|+\|a_2'\|) - (\|b_1'\|+\|b_2'\|)\\
&= n_1\|g_1\|-\|a_2'\|-\|b_2'\|
\ge n_1\|g_1\|-\alpha_1\|g_1\|-\beta_1\|g_1\|\\
&=(n_1-\alpha_1-\beta_1)\|g_1\|
\ge (2V+1) \|g_1\|.
\end{align*}
In the same way, we get $\|h_0\|-(A'+B') \ge (2V+1) \|g_2\|$.

\vskip -\baselineskip
\end{proof} 

Notice that $\|a_i'\|\le A'$ and $\|b_i'\|\le B'$.
Let $r'=\max \{ \|g_1\|, \|g_2\| \}$.
Then each $a_i'g_i^{n_i}b_i'$ in \eqref{E:qr} is a
$(A',B',r')$-quasi-root decomposition of $h_0$
such that $\|h_0\|\ge A'+B'+(2V+1)r'$.

Applying Proposition~\ref{pro:qr-wd} to \eqref{E:qr} yields
$a_1'g_1a_1'^{-1}=a_2'g_2a_2'^{-1}$
and $b_1'^{-1}g_1b_1'=b_2'^{-1}g_2b_2'$.
Consequently
\begin{align*}
&a_1g_1a_1^{-1}
=a_0(a_1'g_1a_1'^{-1})a_0^{-1}
=a_0(a_2'g_2a_2'^{-1})a_0^{-1}
=a_2g_2a_2^{-1}, \\
&
b_1^{-1}g_1b_1
= b_0^{-1}(b_1'^{-1}g_1b_1')b_0
= b_0^{-1}(b_2'^{-1}g_2b_2')b_0
=b_2^{-1}g_2b_2.
\end{align*}
\vskip -\baselineskip
\end{myproof}

\EKLnewpage

\section{Acylindricity of the action of $A(\Gamma)$ on $\Gamma^e$}
\label{sec:acyl}

In this section, we prove the following two theorems.

\begin{theorem}\label{thm:st}
Let $\Gamma$ be a finite simplicial graph
such that $|V(\Gamma)|\ge 2$ and $\bar\Gamma$ is connected.
Then the action of\/ $A(\Gamma)$ on $(A(\Gamma), d_*)$ is
$(R, N)$-acylindrical with
\begin{align*}
&R=R(r) =(2V+7)r+8,\\
&N=N(r) =2(V-2)(r-1)-1,
\end{align*}
where $V=\max\{4,|V(\Gamma)|\}$.
Moreover, for any $x,y\in A(\Gamma)$ with $d_*(x,y)\ge R$,
if $\xi(x,y;r)\ne \{1\}$, then
there exists a loxodromic element $g\in A(\Gamma)$ such that
\begin{enumerate}
\item
$\xi(x,y;r)=\{1,g^{\pm 1}, g^{\pm 2},\ldots, g^{\pm k}\}$
for some $1\le k\le (V-2)(r-1)-1$;

\item
the Hausdorff distance between the $\langle g\rangle$-orbit of $x$ and
that of $y$ is at most $2r+3$.
\end{enumerate}
\end{theorem}

\begin{theorem}
\label{thm:ex}
Let $\Gamma$ be a finite simplicial graph
such that $|V(\Gamma)|\ge 4$
and both $\Gamma$ and $\bar\Gamma$ are connected.
Then the action of\/ $A(\Gamma)$ on $\Gamma^e$ is
$(R, N)$-acylindrical with
\begin{align*}
&R=R(r) =D(2V+7)(r+1)+10D,\\
&N=N(r) =2(V-2)r-1,
\end{align*}
where $D=\diam(\Gamma)$ and $V=|V(\Gamma)|$.
Moreover, for any $x,y\in V(\Gamma^e)$ with
$d(x,y)\ge R$, if\/ $\xi(x,y;r)\ne \{1\}$, then
there exists a loxodromic element $g\in A(\Gamma)$
such that
\begin{enumerate}
\item
$\xi(x,y;r)\subset\{1,g^{\pm 1}, g^{\pm 2},\ldots, g^{\pm k}\}$
for some $1\le k\le (V-2)r-1$;
\item
the Hausdorff distance between the $\langle g\rangle$-orbit of $x$ and
that of $y$ is at most $D(2r+7)$.
\end{enumerate}

\end{theorem}

The following lemma connects the acylindricities of
the actions of $A(\Gamma)$ on $(A(\Gamma), d_*)$ and on $(\Gamma^e, d)$.
It is an improvement of the argument of Kim and Koberda in the proof of Theorem 30 in~\cite{KK14}.

\begin{lemma}\label{lem:2act}
Suppose that $|V(\Gamma)|\ge 4$
and both $\Gamma$ and $\bar\Gamma$ are connected.
Let $D=\diam(\Gamma)$.
If\/ the action of $A(\Gamma)$ on $(A(\Gamma), d_*)$
is $(R_1(r), N_1(r))$-acylindrical,
then the action of $A(\Gamma)$ on $(\Gamma^e, d)$
is $(R_2(r), N_2(r))$-acylindrical
with
\begin{align*}
R_2(r) &=D\cdot R_1(r+1)+2D,\\
N_2(r) &=N_1(r+1).
\end{align*}
More precisely, for any $v_1^{w_1}, v_2^{w_2} \in V(\Gamma^e)$,
where $v_1, v_2\in V(\Gamma)$ and $w_1, w_2\in A(\Gamma)$,
if\/ $d(v_1^{w_1}, v_2^{w_2})\ge R_2(r)$,
then
\begin{enumerate}
\item
$d_* (w_1, w_2) \ge R_1(r+1)$;

\item
$\xi_{(\Gamma^e,d)}(v_1^{w_1}, v_2^{w_2};r)$ is contained in
$\xi_{(A(\Gamma),d_*)}(w_1, w_2; r+1)$.
\end{enumerate}
\end{lemma}

\begin{myproof}
Note that $D=\diam(\Gamma)\ne 0$.
Let $d(v_1^{w_1}, v_2^{w_2})\ge R_2(r)$
for $v_1^{w_1}, v_2^{w_2} \in V(\Gamma^e)$.
Since $\Gamma$ is connected, we can apply Lemma~\ref{lem:KK}
and obtain
\begin{align*}
d_* &(w_1, w_2)
=\|w_2 w_1^{-1}\|_*\\
&\ge \frac{d(v_2, v_2^{w_2 w_1^{-1}})}{D}-1
\ge \frac{d(v_1, v_2^{w_2 w_1^{-1}})-d(v_1, v_2)}{D}-1\\
&\ge \frac{d(v_1^{w_1}, v_2^{w_2})-D}{D}-1
\ge \frac{R_2(r)-2D}D
=R_1(r+1),
\end{align*}
which proves (i).

Let $g\in\xi_{(\Gamma^e,d)}(v_1^{w_1}, v_2^{w_2};r)$.
Then $d(v_i^{w_i g}, v_i^{w_i})\le r$ for $i=1,2$.
By Lemma~\ref{lem:KK} again,
\begin{align*}
d_*(w_i g,w_i ) &
= \|w_i g w_i^{-1}\|_*
\le d(v_i^{w_i g w_i^{-1}}, v_i)+1\\
&= d(v_i^{w_ig}, v_i^{w_i})+1 \le r+1
\end{align*}
for $i=1,2$, hence $g\in \xi_{(A(\Gamma),d_*)}(w_1,w_2;r+1)$.
This shows that
the set $\xi_{(\Gamma^e,d)}(v_1^{w_1}, v_2^{w_2};r)$ is contained in
$\xi_{(A(\Gamma),d_*)}(w_1, w_2; r+1)$,
hence (ii) is proved.

Since
$\xi_{(\Gamma^e,d)}(v_1^{w_1}, v_2^{w_2};r)
\subset \xi_{(A(\Gamma),d_*)}(w_1, w_2; r+1)$
and $d_*(w_1,w_2) \ge R_1(r+1)$,
the $(R_1,N_1)$-acylindricity of
the action of $A(\Gamma)$ on $(A(\Gamma),d_*)$
implies that
\begin{align*}
|\xi_{(\Gamma^e,d)}(v_1^{w_1}, v_2^{w_2};r)|
&\le |\xi_{(A(\Gamma),d_*)}(w_1,w_2;r+1)|
\le N_1(r+1)=N_2(r).
\end{align*}
Therefore the action of $A(\Gamma)$ on $(\Gamma^e, d)$
is $(R_2(r), N_2(r))$-acylindrical.
\end{myproof}

\begin{proposition}\label{pro:A}
Let $g,w\in A(\Gamma)\setminus\{1\}$ and $r,R\ge 1$ such that
$$
\|g\|_*\le r,\quad
\|w^{-1}gw\|_*\le r,\quad
\|w\|_*\ge R,\quad
R\ge 3r+7.
$$
Then there exists a quasi-root decomposition
$$
w=a(g_1^\epsilon)^{n}b,
$$
where $a,b,g_1\in A(\Gamma)$,
$\epsilon\in\{\pm 1\}$ and $n\ge 2$
such that
\begin{itemize}
\item[(i)] $\|a\|_*\le\frac 12 r+1$ and $\|b\,\|_*\le\frac32r+2$;
\item[(ii)] $g_1$ is cyclically reduced and $g=ag_1a^{-1}$ is geodesic.
\end{itemize}
\end{proposition}

Notice that $\|w\|_*\ge R\ge 3r+7\ge 7$,
hence the existence of such an element $w$ implies
that $|V(\Gamma)|\ge 2$ and $\bar \Gamma$ is connected
(as explained in Remark~\ref{rmk:st-ass}).

\begin{myproof}
Let $g=ag_1a^{-1}$ be the geodesic decomposition such that $g_1$ is cyclically reduced.
Let $h=w^{-1}gw$.
Then
$$
h=w^{-1}ag_1a^{-1}w
=(a^{-1}w)^{-1} g_1 (a^{-1}w).
$$
By Theorem~\ref{thm:cnj1},
there exists a geodesic decomposition of $a^{-1}w$
\begin{equation}\label{E:A}
a^{-1}w = w_1w_2w_3
\end{equation}
such that
(i) $w_1\rla g_1$;
(ii) $g_1^{w_2}$ is a cyclic conjugation;
(iii) $h=w_3^{-1}\cdot g_1^{w_2}\cdot w_3$ is geodesic.

\medskip\noindent
\textbf{Claim 1.}\ \
The following hold.
\begin{enumerate}
\item $\|w_2\|_*\ge 3$, hence $w_2$ is strongly non-split.
\item $g_1^{w_2}$ is either a left cyclic conjugation or a right cyclic conjugation.
\end{enumerate}

\begin{proof}[Proof of Claim 1]
Since both $g=ag_1a^{-1}$ and $h=w_3^{-1}g_1^{w_2}w_3$ are geodesic decompositions,
\begin{align*}
\|g_1\|_*+2\|a\|_* -4         &\le \|g\|_* \le \|g_1\|_*+2\|a\|_*,\\
\|g_1^{w_2}\|_*+2\|w_3\|_* -4   &\le \|h\|_* \le \|g_1^{w_2}\|_*+2\|w_3\|_*
\end{align*}
(by Corollary~\ref{cor:st3}), whence
\begin{align*}
\frac{\|g\|_*- \|g_1\|_*}2       &\le \|a\|_*\le  \frac{\|g\|_*- \|g_1\|_*}2+2,\\
\frac{\|h\|_*- \|g_1^{w_2}\|_*}2 &\le \|w_3\|_*\le  \frac{\|h\|_*- \|g_1^{w_2}\|_*}2+2.
\end{align*}

Since $w_1\rla g_1\neq 1$, we have $\|w_1\|_*\le 1$.
Since $g_1\ne 1$ and both $g=ag_1a^{-1}$ and $h=w_3^{-1}g_1^{w_2}w_3$ are geodesic,
we have $1\le \|g_1\|_*\le \|g\|_*\le r$
and $1\le \|g_1^{w_2}\|_*\le \|h\|_* \le r$.
Since $a^{-1}w=w_1w_2w_3$,
\begin{align*}
\|w\|_* &\le \|a\|_* + \|w_1\|_* + \|w_2\|_* + \|w_3\|_*\\
&\le
\left(\frac{\|g\|_*- \|g_1\|_*}2+2\right) + 1 + \|w_2\|_* +
\left(\frac{\|h\|_*- \|g_1^{w_2}\|_*}2+2\right)\\
&\le
\left(\frac{r- 1}2+2\right) + 1 + \|w_2\|_* +
\left(\frac{r- 1}2+2\right)\\
&= \|w_2\|_* + r+4.
\end{align*}
Therefore $\|w_2\|_*\ge \|w\|_*-r-4\ge R-r-4\ge 2r+3\ge 3$,
hence $w_2$ is strongly non-split (by Lemma~\ref{lem:SNS}).
This proves (i).

\smallskip
Assume that the cyclic conjugation $g_1^{w_2}$
is neither a left cyclic conjugation nor a right cyclic conjugation.
Then, by Proposition~\ref{lem:cy2}(v), $w_2=w_2'w_2''$ is geodesic
for some $w_2',w_2''\in A(\Gamma)\setminus\{1\}$
such that
$g_1^{w_2'}$ (resp.\ $g_1^{w_2''}$) is a left (resp.\ right) cyclic conjugation
and $w_2'\rla w_2''$.
Since both $w_2'$ and $w_2''$ are nontrivial, we have $\|w_2'\|_*=\|w_2''\|_*=1$,
hence $\|w_2\|_*=\|w_2'w_2''\|_*\le  \|w_2'\|_* +\|w_2''\|_*=2$,
which contradicts $\|w_2\|_*\ge 3$.
Therefore $g_1^{w_2}$
is either a left cyclic conjugation or a right cyclic conjugation.
This proves (ii).
\end{proof}

\medskip\noindent
\textbf{Claim 2.}\ \
The following hold.
\begin{enumerate}
\item $g_1$ is strongly non-split with $|\supp(g_1)| \ge 2$
and $2\le \|g_1\|_*\le r$.
\item $\|a\|_*\le \frac 12 r+1$,
$w_1=1$, $\|w_2\|_*\ge R-r-2$ and $\|w_3\|_*\le \frac 12 r+1$.
\end{enumerate}

\begin{proof}[Proof of Claim 2]
(i)\ \
Since $g_1^{w_2}$ is either a left or a right cyclic conjugation (by Claim 1),
$$
w_2\le_L g_1^n\quad\mbox{or}\quad
w_2^{-1}\le_R g_1^n
$$
for some $n\ge 1$ (by Proposition~\ref{prop:cc1}).
Since $\|w_2\|_*\ge 3$, we have $\|g_1^n\|_*\ge 3$.
By Lemma~\ref{lem:loxo},
$g_1$ is strongly non-split with $|\supp(g_1)|\ge 2$ and
$\|g_1\|_*\ge 2$.
On the other hand, $\|g_1\|_*\le \|g\|_*\le r$
because $g=ag_1a^{-1}$ is geodesic.

\smallskip
(ii)\ \
If $w_1\ne 1$, then $\|g_1\|_*\le 1$ because $g_1\rla w_1$,
which contradicts $\|g_1\|_*\ge 2$.
Therefore $w_1=1$.

Since $g_1$ is strongly non-split and $g_1^{w_2}$ is a cyclic conjugation,
$g_1^{w_2}$ is also strongly non-split and $| \supp(g_1^{w_2})| = | \supp(g_1)| \ge 2$.
Therefore $\|g_1^{w_2} \|_* \ge 2$ (by Lemma~\ref{lem:SNS}).

Since $w_1=1$, $\|g\|_*\le r$, $\|h\|_* \le r$,
$\|g_1\|_* \ge 2$ and $\|g_1^{w_2}\|_*\ge 2$,
using the inequalities in the proof of Claim 1,
we have
\begin{align*}
\|a\|_*
& \le  \frac{\|g\|_*- \|g_1\|_*}2+2
 \le  \frac{r-2}2+2=\frac r2+1,\\
\|w_3\|_*
& \le  \frac{\|h\|_*- \|g_1^{w_2}\|_*}2+2
\le  \frac{r-2}2+2=\frac r2+1,\\
\|w\|_* &\le \|a\|_* + \|w_1\|_*
+ \|w_2\|_* + \|w_3\|_*\\
&\le \left(\frac r2+1\right) + 0 + \|w_2\|_* + \left(\frac r2+1\right)
=\|w_2\|_*+r+2.
\end{align*}
Therefore
$\|a\|_*\le \frac 12 r+1$,
$\|w_3\|_*\le \frac 12 r+1$
and
$\|w_2\|_*\ge \|w\|_*-r-2\ge R-r-2$.
\end{proof}

\medskip\noindent
\textbf{Claim 3.}\ \
Let $\epsilon=1$ (resp.\ $\epsilon=-1$) if $g_1^{w_2}$ is a left (resp.\ right) cyclic conjugation.
Then there exists a quasi-root decomposition
$$
w=a(g_1^\epsilon)^{n}b
$$
such that $n\ge 2$ and $\|b\,\|_*\le\frac32r+2$.

\begin{proof}[Proof of Claim 3]
\noindent\smallskip
Suppose that $g_1^{w_2}$ is a left cyclic conjugation.
Then $w_2\le_L g_1^k$ for some $k\ge 1$ (by Proposition~\ref{prop:cc1}).
Hence
$$w_2=g_1^{n}d$$
is geodesic for some $0\le n\le k$ and $d\in A(\Gamma)$
with $\|d\|_*\le \|g_1\|_*+1$ (by Corollary~\ref{cor:po3}).
Notice that $n\ge 2$ because
$\|g_1 \|_* \le r$ whereas
$$
\|g_1^{n}\|_*\ge \|w_2\|_*-\|d\|_*
\ge (R-r-2)-(r+1)=R-2r-3 \ge r+4.
$$
Since $a^{-1}w=w_1w_2w_3$ and $w_1=1$, we have
\begin{align*}
w &= aw_2w_3
= a g_1^{n}d w_3.
\end{align*}

We will now prove that the decomposition $w=a g_1^{n}d w_3$ is geodesic.
Since $g=ag_1a^{-1}$ is geodesic and $g_1$ is cyclically reduced,
$ag_1^{n}$ is geodesic (by Lemmas~\ref{lem:cc} and~\ref{lem:cred1}).
Since both  $w_2w_3$ and $w_2=g_1^{n}d$ are geodesic,
$g_1^{n}dw_3$ is also geodesic.
Recall $\|g_1^n\|_*\ge r+4\ge 2$.
Therefore $w = a g_1^{n}dw_3$ is geodesic (by Lemma~\ref{lem:st5}).

Let $b=dw_3$. Then $w=ag_1^{n}b$ is geodesic and
\begin{align*}
\|b\|_*
& \le \|d\|_*+\|w_3\|_*
\le \Bigl(\|g_1\|_*+1\Bigr) + \left(\frac r2+1\right)
\le r+1+\frac r2+1=\frac 32r+2.
\end{align*}
Therefore $w=ag_1^n b$ is a quasi-root decomposition
with the desired properties.

\smallskip
Now suppose that $g_1^{w_2}$ is a right cyclic conjugation.
Then $\left(g_1^{-1} \right)^{w_2}$ is a left cyclic conjugation
(by Proposition~\ref{prop:cc1}).
From the above argument,
$w=a(g_1^{-1})^{n}b$ is a quasi-root decomposition
with the desired properties.
\end{proof}

The proof is now completed.
\end{myproof}

\begin{remark}\label{rmk:proA}
In Proposition~\ref{pro:A}, notice that
\begin{align*}
&w=a(g_1^\epsilon)^n b=(a(g_1^\epsilon)^na^{-1}) ab=(g^\epsilon)^n ab,\\
&\textstyle \|ab\|_* \le \|a\|_*+\|b\|_*\le (\frac12r+1)+(\frac32r+2)=2r+3.
\end{align*}
Thus one could understand Proposition~\ref{pro:A} as follows:
if\/ $\|w\|_*$ is large but both $\|g\|_*$ and $\|w^{-1}gw\|_*$ are small,
then $w=g^nc$  for some integer $n$
and $c\in A(\Gamma)$ with $\|c\|_*$ small.
\end{remark}

If $\|w^{-1}gw\|_*\le r$ and $g=ag_1a^{-1}$ in the statement of
Proposition~\ref{pro:A} are respectively replaced with
$\|wgw^{-1}\|_*\le r$ and $g=a^{-1}g_1a$,
then we have the following corollary.

\begin{corollary}\label{cor:Ac}
Let $g,w\in A(\Gamma)\setminus\{1\}$ and $r,R\ge 1$ such that
$$
\|g\|_*\le r,\quad
\|wgw^{-1}\|_*\le r,\quad
\|w\|_*\ge R,\quad
R\ge 3r+7.
$$
Then there exists a quasi-root decomposition
$$
w=b(g_1^\epsilon)^{n}a,
$$
where $a,b,g_1\in A(\Gamma)$,
$\epsilon\in\{\pm 1\}$ and $n\ge 2$
such that
\begin{itemize}
\item[(i)] $\|a\|_*\le\frac 12 r+1$ and $\|b\,\|_*\le\frac32r+2$;
\item[(ii)] $g_1$ is cyclically reduced and $g=a^{-1}g_1a$ is geodesic.
\end{itemize}
\end{corollary}

\medskip

We will now prove Theorem~\ref{thm:st}.

\medskip

\begin{myproof}[Proof of Theorem~\ref{thm:st}]
Choose $x,y\in A(\Gamma)$ with $d_*(x,y)\ge R$.
Let $w=yx^{-1}$, hence $\|w\|_*=d_*(x,y)\ge R$.

We may assume $\xi(1,w;r)\neq\{ 1\}$
because otherwise
$\xi(x,y;r)=x^{-1}\xi(1,w;r)x=\{1\}$
and there is nothing to prove.

By Lemma~\ref{lem:st0}(ii),
the set $\xi(1,w;r)$ is closed under taking a root,
i.e.~if $h^k\in \xi(1,w;r)$ for some $h\in A(\Gamma)$ and $k\ge 1$,
then $h\in \xi(1,w;r)$.
Therefore there exists a primitive element $g_0$
in $\xi(1,w;r)\setminus\{1\}$,
hence $\|g_0\|_* = d_*(g_0, 1)\le r$ and $\|wg_0w^{-1}\|_* = d_*(wg_0, w)\le r$.

We will now show that $g_0^{\pm 1}$ is uniquely determined from $w=yx^{-1}$.
Let
$$g_0=a^{-1}g_1a$$
be the geodesic decomposition such that $g_1$ is cyclically reduced.
Then $g_1$ is also primitive and $\|g_1\|_* \le \|g_0\|_* \le r$.
Since $R=(2V+7)r+8 \ge 3r+7$,
$(w,g_0,g_1,a,R,r)$ satisfies
the conditions on $(w,g,g_1,a,R,r)$ in Corollary~\ref{cor:Ac},
hence there exists a quasi-root decomposition
$$
w=b (g_1^\epsilon)^{n}a,
$$
where $b\in A(\Gamma)$, $\epsilon\in\{\pm 1\}$, $n \ge 2$,
$\|a\|_*\le\frac 12r+1$ and $\|b\|_*\le\frac 32r+2$.

Let $A= \frac 12r+1$ and $B=\frac 32r+2$.
Then $g_1^\epsilon$ is a primitive
$(B, A, r)^*$-quasi-root of $w$.
Observe that
$2A+2B+(2V+3)r+2
=(r+2)+(3r+4)+(2V+3)r+2
=(2V+7)r+8=R$,
hence
\begin{equation*} 
\|w\|_* \ge R=2A+2B+(2V+3)r+2.
\end{equation*}

The tuple $(w, g_1^{\epsilon}, b, a, B, A,r)$ now
satisfies the conditions on $(h, g_1, a_1, b_1, A, B,r)$
in Theorem~\ref{thm:qr}.
Therefore the primitive element $g_0^{\epsilon}$ is uniquely
determined from $w$
because $g_0^{\epsilon}=a^{-1}g_1^\epsilon a$ is the
rightward-extraction of the $(B, A, r)^*$-quasi-root $g_1^\epsilon$.
This means that each element of $\xi(1,w;r)$
is a power of $g_0$, hence
$\xi(1,w;r)\subset \langle g_0\rangle$.
Since
\begin{align*}
\|g_1^n\|_*
&=\|(g_1^\epsilon)^n\|_*\ge \|w\|_*- \|b\|_*-\|a\|_*
\ge R-B-A\\
&\ge A+B+(2V+3)r+2\ge 3,
\end{align*}
the cyclically reduced element $g_1$ is loxodromic (by Lemma~\ref{lem:loxo}),
hence $\|g_1^{(V-2)j} \|_* \ge j+2$ for all $j\ge 1$
(by Lemma~\ref{lem:loxo}, Proposition~\ref{pro:po4} and Corollary~\ref{cor:st4}).
Since $g_0$ is conjugate to $g_1$, $g_0$ is also loxodromic.
Since $g_0^{(V-2)j}=a^{-1}g_1^{(V-2)j}a$ is a geodesic decomposition
(by Lemma~\ref{lem:cc}(iii)),
$\|g_0^{(V-2)j} \|_* \ge \|g_1^{(V-2)j} \|_* \ge j+2$ for all $j\ge 1$.

If $k\ge (V-2)(r-1)$, then
$\|g_0^k\|_*\ge \|g_0^{(V-2)(r-1)}\|_*\ge (r-1)+2=r+1$,
hence $g_0^k\not\in \xi(1,w;r)$.
From this fact and Lemma~\ref{lem:st0},
it follows that
$$\xi(1,w;r)=\{1,g_0^{\pm 1},\ldots, g_0^{\pm k}\}$$
for some $1\le k\le (V-2)(r-1)-1$.

Let $g=x^{-1}g_0x$. Then $g$ is also loxodromic.
Since $\xi(x,y;r)
= x^{-1} \cdot \xi(1,w;r) \cdot x$,
$$
\xi(x,y;r) = \{1,g^{\pm 1},\ldots,g^{\pm k}\},
$$
hence (i) is proved.

Let $N(r)=2(V-2)(r-1)-1$.
Since $|\xi(x,y;r)|=2k+1\le 2(V-2)(r-1)-1=N(r)$,
the action of $A(\Gamma)$ on $(A(\Gamma),d_*)$
is $(R(r),N(r))$-acylindrical.

Since $g_0=a^{-1}g_1a$, $g=x^{-1}g_0x$ and $yx^{-1}=w=b(g_1^\epsilon)^n a$,
we get
$$
y = wx = b(g_1^\epsilon)^n ax
= bax\cdot x^{-1} a^{-1} (g_1^\epsilon)^n a x
= bax (g^\epsilon)^n.
$$
Hence
$d_*(y,x(g^\epsilon)^n)
=d_*(bax (g^\epsilon)^n,x(g^\epsilon)^n)
=\|ba\|_* \le \|b\|_* +\|a\|_* \le 2r+3$.
Therefore the Hausdorff distance between
the $\langle g\rangle$-orbits
$x\langle g\rangle$ and
$y\langle g\rangle$ is at most $2r+3$,
hence (ii) is proved.
\end{myproof}

\begin{remark}
The above proof shows
that $g_1^\epsilon$ is a primitive
$\left( \frac32 r+2, \frac12 r+1, r \right)^*$-quasi-root of
$w=b (g_1^\epsilon)^{n}a$.
Notice that the rightward-extraction of $g_1^\epsilon$ is $a^{-1}g_1^\epsilon a$,
and that $xgx^{-1} = a^{-1}g_1a$.
Therefore either $xg x^{-1}$ or $xg^{-1} x^{-1}$ is the rightward-extraction of
a primitive $(\frac32r+2,\frac12r+1,r)^*$-quasi-root of $yx^{-1}$.
\end{remark}

We are now ready to prove Theorem~\ref{thm:ex}.
\begin{myproof}[Proof of Theorem~\ref{thm:ex}]
By Theorem~\ref{thm:st}, the action of\/ $A(\Gamma)$ on $(A(\Gamma), d_*)$ is
$(R_1(r),N_1(r))$-acylindrical with
\begin{align*}
R_1(r) &=(2V+7)r+8,\\
N_1(r) &=2(V-2)(r-1)-1.
\end{align*}
Applying Lemma~\ref{lem:2act} to the above,
the action of\/ $A(\Gamma)$ on $(\Gamma^e, d)$ is
$(R(r), N(r))$-acylindrical with
\begin{align*}
R(r) &=D\cdot R_1(r+1)+2D=D (2V+7)(r+1)+ 10D,\\
N(r) &=N_1(r+1)=2(V-2)r-1.
\end{align*}

Choose $x,y\in V(\Gamma^e)$ with
$d(x,y)\ge R(r)$  and $\xi_{(\Gamma^e, d)}(x,y;r)\ne \{1\}$.
Then there exist $v_1, v_2\in V(\Gamma)$ and $w_1, w_2\in A(\Gamma)$ such that
$x = v_1^{w_1}$ and $y = v_2^{w_2}$.
By Lemma~\ref{lem:2act},
\begin{eqnarray*}
& d_* (w_1, w_2) \ge R_1(r+1),\\
& \xi_{(\Gamma^e,d)}(v_1^{w_1},v_2^{w_2}; r)
\subset \xi_{(A(\Gamma),d_*)}(w_1,w_2;r+1).
\end{eqnarray*}

Since $\xi_{(\Gamma^e,d)}(v_1^{w_1},v_2^{w_2}; r)\neq\{ 1\}$,
we have $\xi_{(A(\Gamma),d_*)}(w_1,w_2;r+1)\neq\{ 1\}$.
Hence $(w_1, w_2)$ satisfies all the conditions on $(x, y)$ in Theorem~\ref{thm:st}.
Therefore, by Theorem~\ref{thm:st}(i),
$$
\xi_{(\Gamma^e,d)}(x,y;r)\subset
\xi_{(A(\Gamma),d_*)}(w_1,w_2;r+1)
=\{1,g^{\pm 1}, g^{\pm 2},\ldots, g^{\pm k}\}
$$
for some loxodromic element $g\in A(\Gamma)$ and
$1\le k\le (V-2)r-1$,
hence (i) is proved.

Since the Hausdorff distance between the $\langle g\rangle$-orbits
of $w_1$ and $w_2$ is at most $2(r+1)+3=2r+5$ (by Theorem~\ref{thm:st}(ii)),
$w_2 = cw_1 g^n$ for some $n\in\Z$ and $c\in A(\Gamma)$ with
$\|c\|_*\le 2r+5$.
Hence we get (by Lemma~\ref{lem:KK})
\begin{align*}
d(x^{g^{n}}, y) & = d(v_1^{w_1 g^{n}}, v_2^{w_2})
 = d(v_1, v_2^{w_2 g^{-n} w_1^{-1}}) = d(v_1, v_2^{c})  \\
& \le d(v_1, v_2) + d(v_2, v_2^c) \le D + D(\|c\|_* +1) \\
& = D(\|c\|_* +2)
\le D(2r+7).
\end{align*}
Therefore the Hausdorff distance between
$x^{\langle g\rangle}$ and
$y^{\langle g\rangle}$ is at most $D(2r+7)$,
hence (ii) is proved.
\end{myproof}


\end{document}